\documentclass{article}

\usepackage{lmodern}
\usepackage[margin=1.5in]{geometry}
\usepackage[utf8]{inputenc} 
\usepackage[T1]{fontenc}    

\usepackage{tkz-tab}
\usepackage{float} 
\usepackage[utf8]{inputenc} 
\usepackage[T1]{fontenc}    
\usepackage{hyperref}       
\hypersetup{
	colorlinks=true,
	urlcolor=blue,
	citecolor=blue,
	linktoc=all,
	linkcolor=blue}
\usepackage{url}            
\usepackage{booktabs}       
\usepackage{amsfonts}       
\usepackage{nicefrac}       
\usepackage{microtype}      
\usepackage{amsmath}
\usepackage{lipsum}
\usepackage{graphicx}
\usepackage{xcolor}
\usepackage{cancel}
\usepackage{mathtools} 
\usepackage{booktabs} 
\usepackage{multirow}
\usepackage{diagbox}

\usepackage{array}

\usepackage{adjustbox}

\usepackage{caption}
\usepackage{subcaption}

\usepackage{tkz-tab}

\usepackage[title]{appendix}

\graphicspath{ {./images/} }

\usepackage{natbib}
\newcommand{\thickbar}[1]{\mathbf{\bar{\text{$#1$}}}}

\newcommand{\tod}{\stackrel{d}{\longrightarrow}}
\newcommand{\toP}{\stackrel{P}{\longrightarrow}}

\newtheorem{Lem}{Lemma}

\counterwithin*{Lem}{section}

\newtheorem{LemApp}{Lemma}

\counterwithin*{LemApp}{section}

\counterwithin*{Theo}{section}

\newtheorem{Prop}{Proposition}

\counterwithin*{Prop}{section}

\counterwithin*{Def}{section}

\newtheorem{Cor}{Corollary}

\counterwithin*{Cor}{section}

\definecolor{forestgreen}{rgb}{0.13, 0.55, 0.13}
\definecolor{darkblue}{rgb}{0.0, 0.0, 0.65}
\definecolor{violet2}{rgb}{0.5,0,0.5}
\definecolor{orange2}{rgb}{0.8, 0.1, 0.1}       
\definecolor{red2}{rgb}{1, 0.4, 0}

\newcommand{\gr}[1]{{\textcolor{black}{#1}}}

\title{On the use of a local $\hat{R}$ \\to improve MCMC convergence diagnostic}

\author{Théo Moins \thanks{Univ. Grenoble Alpes, Inria, CNRS, Grenoble INP, LJK, 38000 Grenoble, France.}
\and Julyan Arbel \footnotemark[1]
\and Anne Dutfoy \thanks{EDF R\&D dept. Périclès, 91120 Palaiseau, France.}
\and Stéphane Girard \footnotemark[1]}

\begin{document}
\maketitle

\begin{abstract}
Diagnosing convergence of Markov chain Monte Carlo is crucial and remains an essentially unsolved problem.
Among the most popular methods, the potential scale reduction factor, commonly named $\hat{R}$, is an indicator that monitors the convergence of output chains to a target distribution, based on a comparison of the between- and within-variances.
Several improvements have been suggested since its introduction in the 90s. 
Here, we aim at better understanding the $\hat{R}$ behavior by proposing a localized version that focuses on quantiles of the target distribution. 
This new version relies on key theoretical properties of the associated population value.
It naturally leads to proposing a new indicator $\hat{R}_\infty$, which is shown to allow both for localizing the  Markov chain Monte Carlo convergence in different quantiles of the target distribution, and at the same time for handling some convergence issues not detected by other $\hat{R}$ versions.
\end{abstract}


\section{Introduction}

Markov chain Monte Carlo (MCMC) algorithms have strongly contributed to the popularity of Bayesian models to sample from posterior distributions, especially in high-dimensional or high computational settings.
This success results in a variety of softwares increasingly used for a wide range of applications: Stan \citep{carpenter2017stan}, PyMC3 \citep{PyMC3}, NIMBLE \citep{nimble}, or Pyro \citep{bingham2019pyro}, to cite a few.
The fundamental idea behind these algorithms is the convergence of the sampling distribution to the target (typically the posterior) when the number of samples goes to infinity. 
A major challenge is therefore to know if the behaviour for a finite number of draws is satisfactory or not.
This allows for a handle on the number of iterations to be drawn, which is all the more crucial in complex models with costly sampling schemes.
See \citet{roy2020convergence} for a recent literature review on convergence diagnostics.

\subsection{Diagnosing MCMC convergence}

    \begin{figure}
        \hspace{-0.5cm}
        \setlength\tabcolsep{2pt}
        \begin{tabular}{ccc}
        \includegraphics[trim=0cm 1.5cm 0cm 0cm, clip, width = .32\textwidth]{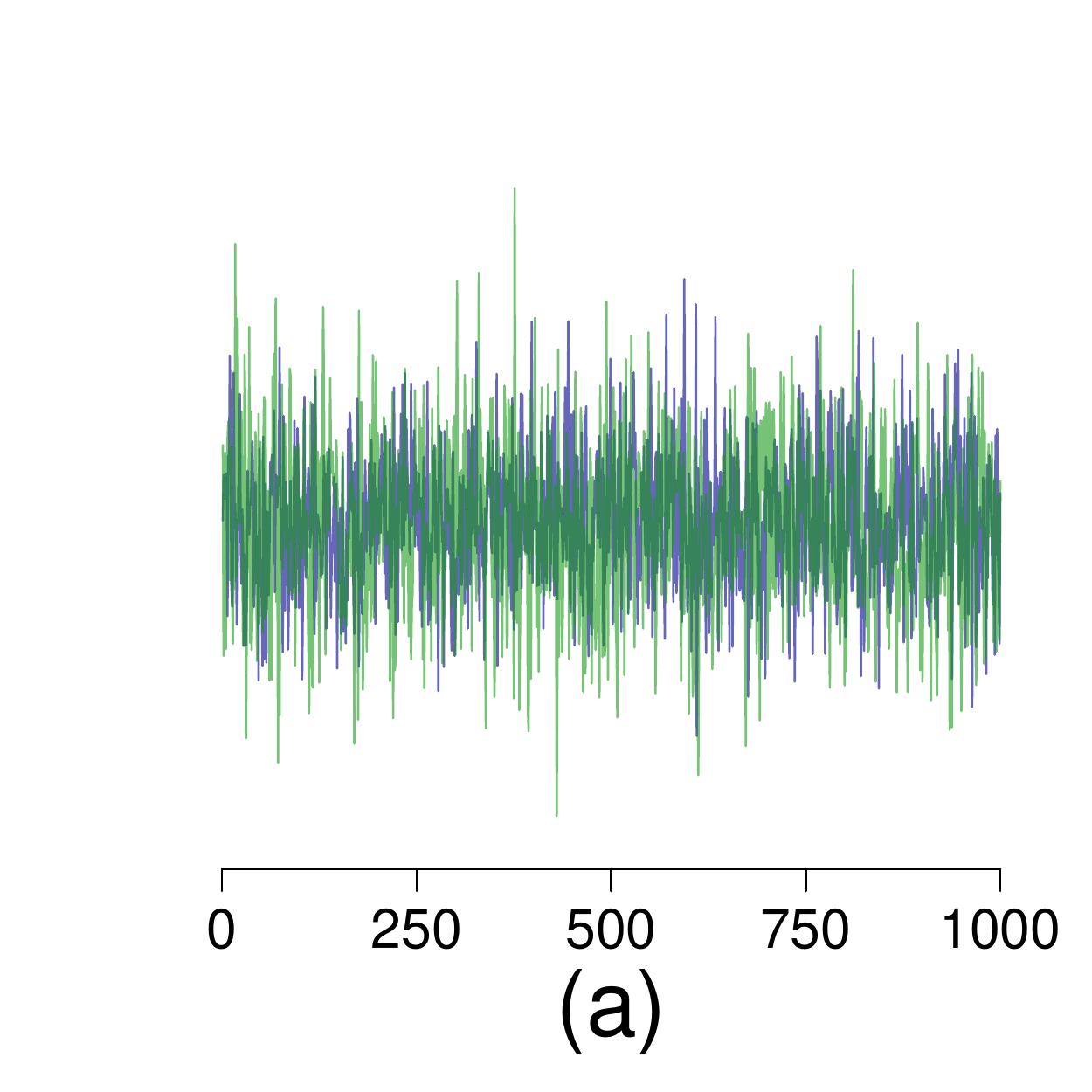} &
        \includegraphics[trim=0cm 1.5cm 0cm 0cm, clip, width = .32\textwidth]{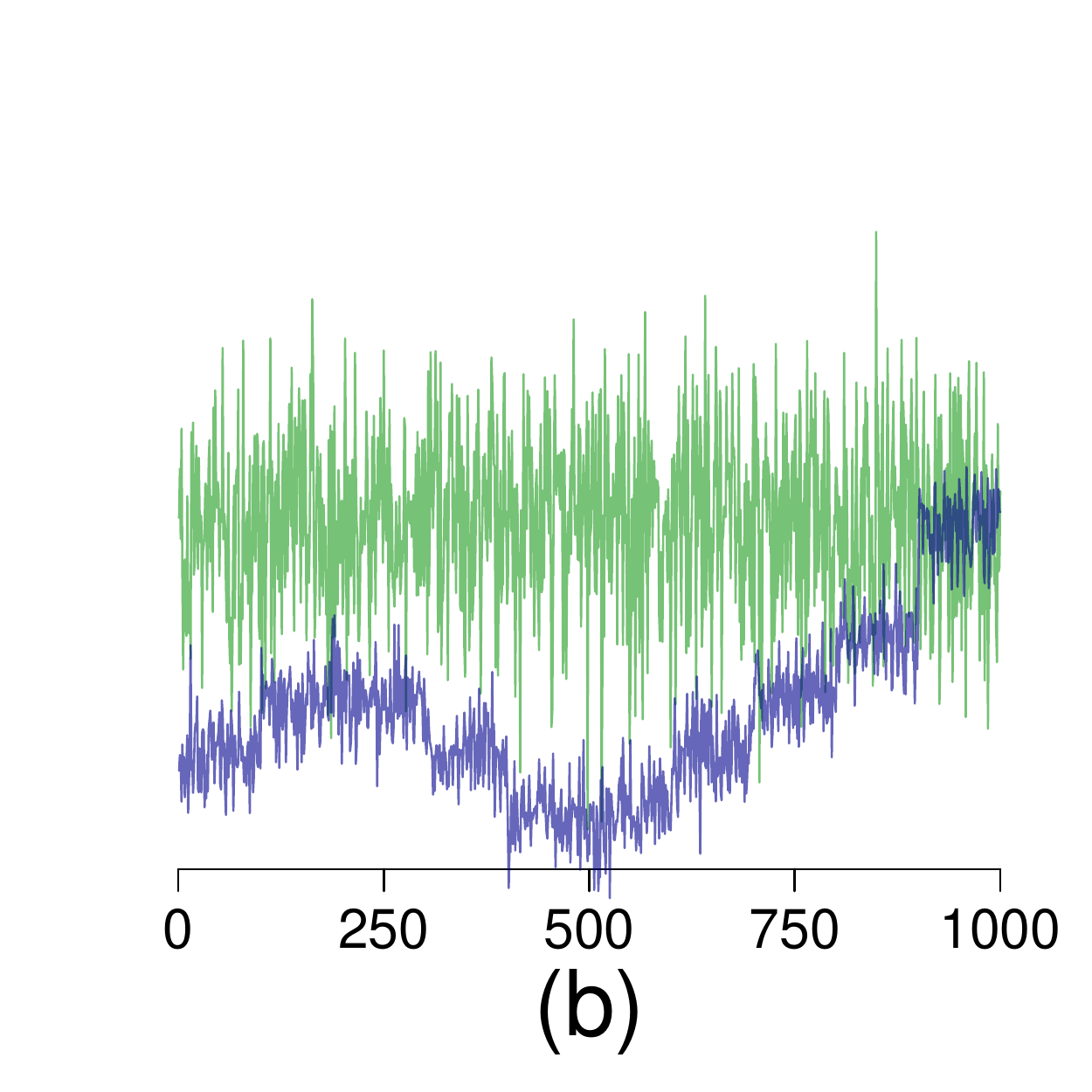} &
        \includegraphics[trim=0cm 1.5cm 0cm 0cm, clip, width = .32\textwidth]{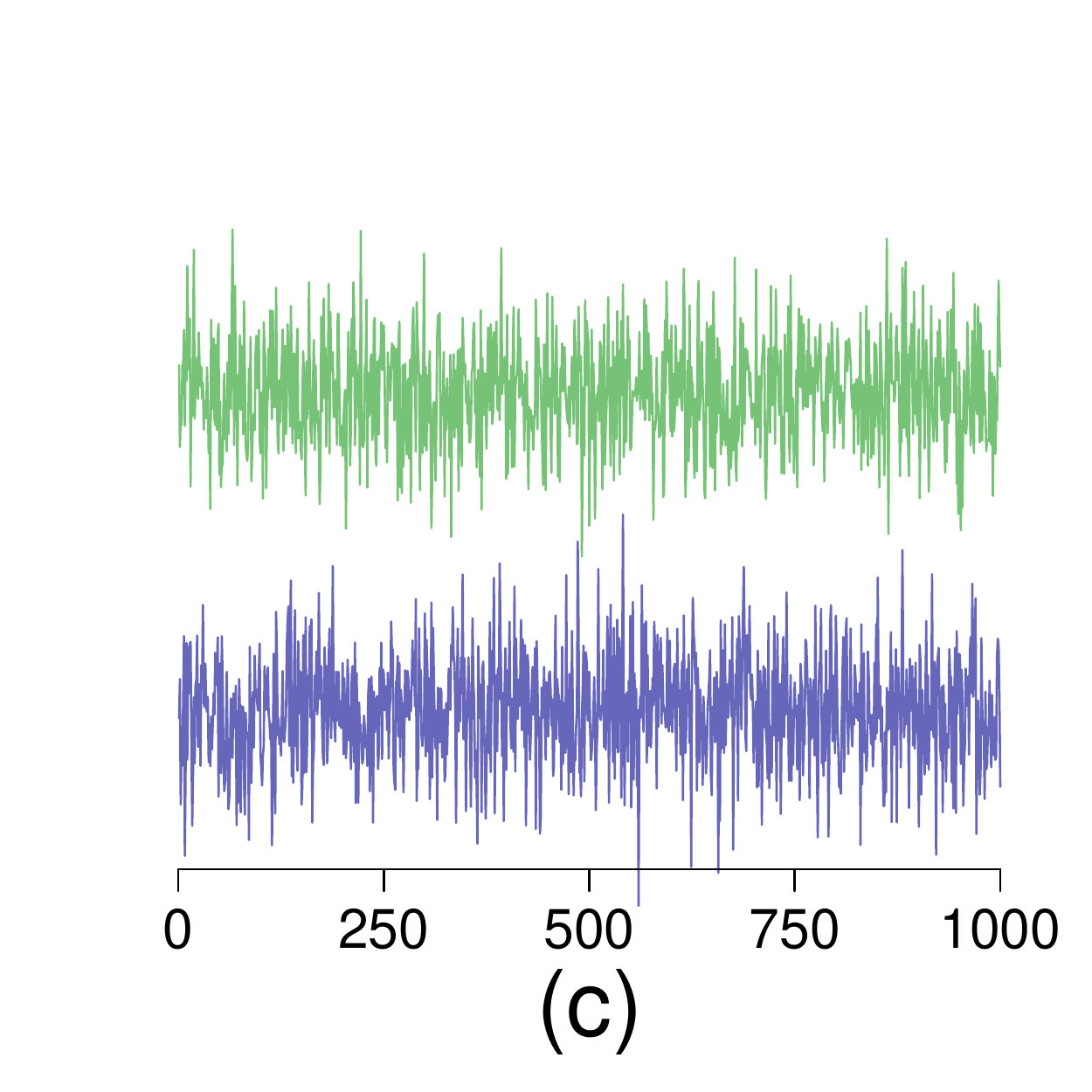}\\
        \small{iterations} & \small{iterations} & \small{iterations}
        \end{tabular}
        \caption{Traceplots illustrating convergence and two types of non-convergence of MCMC.
        \gr{Left:} nothing indicates a convergence issue, as the two chains seem to have the same stationary distribution. \gr{Middle:} the blue chain is still in an exploration phase and therefore is not stationary. 
        \gr{Right:} example where having multiple chains helps detecting a mixing issue despite a stationarity appearance of each.}
        \label{fig:mixing_chains}
    \end{figure}

Two frequently used properties to verify chains convergence are stationarity and mixing \citep[see][for a discussion]{Vats2021discussion}. 
\textit{Stationarity} is related to the invariance property of the target distribution $F$ for standard MCMC algorithms like Metropolis--Hastings or Gibbs sampling \citep{bookRobert}: if $\theta^{(i)}$ is the $i$th element of an MCMC chain, then $\theta^{(i)} \sim F$ implies $\theta^{(i+1)} \sim F$, so that as soon as an element of the chain is distributed according to $F$, all the following ones will be too.
Thus, a chain whose distribution changes drastically during iterations is still in the exploration phase and is therefore not stationary (see middle panel in Figure~\ref{fig:mixing_chains}).
\textit{Mixing} refers in practice to the exploration of the support of $F$: slow mixing chains correspond to chains that only explore a subset of the parameter space, which can lead to strong bias in the distribution \citep[see][for a more rigorous definition]{robert1995convergence}.
A common way to limitate mixing issues is to run several chains in parallel with different starting points, which also allows comparing the chains together.
Stationarity and mixing are two properties that can be treated separately: in principle, being stationary implies convergence to the target distribution and thus necessarily also mixing, but in practice there are examples of chains that seem to have reached stationarity but are not mixing (see right panel in Figure~\ref{fig:mixing_chains}), hence the need for comparing multiple chains.

We place ourselves in the case of several chains: consider $m$ chains of size $n$, with $\theta^{(i,j)}$ denoting the $i$th draw from chain $j$.
We focus here on the Gelman--Rubin diagnostic \citep{gelman1992inference}, named potential reduction scale factor and commonly denoted by $\hat{R}$. It is by far one of the most popular methods to assess MCMC convergence, used in particular in Stan, PyMC3, or NIMBLE.
The original heuristic for $\hat{R}$ construction is the comparison between two estimators that converge to the target variance $\text{Var}[\theta]$, based on $\hat{W}$ and $\hat{B}$, respectively the estimated within- and between-variances.
This diagnostic has the advantage of being scalar even in the case of a huge number of chains and comes with a rule of thumb that makes it very easy to use: generally $\hat{R} \geq 1$, and if it is greater than a given threshold (for example $1.01$), then a convergence issue is raised.
This was originally constructed to diagnose mixing issues only, but \citet[Section~11.4]{gelman2013bayesian} suggest splitting the chains in two before computing $\hat{R}$ to check for stationarity at the same time.
We will also always consider this split version of $\hat{R}$ throughout this paper, thus focusing only on the problem of mixing diagnostic.

\subsection{Different $\hat{R}$ versions and their limitations}
\label{subsec:intro_comparison_rhat}

The original $\hat{R}$ of \citet{gelman1992inference} has some limitations that are listed here with associated improvements suggested in the literature.


\newcounter{limitations}
\newenvironment{limitations}[1][]{\refstepcounter{limitations}\par
   \noindent\textbf{\sffamily L\thelimitations.} #1}

\begin{limitations}{\textbf{\sffamily It must be compared to an arbitrary chosen threshold.}}\label{L1}
To use $\hat{R}$, a threshold must be set to determine a convergence issue.
Originally set to $1.1$, \citet{vats2018revisiting} note that this choice is arbitrary and usually too optimistic.
Thus, the authors propose a threshold according to a confidence level based on a relationship made with effective sample size (ESS).
This observation was then shared by \citet{Vehtari} who suggest dropping the threshold to $1.01$.
Driven by \gr{practical arguments}, this choice remains unprincipled nor theoretically justified, which is related to the next limitation.
\end{limitations}

\begin{limitations}{\textbf{\sffamily  It suffers from a lack of interpretability.}}\label{L2}
How to interpret a given value of $\hat{R}$? 
By construction, $\hat{R}$ is a ratio of two quantities that must estimate the posterior variance.
Therefore, having a value close to one can be seen as having two correct estimations of the same quantity, which is an indication of convergence. 
However to our knowledge, no study investigates the theoretical or population value $R$ that $\hat{R}$ aims at estimating, which would shed light on what is actually diagnosed. 
Typically, chains such that $\hat{R} \approx 1$ do not necessarily correspond to mixing chains: \citet{Vehtari} exhibit some counter-examples in order to motivate a more robust version called rank-$\hat{R}$.
Still, the different versions of $\hat{R}$ only allow to draw conclusions when they are significantly greater than 1, and the common properties of chains producing $\hat{R} \approx 1$ are not well known as they are constructed at the estimator level.
\end{limitations}

\begin{limitations}{\textbf{\sffamily   It is not robust to certain types of non-convergence.}}\label{L3}
Traditional $\hat{R}$ can be fooled, in the sense that $\hat{R} \approx 1$ without convergence. This motivates the construction of rank-$\hat{R}$ \citep{Vehtari}, based on two cases where the original $\hat{R}$ is not robust:
\begin{enumerate}
    \item[(i)] When the mean of the target distribution is infinite: in that case $\hat{W}$ and $\hat{B}$ are ill-defined and $\hat{R} \approx 1$ even though the chains follow different distributions.
    One solution is to apply rank transformation on the chains before computing $\hat{R}$ (this version is named bulk-$\hat{R}$ by \citealp{Vehtari}).
    \item[(ii)] When the means of the chains are equal: in that case, the variance of means $\hat{B}$ is zero, and so $\hat{R} \approx 1$ even if the variances of chains are different.
    Here in addition to the rank-transformation, transforming the chains to get the deviation from their median allows to overcome this problem (this version is named tail-$\hat{R}$ by \citealp{Vehtari}).
\end{enumerate}
Defining rank-$\hat{R} = \max \{\text{bulk-}\hat{R}, \text{tail-}\hat{R}\}$ overcomes the two issues at the same time. 
However, this robustness can be seen as very specific and can easily be fooled by simple examples. 
One way is to consider chains with different distributions, but with (i) same mean (to fool bulk-$\hat{R}$), and (ii) same mean over the median (to fool tail-$\hat{R}$).
For example, uniform $\mathcal{U}(\mu-2\sigma, \mu+2\sigma)$, normal $\mathcal{N}(\mu, \frac{\pi}{2}\sigma^2)$, or Laplace $\mathcal{L}(\mu, \sigma)$ distributions share the same mean (equal to $\mu$) and same mean over the median (equal to $\sigma$), and thus mixing them yields rank-$\hat{R} \approx 1$.
We provide an example and a more general framework to construct such \gr{cases} in Appendix~\ref{app:counter_example} of the supplementary material. 
One illustration can also be found in the right column of Figure~\ref{fig:example_known_dist}.
Although these counter-examples may never appear in practice, they do show some fairly counter-intuitive results \gr{that the additional layer of computation carried by rank-$\hat{R}$ makes even more difficult to analyse}.
\end{limitations}

\begin{limitations}{\textbf{\sffamily It does not target a specific quantity of interest.}}\label{L4}
\label{paragraph:L4}
Another point raised by \citet{Vehtari} is that the convergence diagnostic does not depend on inferential features of interest. 
It might be more precise to speak of convergence for a given posterior quantity, typically a mean, higher order moment, or quantile. 
\gr{Typically, practitioners apply $\hat R$ on quantities of interest such as the log-likelihood, the posterior density, or quantiles. 
On their side, \citet{Vehtari} suggest a local transformation on ESS to obtain a tail-ESS associated with $5\%$ and $95\%$ quantiles.}
\end{limitations}

\begin{limitations}{\textbf{\sffamily   It is associated with a univariate parameter.}}\label{L5}
Although the vast majority of Bayesian models have multivariate parameters, $\hat{R}$ focuses on univariate convergence (i.e. convergence of margins). 
Some multivariate extensions exist, like \citet{brooks1998general} or \citet{vats2019multi}, but do not seem to be universally accepted: for example Stan or PyMC3 use instead a table containing univariate $\hat{R}$ with one value per parameter. 
However, assessing convergence on margins misses the point of dependence among parameter components, and does not guarantee the convergence of the joint distribution.
Another version of $\hat{R}$ called $R^*$ is suggested by \citet{lambert2021r} and can deal with multivariate parameters: the idea is to use a classification algorithm which, in the case of converging chains, would not be able to identify to which chain a sample belongs.
To avoid a result depending on the seed of the experiment, the authors suggest to draw several samples from the simplex obtained with the classification algorithm. 
In addition to the interpretability issues mentioned previously, this method has the constraint of not being able to study only a scalar value but a histogram, to check to what extent it contains or not the value $1$. 
\end{limitations}
    
We take a step forward in addressing \gr{all} these limitations with a localized version of $\hat{R}$ briefly introduced in \cite{moins2021discussion} and developed here: we analyze $\hat{R}(x)$, a local version of $\hat{R}$ associated with a given quantile $x$, and the corresponding population value $R(x)$.
This study leads us to propose a new indicator $\hat{R}_\infty$. In addition to being more interpretable, $\hat{R}_\infty$ shows better results than $\hat{R}$ in terms of MCMC convergence diagnostic, both on simulated experiments and on Bayesian models.
As with all other versions of $\hat{R}$, this one can be applied to any MCMC algorithm: Metropolis--Hastings, HMC \citep{neal2011mcmc}, NUTS \citep{hoffman2014no}, etc.

The rest of the paper is organized as follows: we introduce in Section~\ref{sec:r-hat-def} the population version $R(x)$ and the corresponding sample version $\hat{R}(x)$, as well as their scalar counterparts $R_\infty$ and $\hat{R}_\infty$. \gr{Since this proposed version depends on a quantile $x$ and is constructed at a population level, it is both targeting a specific quantity of interest and interpretable, addressing respectively limitations L\ref{L4} and L\ref{L2}. We also establish several properties on the behaviour of $R(x)$ function and on the convergence of the estimator $\hat{R}(x)$, helping in establishing a threshold and addressing limitation L\ref{L1}.}
Our proposed approach to deal with the multivariate case of limitation L\ref{L5} is described in Section~\ref{sec:multivariate}. Some empirical results are given in Section~\ref{sec:simulations}, \gr{showing that our proposed solution helps overcoming many of convergence issues identified in limitation L\ref{L3}}. We conclude in Section~\ref{sec:conclusion}.
All proofs and details of the calculations are provided in the supplementary material, and experiments are available online\footnote{\label{note1}\url{https://theomoins.github.io/localrhat/Simulations.html}}  as well as the R package \texttt{localrhat} \citep{moinsLocalrhat} containing our diagnostic implementation.

\section{Local version of $\hat{R}$}
\label{sec:r-hat-def}

    Since the original version of \cite{gelman1992inference}, the heuristic for the construction of $\hat{R}$ \gr{was based on an analysis of variance}.
    It consists in comparing two estimators of the posterior variance Var$[\theta]$. 
    The first one is the within-variance $\hat{W}$, which underestimates Var$[\theta]$ as the bias of the estimator is (most of the time) strictly negative if the elements of the chains are not i.i.d, see \cite{vats2018revisiting}.
    The second one adds the between-variance $\hat{B}$ as a bias correction. This typically overestimates Var$[\theta]$ if the initial values are chosen over-dispersed. 
    As pointed out by \cite{vats2018revisiting}, following this heuristic does not exclude the use of other estimators of the bias than $\hat{B}$.
    Moreover, defining $\hat{R}$ at the sample level hinders a theoretical study of a population version to be conducted.
    Another justification can start with the law of total variance: 
    assume that a univariate $\theta$ is sampled using $m$ chains, and let $Z\in\{1,\dots,m\}$ be the corresponding index of the chain. Then,
    \begin{equation}
        \label{eq:LTV}
        \text{Var}[\theta] = \mathbb{E}_Z[\text{Var}_{\theta\mid Z}[\theta \mid Z]] + 
        \text{Var}_Z[\mathbb{E}_{\theta\mid Z}[\theta \mid Z]].
    \end{equation}
    The two terms in the right-hand side correspond respectively to the population versions of the within-variance $W$ and the between-variance $B$. 
    Replacing them by their estimated versions yields the original $\hat{R}$ formula of \cite{gelman1992inference}.
    In the following, we use (\ref{eq:LTV}) on a chains transformation which allows to localise convergence at a given quantile.
    For the theoretical study, we suppose stationarity of the chains to focus only on chain mixing issues.
    Thus, samples within a chain $j \in \{1, \ldots, m\}$ may be correlated but are all distributed according to the same distribution $F_j$ which may vary with $j$.\\
    
    \subsection{Population version}
    \label{subsec:theorical_def}

For all $x\in{\mathbb R}$, introduce the Bernoulli random variable $I_x = {\mathbb I}\{\theta \leq x\}$, where ${\mathbb I}\{\cdot\}$ denotes the indicator function. 
Similarly to the Raftery--Lewis diagnostic \citep{raftery1991many}, the idea of our local convergence estimate is decidedly simple: we use $I_x$ in place of $\theta$ in the original Gelman--Rubin construction.
The population within-chain and between-chain variances at point $x$ are then defined respectively as
$W(x) = \mathbb{E}\left[\text{Var}[I_x \mid Z ]\right]$ and $B(x) = \text{Var}[\mathbb{E}[I_x \mid Z ]]$.
Note that both quantities exist whatever the tail heaviness of $\theta$ distribution thanks to introduction of the indicator function, thus relaxing moment conditions of the original $\hat{R}$.
We define the associated population $R(x)$ as
$$
 R(x) = \sqrt{\frac{W(x)+B(x)}{W(x)}}.
$$
It turns out that under the assumption of stationarity for each chain, $R(x)$ can be expressed in closed-form with respect to the chains' distribution.
\begin{Prop}
\label{prop-calcul-R}
    Suppose that, for any $j\in\{1,\dots,m\}$, $\mathbb{P}(Z = j) = 1/m$ and $\theta$ given $Z = j$ has cumulative distribution function (cdf) $F_j$. Then, one has for any $x\in\mathbb{R}$:
    \begin{equation}
        \label{eq:R_theorique}
        R(x) = \sqrt{1 + \frac{\sum_{j=1}^m\sum_{k=j+1}^m \left(F_j(x)-F_k(x)\right)^2}{m\sum_{j=1}^m F_j(x)(1-F_j(x))}}.
    \end{equation}
\end{Prop}
Thus, using $I_x$ instead of $\theta$ defines a local convergence estimate at any point $x$ which quantifies a distance between the $F_j$'s.
This allows for \gr{diagnosing} convergence relatively to a quantile one wants to estimate (for a posterior credible interval for example).
The following proposition states straightforward properties of $R(x)$ emanating from~(\ref{eq:R_theorique}):
\begin{Prop} \gr{The population $R(x)$ satisfies the following properties:} \hfill
\label{prop-first}
\begin{itemize}
    \item [\normalfont (i)] $R(x) \geq 1$ for all $x\in \mathbb{R}$.
    \item [\normalfont (ii)]  $R(x) =1$ for all $x\in \mathbb{R}$ if and only if $F_1=\dots=F_m$.
    \item [\normalfont (iii)] $R(x) \to 1$ as $|x| \to \infty$.
    \item [\normalfont (iv)] $R(x)$ inherits continuity property of $F_1, \ldots, F_m$ if the support of the $F_j$'s are overlapping.
\end{itemize}
\end{Prop}
Based on these results and in order to summarize this continuous index into a scalar one, we may also consider its supremum over $\mathbb{R}$:
\begin{equation}
    R_\infty = \sup_{x \in\mathbb{R}} R(x).
\end{equation}
Note that, in view of Proposition~\ref{prop-first}(iv), $R_\infty$ is finite simply as soon as the $F_j$'s are continuous with overlapping supports. 
Considering $R_\infty$ amounts to considering the local version $R(x)$ corresponding to the quantile $x$ with the poorest convergence when no information is given on the posterior interval used for inference.

\subsection{Sample version}

    Population version $R(x)$ can be estimated by replacing the $F_j(x)$'s in~(\ref{eq:R_theorique}) by their empirical counterparts $\hat{F}_j(x) = \frac{1}{n} \sum_{i=1}^{n} {\mathbb I}\{\theta^{(i,j)} \leq x\}$.
    This is equivalent to computing the original version of $\hat{R}$ on indicator variables $I_x^{(i,j)} = \mathbb{I}\{\theta^{(i,j)} \leq x\}$ instead of $\theta^{(i,j)}$.
    This connects with the Raftery--Lewis diagnostic \citep{raftery1991many} and more recently with \cite{Vehtari} who suggest this transformation for effective sample size (ESS) to construct graphical diagnostics or ``tail-versions'' of this diagnostic. 
    Moreover, a rank-normalization step is added in \cite{Vehtari}'s to prevent from infinite moments, although using $I_x^{(i,j)}$ ensures the index existence whatever the $\theta^{(i,j)}$ distribution is.
    Skipping this step for $\hat{R}$ yields an explicit expression of what is estimated in the stationary case with~\eqref{eq:R_theorique}. 
    This makes the diagnostic more interpretable and allows us to obtain key theoretical results for the associated theoretical $R$ and $R_\infty$.
    
    Note that for a given number of chains $m$ and chain length $n$, $\hat{R}(x)$ can only take $m(n+1)$ different values, as the computation is based on $n m$ indicator variables.
    Thus, the best accuracy we can obtain for $\hat{R}_\infty$ for a given $n$ and $m$ consists in evaluating $\hat{R}(x)$ at all the $\theta^{(i,j)}$'s.
    This can be accelerated by subsampling, often with limited decrease in accuracy.
    
\subsection{Convergence properties} 
    \label{subsec:convergence-rhat}
    
    Let us assume that all $m$ chains are mutually independent and have converged to a common distribution so that $F_1=\dots=F_m\eqqcolon F$.
    Assume, moreover, that a Markov chain central limit theorem holds \citep[see for instance][Theorem~6.65]{bookRobert}, so that we can write
    \begin{equation}
        \sqrt{n}(\hat F_j(x)-F(x)) \tod {\mathcal N}\left(0, \sigma^2(x) \right),
        \label{eq:CLT_MC_j}
    \end{equation}
    as $n\to\infty$, for all $j\in\{1,\dots,m\}$ and where $\sigma^2(x)$ is some asymptotic variance.   In particular in the i.i.d. setting, $\sigma^2(x) = F(x) (1-F(x))$.
    Letting $\hat{F}(x) = \frac{1}{nm} \sum_{j=1}^m\sum_{i=1}^{n} {\mathbb I}\{\theta^{(i,j)} \leq x\} = \frac{1}{m} \sum_{j=1}^m \hat F_j(x)$
    and taking into account of the independence between chains yield
    \begin{equation}
        \sqrt{nm}(\hat F(x)-F(x)) \tod {\mathcal N}\left(0, \sigma^2(x) \right),
        \label{eq:CLT_MC}
    \end{equation}
    as $n\to\infty$, and $\sigma(x) / \sqrt{nm}$ can be interpreted as the Monte Carlo standard error (MCSE) associated with $\hat F(x)$. 
    Following the definition of the ESS used in \cite{Gong2016} or \cite{vats2019multi}, we can define a local-$\text{ESS}$ as the ratio of the target variance to the squared MCSE:
    \begin{equation}
    \label{eq-def-local-ESS}
        \text{ESS}(x) = nm \frac{F(x) (1- F(x))}{\sigma^2(x)}.
    \end{equation}
    This quantity is in line with the definition of ESS for quantile of \cite{Vehtari}, and has already been studied by \cite{raftery1991many} who focus on this indicator transformation and approximate the resulting process as a two-state Markov chain.
    This yields an explicit expression of the stationary distribution $F$, which 
    can be used to obtain an expression of ESS$(x)$ as a function of the transition probabilities. 
    Several limitations of this two-state Markov chain approximation are raised by  \cite{brooks1999, doss2014markov}, for example.
    A more general way to estimate ESS$(x)$ is to apply the same idea as in the definition of the local $\hat{R}(x)$: use any estimator of ESS \citep{bookRobert, gelman2013bayesian} on indicator variables $I_x^{(i,j)}$ instead of $\theta^{(i,j)}$.
    
    Combining the asymptotic result \eqref{eq:CLT_MC} with expression \eqref{eq-def-local-ESS} yields the following large $n$ limiting distribution result on $\hat{R}(x)$ ($\chi^2_{m-1}$ denotes the chi-square distribution with $m-1$ degrees of freedom).
    \begin{Prop} 
        \label{prop:asymptotic_R-hat}
    Assume that all $m$ chains are mutually independent and have converged to a common distribution $F\coloneqq F_1=\dots=F_m$. Then:
    \begin{itemize}
        \item [\normalfont (i)] The distribution of $\hat{R}_\infty$ does not depend on the underlying distribution $F$.
        \item [\normalfont (ii)] For any $x\in{\mathbb R}$,
        $\mathrm{ESS}(x)(\hat{R}^2(x) - 1) \tod \chi^2_{m-1}$ as $n\to\infty$.
    \end{itemize}
    \end{Prop}
    \gr{Note that casting the problem of convergence monitoring in terms of analysing components of variance from multiple sequences dates back to \cite{gelman1992inference}, Section 2.2, and earlier works by \cite{fosdick1959calculation,gelfand1990sampling}. 
    Let us highlight that the assumption $F_1(x)=\cdots=F_m(x)$ is equivalent to the ANOVA hypothesis $\mathbb{E}(I_x^{(\cdot,1)})=\cdots=\mathbb{E}(I_x^{(\cdot,m)})$ and that the statistics studied in Proposition~\ref{prop:asymptotic_R-hat}(ii)
    can similarly be rewritten in terms of the ANOVA test statistics: $\hat{R}^2(x) - 1=\hat B(x)/\hat W(x)$, where
    $\hat B(x)$ and $\hat W(x)$ are the respective empirical counterparts of $B(x)$ and $W(x)$. These interpretations can then be
    used to derive a statistical test on the convergence of the chains. To this end, note also} that the limit in distribution of Proposition~\ref{prop:asymptotic_R-hat}(ii) still holds when $\text{ESS}(x)$ is replaced by a consistent estimator $\widehat{\text{ESS}}(x)$.
    This result allows computing the type I error associated with the null hypothesis that $\hat{R}(x) = 1$, in other terms that all the chains have converged to a common distribution at $x$.
    Let $z_{m-1, 1-\alpha}$ be the quantile of level $1-\alpha$ of the $\chi^2_{m-1}$ distribution,
    and introduce the associated threshold
    \begin{equation}
        \label{eq-Rlimite}
    {R}_{\text{lim},\alpha}(x) \coloneqq  \sqrt{1 + \frac{z_{m-1, 1-\alpha}}{\text{ESS}(x)}}.
    \end{equation}
    The type I error is then given by ${\mathbb P}(\hat R(x)\geq R_{\lim,\alpha}(x)) \simeq \alpha$.
    As an illustration, some values of $\alpha$ are reported for the threshold ${R}_{\lim,\alpha}(x) = 1.01$, $m=4$ chains and different values of $\text{ESS}(x)$ in the left panel of Table~\ref{tab:rhat_tab}.
    For example, it appears that the probability of having $\hat{R}(x) > 1.01$ and $\text{ESS}(x) = 400$ when convergence is reached is $0.04$, and decreases quickly for larger values of $\text{ESS}(x)$.
    
    \begin{table}
        \centering
        \begin{tabular}{cccc}
            \toprule
            $m$ & ${R}_{\lim,\alpha}(x)$ & $\text{ESS}(x)$ & $\alpha$ \\
            \midrule
            \multirow{6}{*}{4} & \multirow{6}{*}{1.01} & 50 & 0.80\\
            & & 100 & 0.57\\
            & & 200 & 0.26\\
            & & 400 & 0.04\\
            & & 800 & $< 10^{-3}$\\
            & & 1500 & $< 10^{-6}$\\
            \bottomrule
        \end{tabular}\qquad \qquad 
        \begin{tabular}{cccc}
            \toprule
            $m$ & ${R}_{\lim,\alpha}(x)$ & $\text{ESS}(x)$ & $\alpha$ \\
            \midrule
            2 & 1.005 & \multirow{6}{*}{400} & \multirow{6}{*}{0.05} \\
            4 & 1.010 & & \\
            8 & 1.017 & & \\
            15 & 1.029 & &  \\
            50 & 1.080 & &  \\
            100 & 1.144 & & \\
            \bottomrule
        \end{tabular}
        \caption{Left: Type I error $\alpha$ as a function of $\text{ESS}(x)$ when
         $R_{\lim,\alpha}(x) = 1.01$ and $m = 4$. Right: ${R}_{\lim,\alpha}(x)$ as a function of $m$ when $\text{ESS}(x) = 400$ and $\alpha = 0.05$.}
        \label{tab:rhat_tab}
    \end{table}
    
    \subsection{Threshold elicitation}
    \label{subsec:threshold}
    
    \paragraph{Threshold for the local $\hat{R}(x)$.}
    Proposition~\ref{prop:asymptotic_R-hat}(ii) allows us to associate a threshold for $\hat{R}(x)$ to a type I error $\alpha$, using the definition of $R_{\lim,\alpha}(x)$ in (\ref{eq-Rlimite}).
    Some values are displayed in the right panel of Table~\ref{tab:rhat_tab} for a fixed $\text{ESS}(x) = 400$ and $\alpha = 0.05$.
    It appears that the value of $1.01$, the recent recommendation of \cite{Vehtari}, seems to be coherent for $\hat{R}(x)$ and a moderate number of chains, typically the default configuration in Stan ($m=4$), JAGS ($m=3$) or PyMC3 ($m= \max\{n_{c}, 2\}$ with $n_c$ the number of cores).
    However, the value of $m$ must be doubled if a split version is used, 
    and when $m$ increases the threshold becomes more severe and it may be appropriate to consider a higher (i.e. less stringent) one: for example, a threshold of $1.1$ can be enough provided the number of chains $m$ is larger than $100$.
    \gr{The case of a large number of chains has been recently studied by \cite{margossian} who suggest a new version of $\hat{R}$ for this configuration.}
    Note that a similar observation about the stringency of the threshold can be made with rank-$\hat{R}$, see Appendix~\ref{app:r_inf_threshold} for more details.

    Therefore, we recommend to keep the threshold of $1.01$ as a general rule of thumb for $\hat{R}(x)$, except if the number of chains is too large or if one wants to have a more precise threshold.
    In such a case it only requires to provide $\alpha$, $m$ and a target value ESS$(x)$ to compute $R_{\lim,\alpha}(x)$ using (\ref{eq-Rlimite}).
    
    \paragraph{Threshold for the supremum $\hat{R}_\infty$.}
    Proposition~\ref{prop:asymptotic_R-hat} does not induce any threshold \gr{$R_{\infty,\lim}$} for $\hat{R}_\infty$, since Proposition~\ref{prop:asymptotic_R-hat}(ii) only establishes the pointwise convergence of the empirical process $\hat{R}(\cdot)$. 
    However, Proposition~\ref{prop:asymptotic_R-hat}(i) shows that under the null hypothesis where all chains follow a common distribution $F$, the latter $F$ is irrelevant to the $\hat R_\infty$ statistic. 
    Such an independence to the underlying distribution $F$ makes it possible the use of a quantile of $\hat{R}_\infty$ as a threshold associated with a given probability $\alpha$ and number of chains $m$.
    Table~\ref{tab:rhat_inf_tab} provides estimations \gr{of $R_{\infty,\lim}$} using replications for several values of $\alpha$ and $m$ and a fixed number of effective samples of $400$, as recommended by \cite{Vehtari} (more details are provided in Appendix~\ref{app:r_inf_threshold}).
    Here, we can see that a fixed rule of thumb for a range of $m$ would be too imprecise, as the quantile values increase rapidly with $m$.
    Nevertheless, Table~\ref{tab:rhat_inf_tab} illustrates a linear relationship between $m$ and the appropriate threshold for a given $\alpha$.
    
    In the simulations in Section~\ref{subsec:toy_examples} and in the experiments in Section~\ref{sec:simulations}, we mostly consider $m=4$ and therefore choose a threshold of $1.02$, which is a little more accurate than $1.01$ by looking at Table~\ref{tab:rhat_inf_tab}.
    \gr{Note that if $m=8$ or if a split version of $\hat{R}_\infty$ is used with $m=4$, then a threshold of $1.03$  should be preferred.}
    In the \texttt{localrhat} R package \citep{moinsLocalrhat}, the computation of $\hat{R}_\infty$ comes with the associated threshold at $5\%$ based on the calculations in Table~\ref{tab:rhat_inf_tab}, as well as a p-value associated with the obtained $\hat{R}_\infty$.
    
    \begin{table}
    \centering
    \begin{tabular}{cc}
        \begin{tabular}{ccccc}
            \toprule
            & \multicolumn{4}{c}{\gr{$R_{\infty,\lim}$}}\\
            \midrule
            \diagbox[height=0.6cm, width=1.1cm]{$m$}{$\alpha$} & 0.005 & 0.01 & 0.05 & 0.1 \\
            \midrule
            2  & 1.018 & 1.016 & 1.012 & 1.010 \\
            3  & 1.023 & 1.022 & 1.016 & 1.014 \\
            4  & 1.027 & 1.025 & 1.020 & 1.018 \\
            8  & 1.038 & 1.037 & 1.031 & 1.028 \\
            10 & 1.043 & 1.041 & 1.036 & 1.033 \\
            20 & 1.080 & 1.076 & 1.062 & 1.056 \\
            \bottomrule
        \end{tabular} & \raisebox{-.5\height}{\includegraphics[trim = 0cm 0.5cm 0cm 2.5cm, clip, width = .5\textwidth]{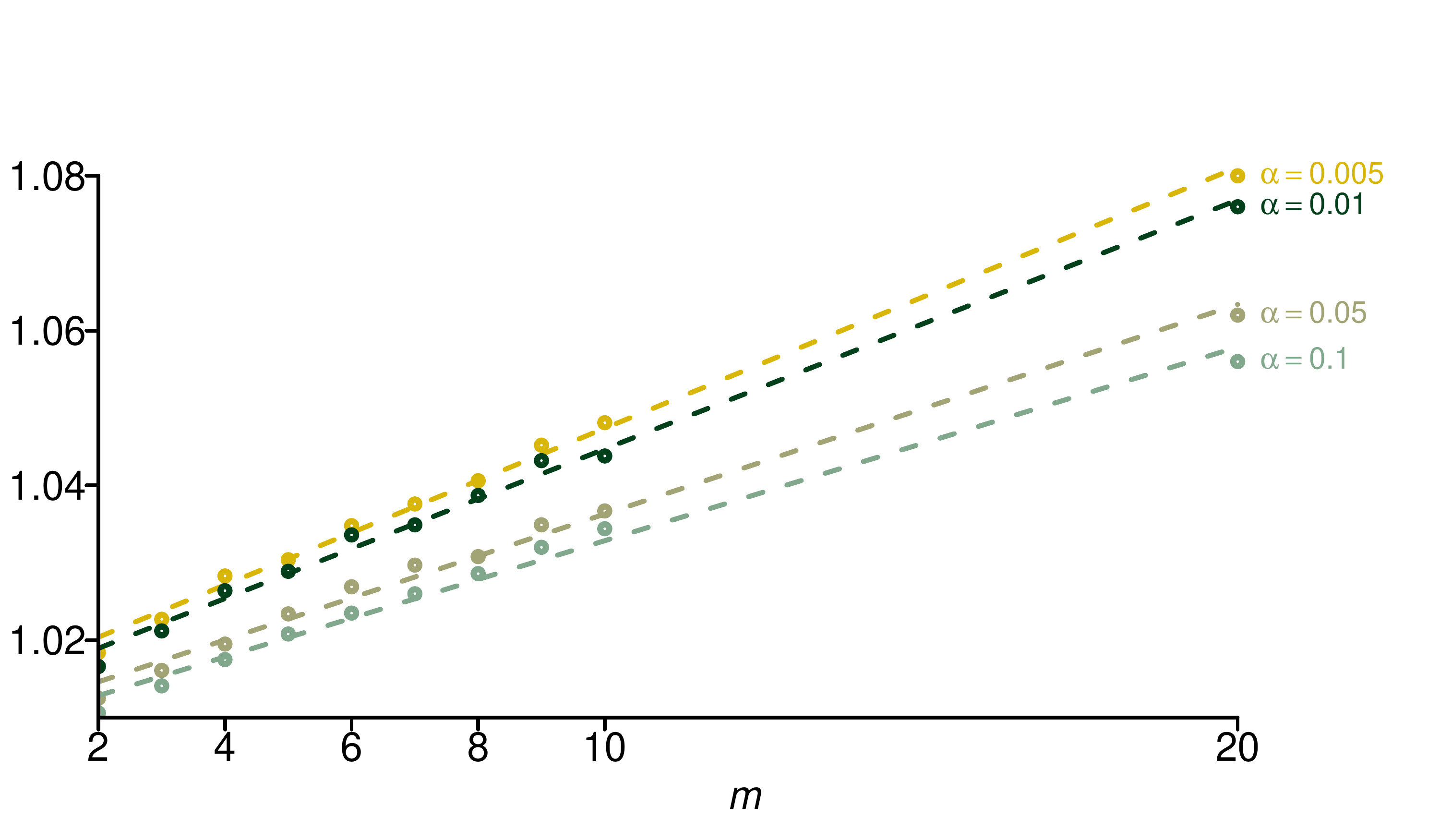}}
    \end{tabular}
    \caption{Left: Empirical quantiles \gr{$R_{\infty,\lim}$} of the $\hat{R}_\infty$ distribution under the null hypothesis (i.e. all chains follow the same distribution) for a target effective sample size of 400, based on 2000 replications.
    This table is illustrated by the figure on the right.}
    \label{tab:rhat_inf_tab}
    \end{table}

     \subsection{\gr{Illustrative} examples}
    \label{subsec:toy_examples}
    
    In this section, we consider toy distributions for the chains, where the computation of $R_\infty$ can be done explicitly.
    In particular, we first focus on two cases raised by \cite{Vehtari} of deficient behaviour of the traditional $\hat{R}$.
    Then, we exhibit a failure situation for rank-$\hat{R}$.
    All these theoretical behaviours are illustrated on a simulation study.
    Further applications to Bayesian inference are provided in  Section~\ref{sec:simulations}, and other examples where $\hat{R}$ and rank-$\hat{R}$ fail in Appendix~\ref{app:example_theorique}.

    \begin{figure}
        \hspace{-0.5cm}
        \setlength\tabcolsep{2pt}
        \begin{tabular}{ccc}
        \quad {\textbf{\textsf{Example~\hyperref[fig:example_known_dist]{1}}}}:&
        \quad {\textbf{\textsf{Example~\hyperref[fig:example_known_dist]{2}}}}:&
        \quad {\textbf{\textsf{Example~\hyperref[fig:example_known_dist]{3}}}}:\\
        \quad $\mathcal{U}(-\sigma,\sigma)$ - $\mathcal{U}(-\sigma_m,\sigma_m)$ &
        \quad Pareto$(\alpha,\eta)$ - Pareto$(\alpha,\eta_m)$ &
        \quad Exp$(1)$ - $\mathcal{U}\left(1-\frac{\lambda}{2}, 1+\frac{\lambda}{2}\right)$\\
        \vspace{-.3cm}\includegraphics[trim=0cm 1cm 0.5cm 1cm, clip, width = .33\textwidth]{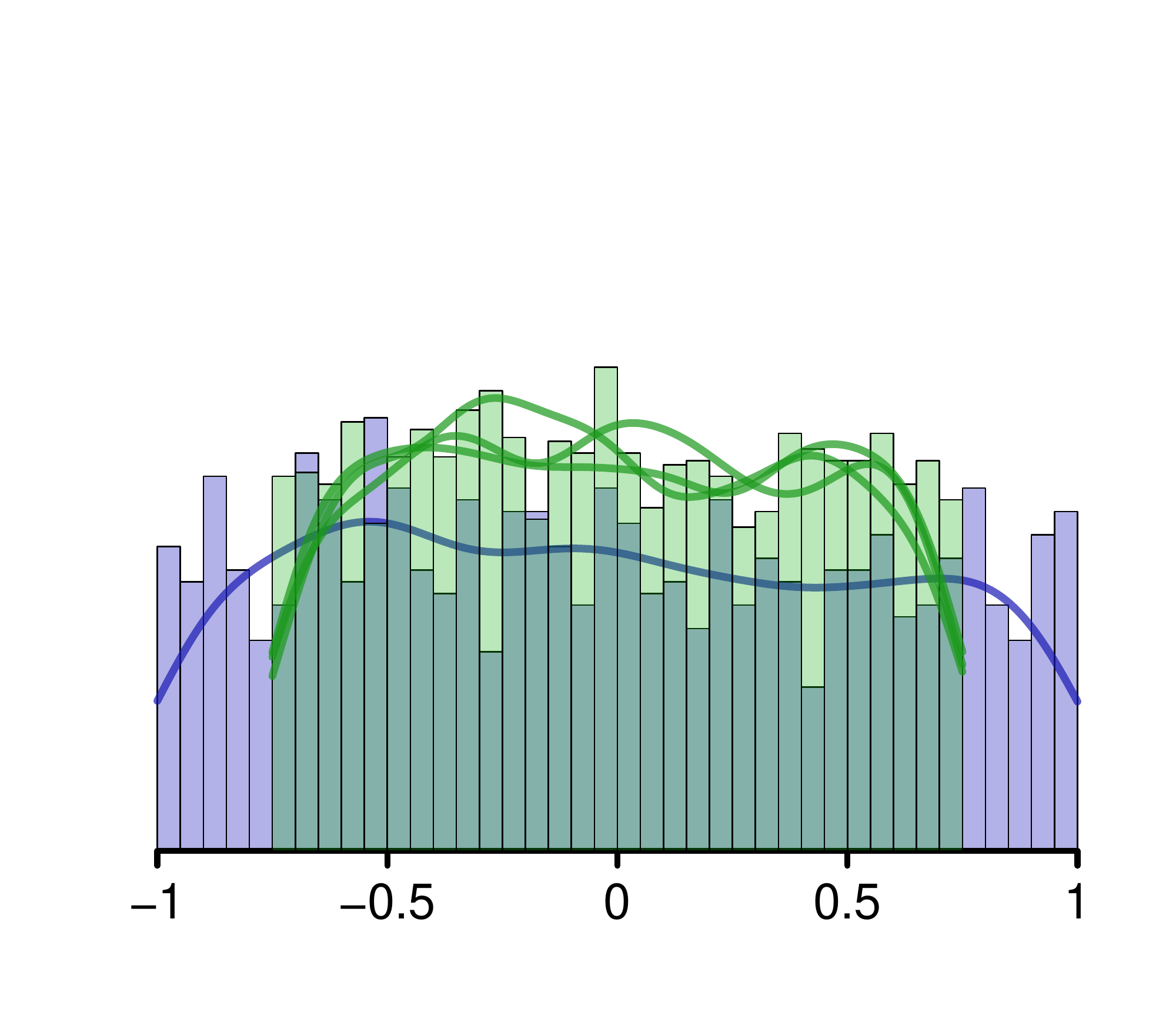} &
        \includegraphics[trim=0cm 1cm 0.5cm 1cm, clip, width = .33\textwidth]{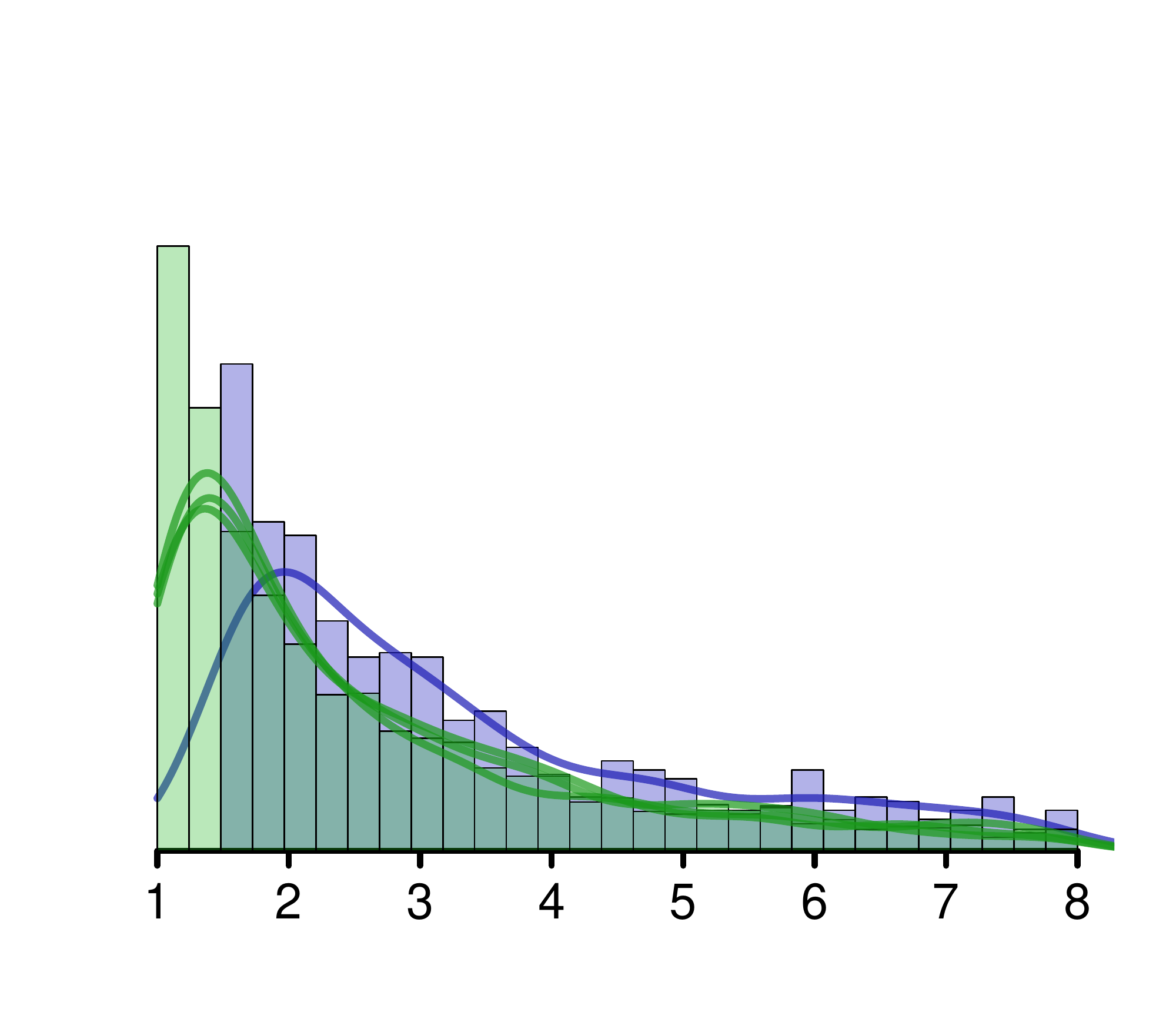} &
        \includegraphics[trim=0cm 1cm 0.5cm 1cm, clip, width = .33\textwidth]{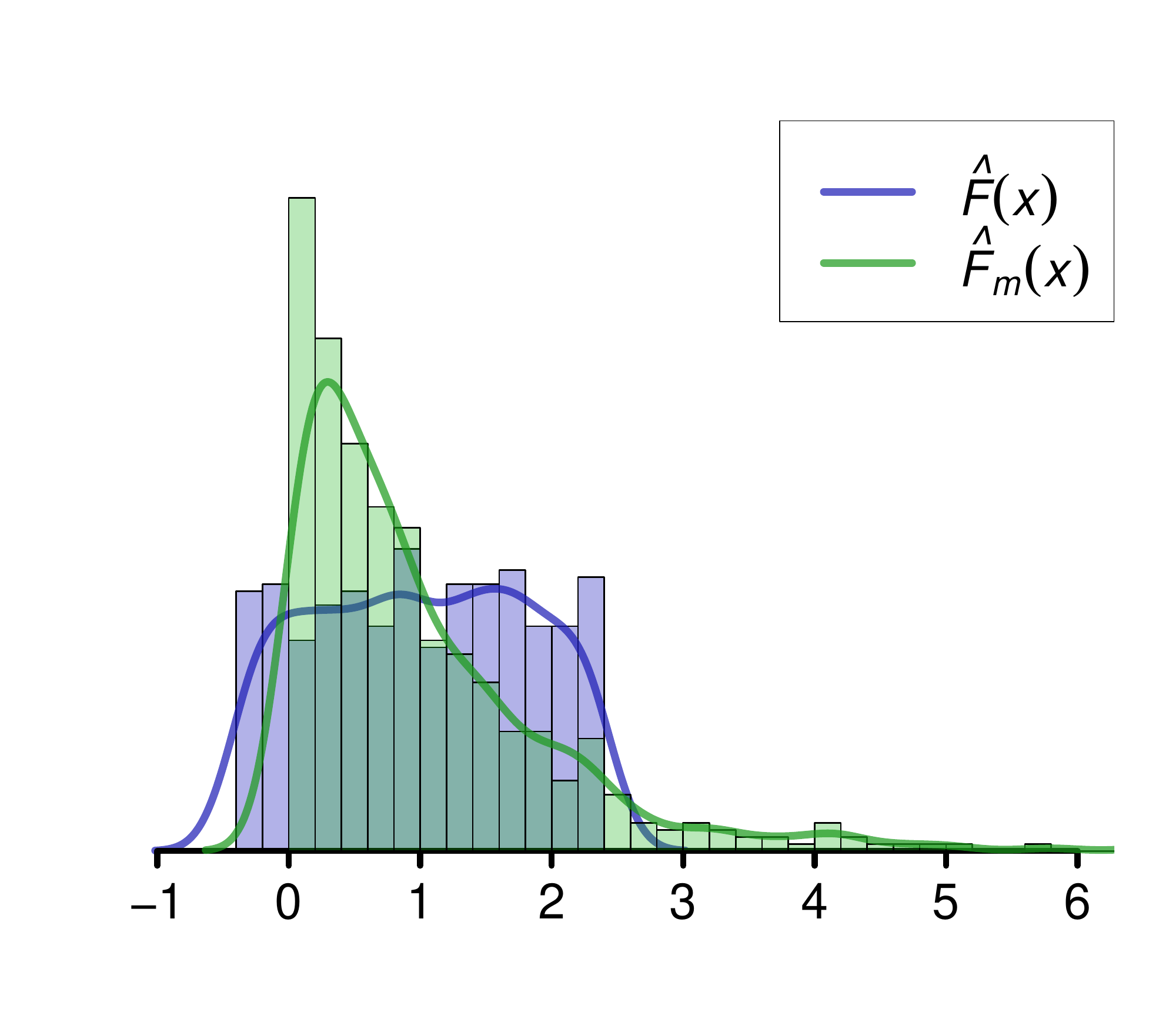}\\
        \vspace{-.3cm}\scriptsize{\textsf{x}} & \scriptsize{\textsf{x}} & \scriptsize{\textsf{x}} \\
        \includegraphics[trim=0 0 0.5 0.5cm, clip, width = .33\textwidth]{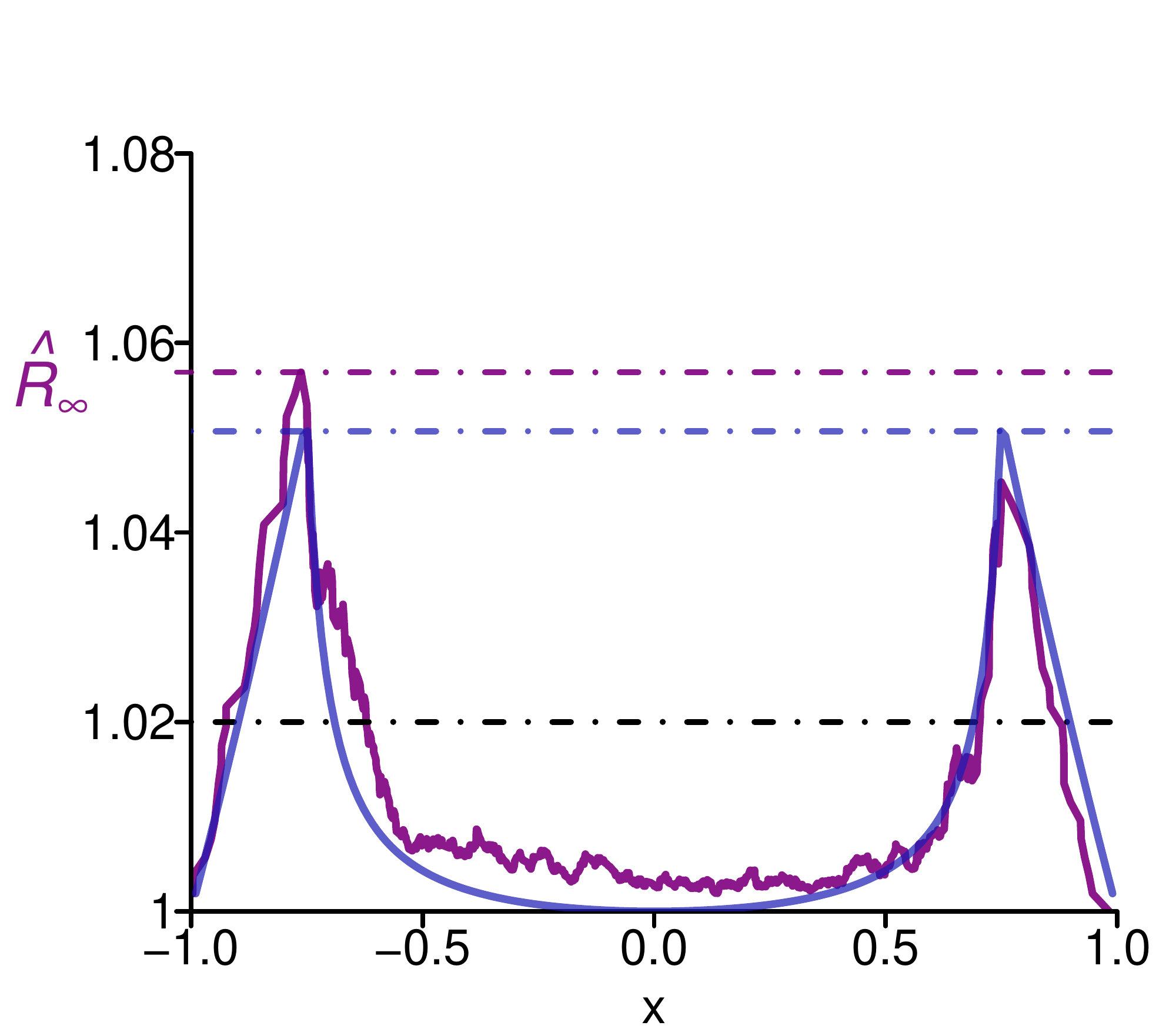} &
        \includegraphics[trim=0 0 0.5 0.5cm, clip, width = .33\textwidth]{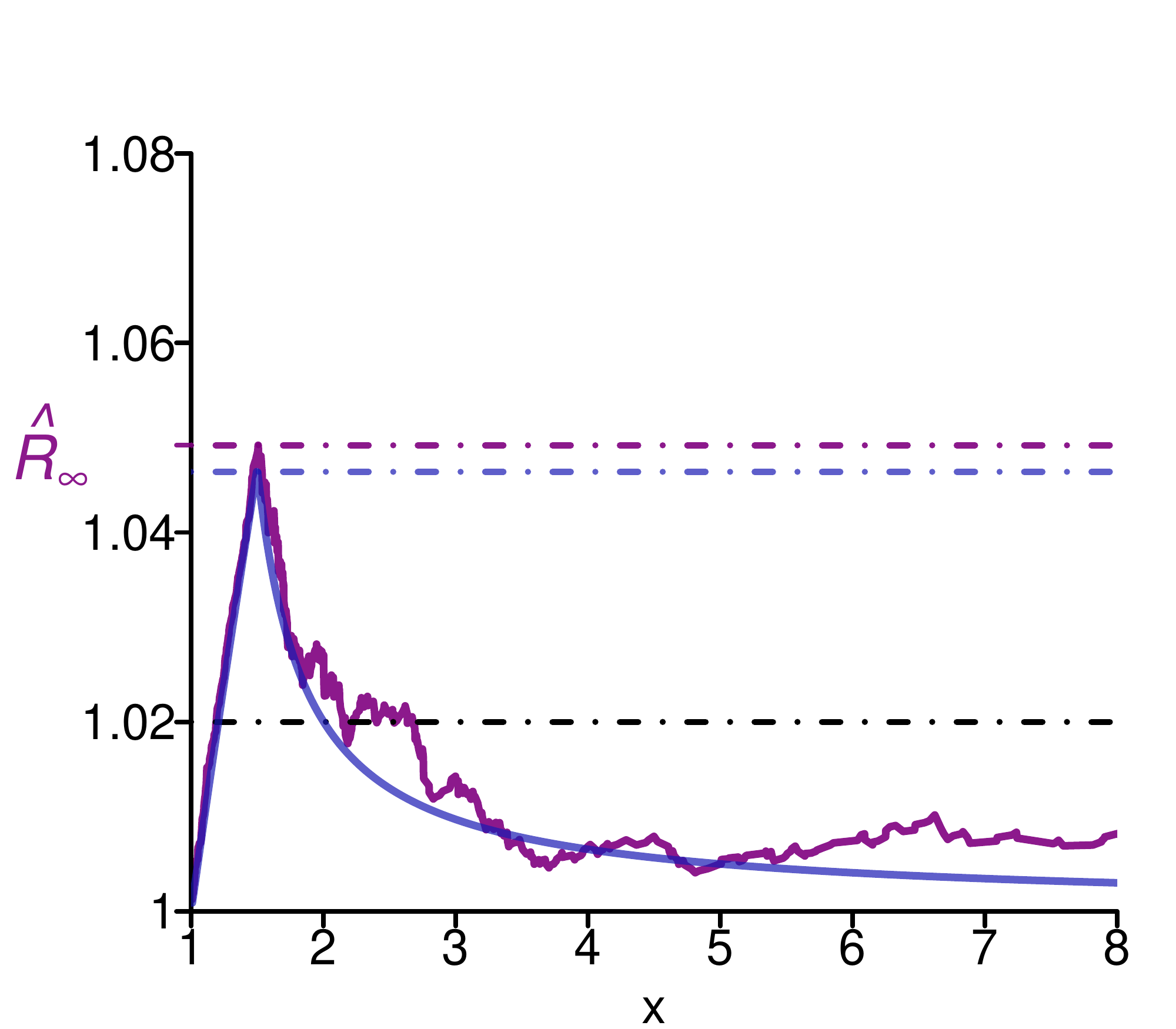} &
        \includegraphics[trim=0 0 0.5 0.5cm, clip, width = .33\textwidth]{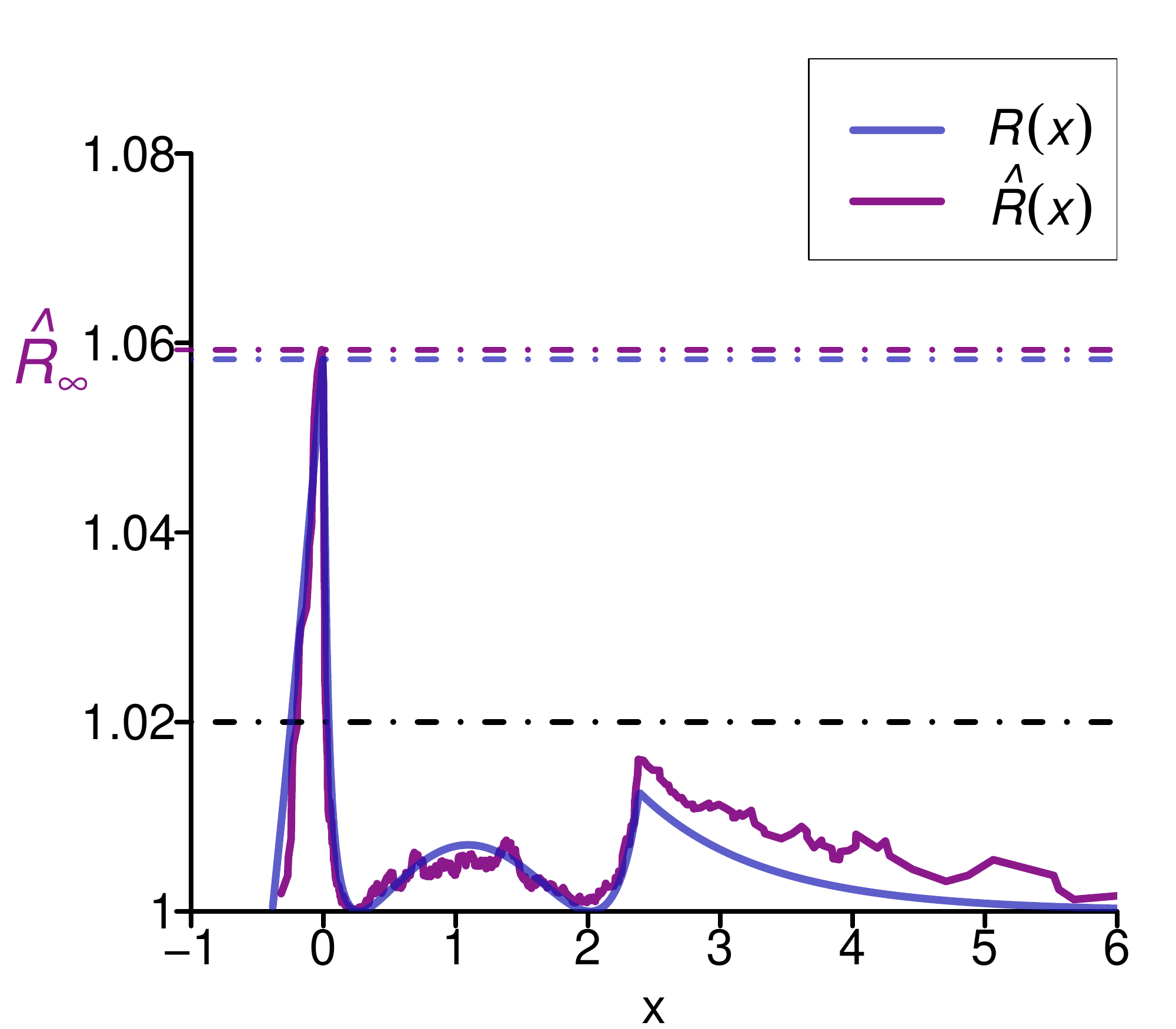}\\
        \vspace{-.2cm}
        \includegraphics[trim=0cm 1.cm 0.5cm 0.5cm, clip, width = .33\textwidth]{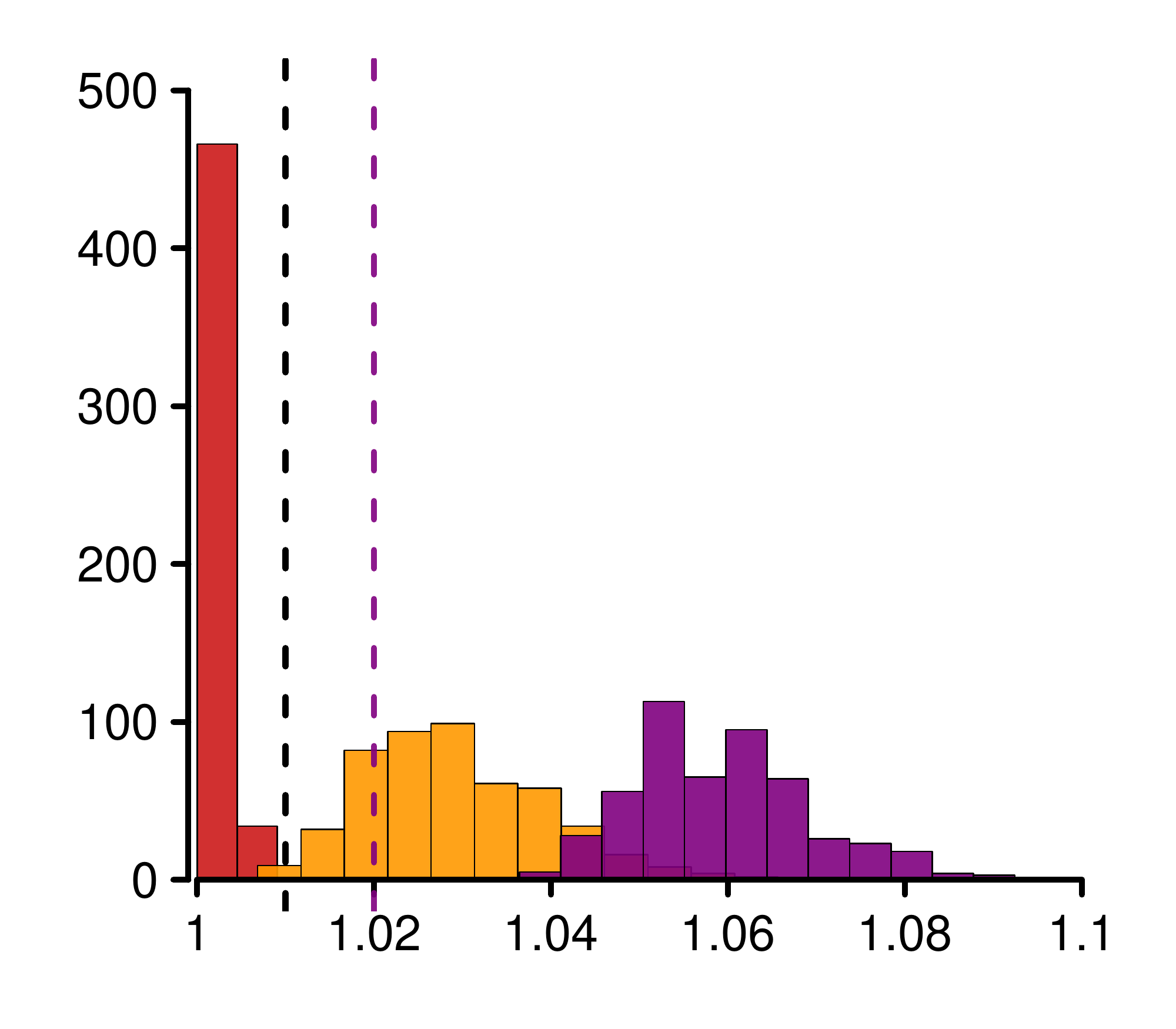} &
        \includegraphics[trim=0cm 1.cm 0.5cm 0.5cm, clip, width = .33\textwidth]{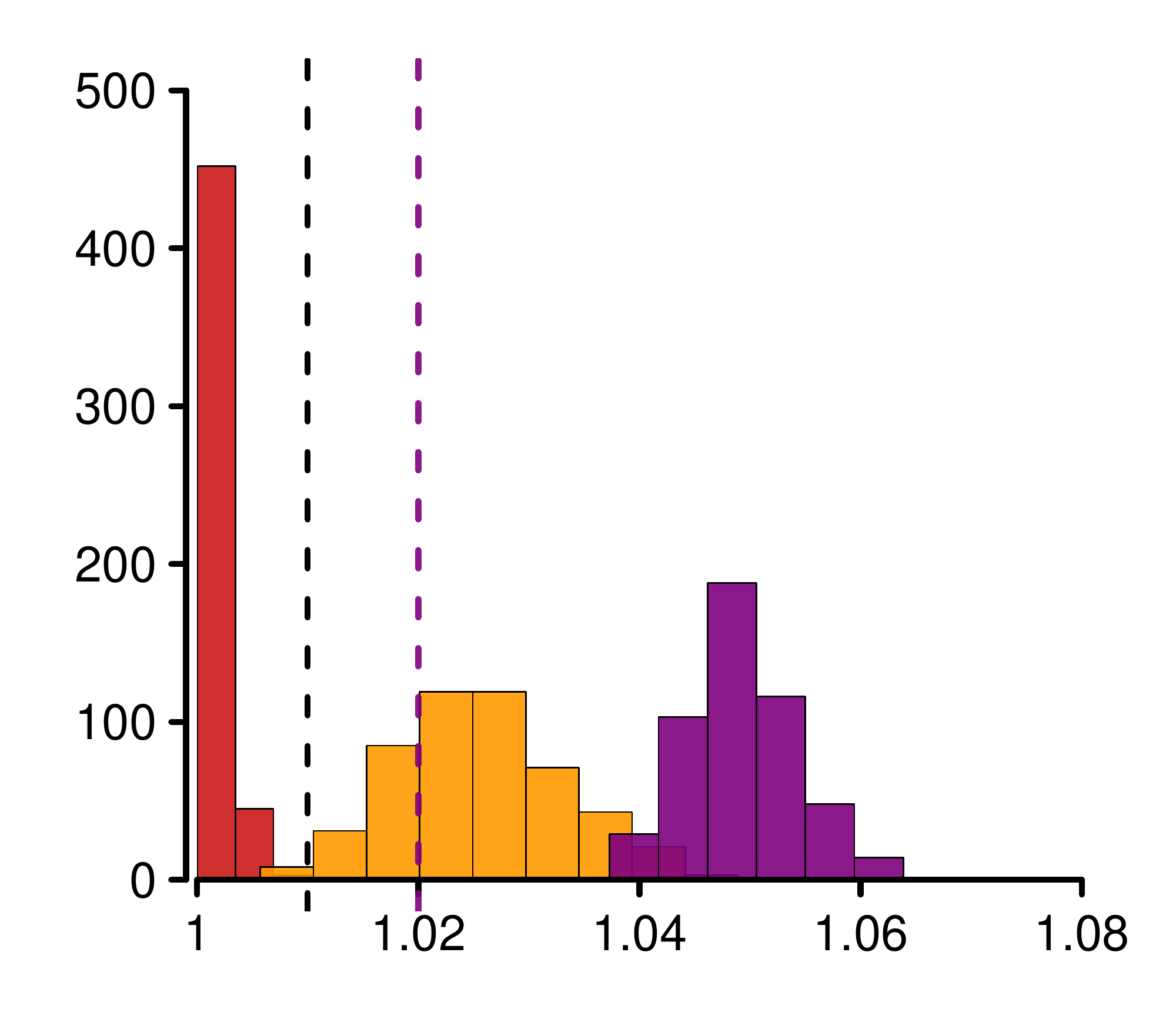} &
        \includegraphics[trim=0cm 1.cm 0.5cm 0.5cm, clip, width = .33\textwidth]{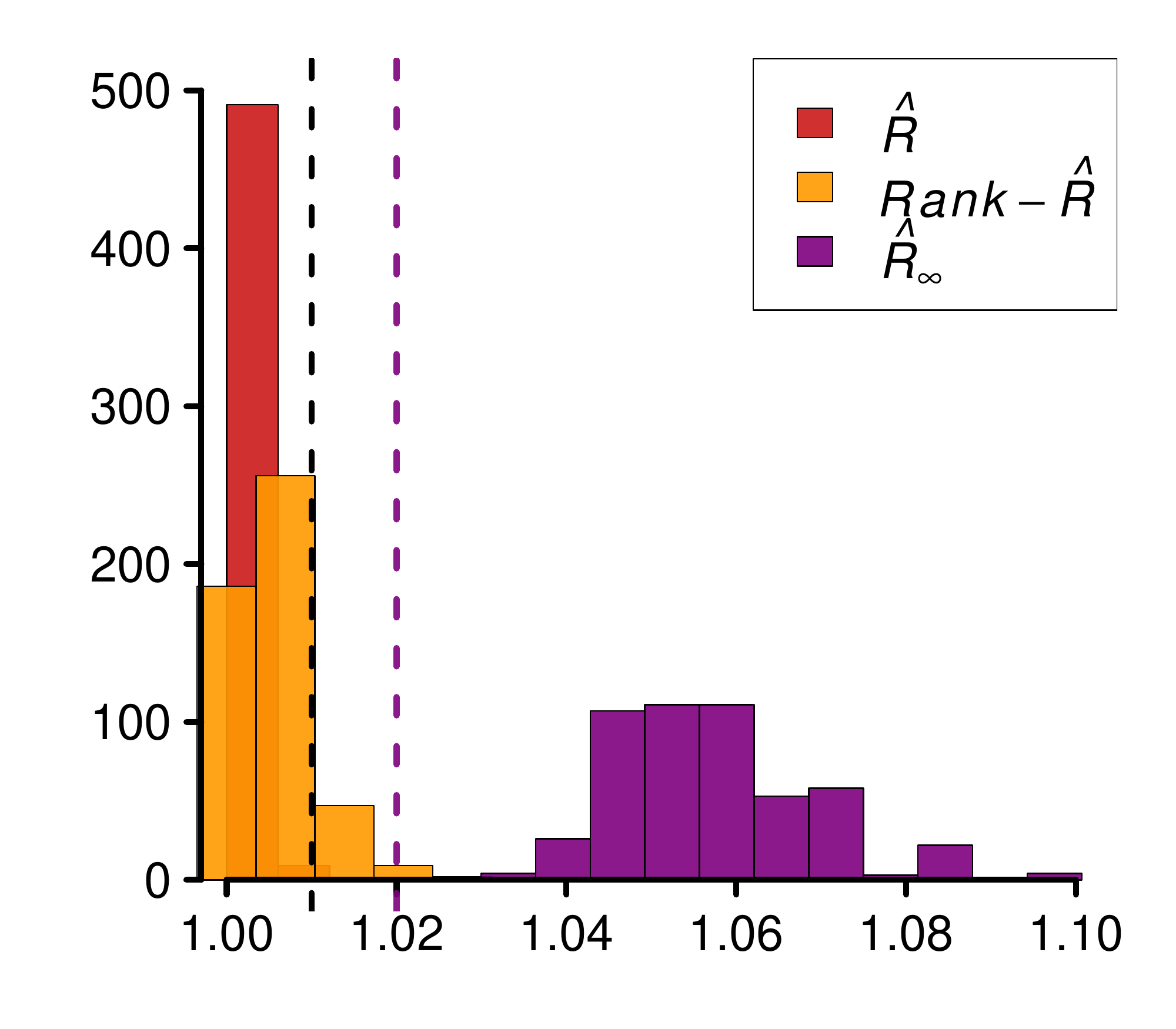}\\
        \vspace{-.3cm}
        \hspace{-.5cm}\footnotesize{\color{orange2}$\hat{R}$} \hspace{.3cm} \footnotesize{\color{red2}rank-$\hat{R}$} \hspace{.3cm} \footnotesize{\color{violet2}$\hat{R}_\infty$} &
        \hspace{-.5cm}\footnotesize{\color{orange2}$\hat{R}$} \hspace{.3cm} \footnotesize{\color{red2}rank-$\hat{R}$} \hspace{.3cm} \footnotesize{\color{violet2}$\hat{R}_\infty$} &
        \hspace{-.7cm}\footnotesize{\color{orange2}$\hat{R}$} \hspace{-.1cm} \footnotesize{\color{red2}rank-$\hat{R}$} \hspace{.3cm} \footnotesize{\color{violet2}$\hat{R}_\infty$} 
        \end{tabular}
        \caption{
        Illustrations with $m=4$ chains, $n=200$ independent iterations each.
        Top row: Simulation of $F_1 = \cdots = F_{m-1}$ in {\color{forestgreen} green} distinct from $F_m$ in {\color{darkblue} blue}. For the uniform example (left), {\color{forestgreen}$\sigma = 3/4$} and {\color{darkblue}$\sigma_m = 1$}, for the Pareto (middle) {\color{forestgreen}$\eta = 1$} and {\color{darkblue}$\eta_m = 1.5$}, and for the \gr{uniform (right) $\lambda=4\log(2)$}.
        Second row: The corresponding population version {\color{darkblue} $R(x)$} and empirical version {\color{violet2} $\hat{R}(x)$} as functions of $x$ for one replication.
        Bottom row: Histograms of $500$ replications of {\color{orange2}$\hat{R}$}, {\color{red2}rank-$\hat{R}$} and {\color{violet2}$\hat{R}_\infty$}. Dashed lines correspond to the threshold of $1.01$ for {\color{orange2}$\hat{R}$} and {\color{red2}rank-$\hat{R}$} and $1.02$ for {\color{violet2}$\hat{R}_\infty$} (see Section~\ref{subsec:convergence-rhat}).
        }
        \label{fig:example_known_dist}
    \end{figure}
    
    \paragraph{Example~\hyperref[fig:example_known_dist]{1}: Chains with same mean and different variances.}
    \label{subsec:uniform_example}
    
    To tackle the first situation of poor behaviour of the traditional $\hat{R}$, we consider $m$ chains following centered uniform distributions
    with different variances. More specifically,
    assume that the $m-1$ first chains have the cdf $F_1=\dots=F_{m-1}$ of the uniform distribution ${\mathcal U}(-\sigma,\sigma)$ while the last chain has the cdf $F_m$ of the uniform distribution ${\mathcal U}(-\sigma_m,\sigma_m)$ with $0<\sigma\leq\sigma_m$.
    In such a case, the between-variance is zero and it is thus expected that $\hat R\approx 1$.
    In contrast, Lemma~\ref{lem-unif} in Appendix~\ref{app:example_theorique} provides an explicit expression for $R(x)$ as well as 
    \begin{equation*}
        R_\infty =
        \sqrt{1 + \frac{m-1}{m}\left(1-\frac{2}{1+{\sigma_m}/{\sigma}}\right)}.
    \end{equation*}
    It appears that $R_\infty$ is an increasing function of $\sigma_m/\sigma $ starting from $R_\infty = 1$ when $\sigma_m/\sigma=1$, and upper-bounded by $\sqrt{2-1/m}$ when $\sigma_m/\sigma\to\infty$.
    Results are illustrated in the left column of Figure~\ref{fig:example_known_dist}.
    In the bottom row, the histograms of replications confirm that $\hat{R}_\infty$ is able to spot the same convergence issue as the one \cite{Vehtari} suggests.
    
    \paragraph{Example~\hyperref[fig:example_known_dist]{2}: Chains with heavy-tails and different locations.}
    \label{subsec:pareto_example}

    As a second example of poor behaviour of $\hat{R}$, we consider chains following Pareto$(\alpha,\eta)$ distributions, with cdf 
    $$
    F(x \mid \alpha, \eta) = 1 - \left({x}/{\eta}\right)^{-\alpha}, \quad \forall x \in [\eta, +\infty),
    $$
    shape parameter $\alpha>0$ and lower bound $\eta>0$. 
    Let us recall that such a distribution is heavy-tailed~\cite[Table~3.4.2]{embrechts2013modelling} and has an infinite first moment when $\alpha \leq 1$.
    We focus on the case where one chain is shifted from the other ones: $F_1(x) = \cdots = F_{m-1}(x) = F(x \mid \alpha, \eta)$ and $F_m(x) = F(x \mid \alpha, \eta_m)$ with $0<\eta \leq \eta_m$ and $\alpha \leq 1$. Here, the within- and between-variances do not exist and  it is expected in practice that $\hat{R} \approx 1$. In contrast, $R_\infty$ can be written as
    \begin{equation*}
        R_\infty = \sqrt{1+  \frac{1}{m}\left(\left(\frac{\eta_m}{\eta}\right)^{\alpha}-1\right)},
    \end{equation*}
    see Lemma~\ref{lem-pareto} in the supplementary material. Clearly, $R_\infty$ is an increasing function of $\eta_m/\eta$ starting from $R_\infty=1$ when $\eta_m=\eta$
    and such that $R_\infty\to \infty$ as $\eta_m/\eta\to\infty$.
    Results are shown in the middle column of Figure~\ref{fig:example_known_dist}. 
    This experiment corresponds to the second example of convergence issue raised by \cite{Vehtari}. 
    The same observations as for Example~\hyperref[fig:example_known_dist]{1} can be made here: $\hat{R}_\infty$ is prone to indicating a convergence issue than rank-$\hat{R}$.
    
    \paragraph{Example~\hyperref[fig:example_known_dist]{3}: Chains with same mean and mean over the median.}
      
    Finally, we come back to the example described in Section~\ref{subsec:intro_comparison_rhat} where both $\hat{R}$ and rank-$\hat{R}$ fail to detect non-convergence.
    Following the method described in Appendix~\ref{app:counter_example}, we consider $m-1$ exponential chains Exp$(1)$ and one uniform $\mathcal{U}(1 - 2\log 2, 1 + 2\log 2)$.
    This results in chains with same mean and mean over the median.
    Results are illustrated in the right panel of Figure~\ref{fig:example_known_dist}: the histograms of replications confirm that $\hat{R}_\infty$ is able to detect the convergence issue that neither $\hat{R}$ nor rank-$\hat{R}$ are able to detect.
    Here, the explicit calculation of $R_\infty$ is not feasible, but Lemma~\ref{lem:laplace_uniform_proof} in the supplementary material provides another example where the computation can be done, with uniform and Laplace distributions.

\section{Multivariate extension}
\label{sec:multivariate}

    \subsection{Population version and algorithm for multivariate diagnosis}
    \label{subsec:multivariate_def}
    
    Our $R(x)$ can naively be adapted to the multivariate case:
    assume now that the parameter is multivariate and write $\boldsymbol{\theta}=(\theta_1, \ldots, \theta_d)\in{\mathbb R}^d$ with $d\geq 2$, and  denote by $\theta_p^{(j)}$ the coordinate $p\in\{1,\dots,d\}$ from chain $j\in\{1,\dots,m\}$.
    Similarly to the univariate case, $\hat R$ can be computed on the indicator variables $I_{\boldsymbol{x}}^{(j)} = \mathbb{I}\{\theta_1^{(j)} \leq x_1, \ldots, \theta_d^{(j)} \leq x_d\}$ for any $\boldsymbol{x} = (x_1, \ldots, x_d) \in \mathbb{R}^d$.
    Under the assumptions of Proposition~\ref{prop-calcul-R}, all calculations remain valid in dimension $d$ and therefore the expression of $R(\boldsymbol{x})$
    is formally the same as in~(\ref{eq:R_theorique}):
    \begin{equation}
        R(\boldsymbol{x}) = \sqrt{\frac{W(\boldsymbol{x})+B(\boldsymbol{x})}{W(\boldsymbol{x})}}
        = \sqrt{1 + \frac{\sum_{j=1}^m\sum_{k=j+1}^m \left(F_j(\boldsymbol{x})-F_k(\boldsymbol{x})\right)^2}{m\sum_{j=1}^m F_j(\boldsymbol{x})(1-F_j(\boldsymbol{x}))}}.
        \label{eq:R_theorique_multi}
    \end{equation}
    The properties listed in Proposition~\ref{prop-first} in the univariate case remain true as well.
    The associated $R_\infty$ is defined as $R_\infty(F_1,\dots,F_m)=\sup_{\boldsymbol{x}\in{\mathbb R}^d} R(\boldsymbol{x})$, while $\hat R(\boldsymbol{x})$ is computed
    by replacing the cumulative distribution functions in~(\ref{eq:R_theorique_multi}) by their empirical counterparts.
    Note also that all values computed in Table~\ref{tab:rhat_tab} and Table~\ref{tab:rhat_inf_tab} remain identical in this multivariate extension. 
    \gr{However, those results are not giving information about the sensitivity to convergence issues, which in the multivariate case can come from margins but also from the dependence structure.}
    
    It is easily seen that, if the \gr{marginal distributions} of $F_1,\dots,F_m$ coincide, then $R_\infty$ is the same as the one associated with uniform margins (see Lemma~\ref{lem-invariance} in the supplementary material).
    In other words, we have $R_\infty(F_1,\dots,F_m)=R_\infty(C_1,\dots,C_m)$ where $C_j$ is the copula defined in $[0,1]^d$ associated with $F_j$, $j\in\{1,\dots,m\}$.
    This suggests that a multivariate diagnosis can be conducted in two steps as follows: 
    \gr{\begin{enumerate}
        \item Compute the univariate $\hat{R}_{\infty,p}$ separately on each of the coordinates $p\in\{1,\dots,d\}$.
        If $\hat{R}_{\infty,p} < R_{\infty,\text{lim}}^{(M)}$ for all $p\in\{1,\dots,d\}$, with $R_{\infty,\text{lim}}^{(M)}$ a choice of margins threshold, then all of them are deemed to have converged and to be identically distributed.
        \item Compute the multivariate $\hat{R}_\infty$ to check the dependence structure convergence.
        If $\hat{R}_{\infty} < R_{\infty,\text{lim}}^{(C)}$, with $R_{\infty,\text{lim}}^{(C)}$ a copula threshold, then the dependence structure is also deemed to have converged, and so has the multivariate distribution. 
    \end{enumerate}
    The test for convergence is now separated in two parts: 1. convergence of the margins, and 2. convergence of the copula knowing that the margins have converged.
    It can easily be shown that, up to a first order approximation, one way to obtain a type I error $\alpha$ for the global two-step test is to consider a  level $\alpha/2$ for each of the two components.
    The first step corresponds to $d$ univariate tests, so for $R_{\infty,\text{lim}}^{(M)}$ one can use the univariate threshold $R_{\infty,\lim}$ defined in Section~\ref{subsec:threshold} with a  level $\alpha/2d$, corresponding to a Bonferroni correction for the error level $\alpha/2$.
    In the following subsections, we focus on the second step of the algorithm: the theoretical properties of the multivariate $\hat{R}_\infty$ in the case of convergence on the margins, which will provide insights for choosing $R_{\infty,\text{lim}}^{(C)}$. 
    Values of $R_{\infty,\text{lim}}^{(M)}$ and $R_{\infty,\text{lim}}^{(C)}$ are then given as functions of $(\alpha, d, m)$ in Table~\ref{tab:rhat_max_tab}. As a general rule, one can reasonably use for $\alpha=0.05$ the values  $(R_{\infty,\text{lim}}^{(M)}, R_{\infty,\text{lim}}^{(C)}) = (1.03, 1.03)$ in the case of $m=4$ chains, and $(R_{\infty,\text{lim}}^{(M)}, R_{\infty,\text{lim}}^{(C)}) = (1.04, 1.05)$ if $m=8$ or if a split version is used with $m=4$, with limited variations around these values for varying dimension $d$.
    }
    
    \subsection{Upper bounds}

    Let us first consider the case of $m=2$ chains with uniform margins and associated copulas $C_1$ and $C_2$.
        For all  $\boldsymbol{u} = (u_1, \ldots, u_d)\in[0,1]^d$, one has
        \begin{equation}
            R(\boldsymbol{u})
            = \sqrt{ 1 + \frac{\left(C_1(\boldsymbol{u})-C_2(\boldsymbol{u})\right)^2}{2\left(C_1(\boldsymbol{u})(1-C_1(\boldsymbol{u})) + C_2(\boldsymbol{u})(1-C_2(\boldsymbol{u}))\right)}}.
            \label{eq:R_theorique_2D}
        \end{equation}
        In addition to having the usual lower bound of $1$, the next lemma allows establishing an upper bound on $R_\infty(C_1,C_2)$.
        \begin{Lem}
        \label{lem:copule_bound}
        Let $C_1, C_2, C_-$ and $C_+$ be copulas such that:
          \begin{equation}
            \text{for all } \boldsymbol{u} \in [0,1]^d,
            \begin{cases}
                C_{-}(\boldsymbol{u}) \leq C_1(\boldsymbol{u}) \leq C_{+}(\boldsymbol{u}), \\
                C_{-}(\boldsymbol{u}) \leq C_2(\boldsymbol{u}) \leq C_{+}(\boldsymbol{u}).
            \end{cases} 
            \label{eq:bound_def}
          \end{equation}
            Then, $R_\infty(C_1, C_2) \leq R_\infty(C_{-}, C_{+})$.
        \end{Lem}
        Let $W_d$ and $M_d$ the lower and upper Fréchet--Hoeffding bounds in dimension $d$ \citep[see][Theorem~2.10.12]{Nelsen2006}:
        \begin{equation*}
            W_d(\boldsymbol{u}) \coloneqq  \max\left\{ 1 - d + \sum_{i=1}^d u_i, 0 \right\} 
            \quad \text{and} \quad 
            M_d(\boldsymbol{u}) \coloneqq  \min\left\{ u_1, \ldots, u_d \right\}.
        \end{equation*}
        Any copula is bounded from below and from above by $W_d$ and $M_d$ respectively, in the sense of (\ref{eq:bound_def}).
        Thus, applying Lemma~\ref{lem:copule_bound} with $(C_-, C_+) = (W_d, M_d)$ yields:
        \begin{Prop}
            \label{prop:global_bound}
            For any $d$-variate copulas $C_1$ and $C_2$, 
            $$R_\infty(C_1,C_2) \leq \sqrt{\frac{d+1}{2}},$$
        \end{Prop}
        Unlike the univariate version (see for instance Example~\hyperref[fig:example_known_dist]{2} in Section~\ref{subsec:toy_examples}), the value of $R_\infty$ associated with the convergence of the dependence structure is upper-bounded, with a bound that grows with the dimension.
        This difference of behaviour could be used for example to tune the threshold for the multivariate case.
        However this bound, although it is the ``best possible'' \citep[Theorem~2.10.13]{Nelsen2006}, is tight only in the case $d=2$ since
        $W_d$ is no more a copula when $d>2$. It may also be too loose since it compares the extreme case of one chain with comonotonic dependence and another one with anti-comonotonic dependence. Some refinements are proposed in Section~\ref{sub-direction}.
        
        In the case of $m>2$ chains, the previous bounding technique does not apply anymore, and we propose the following result based on bounding
        pairwise $R_\infty$'s:
        \begin{Cor}             
        \label{cor:global_bound_d}
            For any $m \geq 2$ and $d$-variate copulas $(C_1, \ldots, C_m)$, 
            \begin{equation*}
                R_\infty(C_1, \ldots, C_m) \leq \sqrt{1 + \frac{m-1}{2}(d-1)}.
            \end{equation*}
        \end{Cor}
        Although this limit is not tight in the general case, it coincides with the upper bound of Proposition~\ref{prop:global_bound} when $m=2$.
        Let us also note that, for any fixed $m\geq 2$, the upper bound of $R_\infty(C_1, \ldots, C_m)$ diverges at a fixed $\sqrt{d}$ rate as the dimension increases.
        
    \subsection{Influence of the dependence direction on the sensitivity of $\hat{R}_\infty$}
    \label{sub-direction}
    
    When $m=2$, one way to refine the upper bound established in Proposition~\ref{prop:global_bound} is to assume that both copulas are modelling 
        either positive of negative dependence. 
    More specifically, let us recall the notions of positive lower orthant dependence (PLOD) and negative lower orthant dependence (NLOD) \citep[see][Section~5.7]{Nelsen2006}. The random vector $(\theta_1, \ldots, \theta_d)$ is
    \begin{itemize}
        \item PLOD if
    $
    \forall \boldsymbol{x} \in \mathbb{R}^d, \quad \mathbb{P}(\theta_1 \leq x_1, \ldots, \theta_d \leq x_d) \geq \prod_{i=1}^d \mathbb{P}(\theta_i \leq x_i),
    $
        \item NLOD if
    $
    \forall \boldsymbol{x} \in \mathbb{R}^d, \quad \mathbb{P}(\theta_1 \leq x_1, \ldots, \theta_d \leq x_d) \leq \prod_{i=1}^d \mathbb{P}(\theta_i \leq x_i).
    $
    \end{itemize}
    Both properties can be characterized in terms of the associated copula. The PLOD (resp. NLOD) property holds if and only if
    $C(\boldsymbol{u}) \geq \Pi_d(\boldsymbol{u})$ (resp. $C(\boldsymbol{u}) \leq \Pi_d(\boldsymbol{u})$) for all 
    $\boldsymbol{u} \in [0,1]^d$ where $\Pi_d$ is the independent copula defined by
     $\Pi_d(\boldsymbol{u}) \coloneqq  \prod_{i=1}^d u_i$.
    Note that this does not define a total order on copulas since some copulas are neither PLOD nor NLOD.
    Nevertheless, it allows us to derive refined bounds for $R_\infty$ in the NLOD and PLOD cases.
    
    For PLOD, the upper bound is in not closed-form for any dimension $d$, but simple bounds can be derived in the two extreme cases $d=2$ and $d\to\infty$.
    \begin{Cor}
        \label{cor:PLOD_bound}
    Let $m=2$.
    For any two PLOD $d$-variate copulas $C_1$ and $C_2$, $R_\infty(C_1,C_2) \leq R_\infty(\Pi_d, M_d)$ with
    $$
    \begin{cases}
        R_\infty(\Pi_2, M_2) = \sqrt{\frac{1}{2} + \frac{1}{\sqrt{3}}} \approx 1.038 \quad \text{if } d=2,\\
      \sqrt{\frac{d}{2\log d}}(1+o(1))  \leq R_\infty(\Pi_d, M_d) \leq  \sqrt{\frac{d+1}{2}}\quad \text{ as } d \to \infty.
    \end{cases}
    $$
    \end{Cor}
    Conversely, the upper bound can be computed explicitly in the NLOD case.
    \begin{Cor}
        \label{cor:NLOD_bound}
    Let $m=2$. For any two NLOD $d$-variate copulas $C_1$ and $C_2$, $R_\infty(C_1,C_2) \leq R_\infty(\Pi_d, W_d)$ with
    $$
    R_\infty(\Pi_d, W_d) = \sqrt{1+\frac{1}{2}\frac{1}{\left(1-\frac{1}{d}\right)^{-d}-1}}.
    $$
    \end{Cor}
    Let us stress that positive and negative dependence are handled differently by $R_\infty$. 
    When $d=2$, the PLOD and NLOD bounds (respectively equal to $1.04$ and $1.08$) are significantly lower than the value $\sqrt{3/2} \approx 1.22$ corresponding to the global bound, with a value higher in the NLOD case than in the PLOD one.
    However, this observation is quickly inverted when $d$ increases: for NLOD, $R_\infty(\Pi_d, W_d)$ is bounded and converges to $\sqrt{1+\frac{1}{2(e-1)}} \approx 1.136$ as $d\to\infty$, which strongly constrains the range of values that can be obtained whatever the dimension.
    In contrast, the upper bound $R_\infty(\Pi_d, M_d)$ in the PLOD case diverges with the dimension, at the same rate (up to a logarithmic factor) as in the general case, see Proposition~\ref{prop:global_bound}.
    Thus, the sensitivity of $R_\infty$ strongly depends on the sign of dependence and asymptotically favours PLOD dependence when $d$ increases.
    
    This difference can be explained by the construction of $R(x)$ itself (and thus $R_\infty$), which favours a dependence direction in ${\mathbb R}^d$ due to the computation of $\mathbb{I}\{\theta_1^{(\cdot)} \leq x_1, \ldots, \theta_d^{(\cdot)} \leq x_d\}$.
    One way to overcome this issue in the bivariate case is to compute two versions of $R_\infty$, denoted respectively by
    $R^{+}_\infty$ and $R^{-}_\infty$, 
    based respectively on  $\mathbb{I}\{\theta_1^{(\cdot)} \leq x_1, \theta_2^{(\cdot)} \leq x_2\}$ and $\mathbb{I}\{\theta_1^{(\cdot)} \leq x_1, \theta_2^{(\cdot)} \geq x_2\}$. 
    Note that $R^{+}_\infty$ coincides with the construction proposed in Section~\ref{subsec:multivariate_def}.
    \begin{Cor}
        \label{cor:changement_bound}
        Let $m=2$. Then,
            $R^{+}_\infty(\Pi_2, M_2) = R^{-}_\infty(W_2, \Pi_2)$ and $R^{+}_\infty(W_2, \Pi_2) = R^{-}_\infty(\Pi_2, M_2)$.
    \end{Cor}
    It appears that PLOD and NLOD upper bounds are exchanged by computing $R^{-}_\infty$ instead of $R^{+}_\infty$, which makes $R^{-}_\infty$ 
    more sensitive to negative dependence than positive dependence (in the bivariate case). 
    One way to consider symmetrically both dependencies would be to consider $\hat{R}^{(\text{max})}_\infty = \max(R^{+}_\infty, R^{-}_\infty)$.
    However, in dimension $d$, considering  all directions would imply the computation of $2^{d-1}$ different $R_\infty$, which would be too expensive for large $d$.
    Similar  curse of dimensionality occurs in the multivariate extension of the Kolmogorov–Smirnov test, see for example \cite{lopes2007two} for improvements of the naive multidimensional version of the test.
    Computing $\hat{R}^{(\text{max})}_\infty$ is still feasible for small values of $d$: typically for $d \leq 6$ we were able to replicate values in our experiments.
    Therefore, we provide in Table~\ref{tab:rhat_max_tab} (Appendix \ref{app:r_inf_threshold}) the estimated threshold \gr{$R_{\infty,\text{lim}}^{(C)}$} associated with the maximum of $\hat{R}_\infty$ in all possible directions when $d \leq 6$.

    One alternative in the high-dimensional case could be to apply $\hat{R}$ on an indicator function associated with a univariate function of the parameters, to return to the case described in Section~\ref{sec:r-hat-def}.
    Typically in a Bayesian model, one could use the log-likelihood $l_{\boldsymbol{\theta}} = \log p(y \mid \boldsymbol{\theta})$ when it is available, and compute $\hat{R}_\infty$ with $\mathbb{I}\{l_{ \boldsymbol{\theta}} \leq x\}$.
    Similarly, the log posterior as implemented in Stan can also be used, as suggested in the Stan reference manual \citep{carpenter2017stan}. 
    Ensuring convergence for all $x$ on the log posterior may be satisfying for multivariate diagnosis, as it is illustrated in Example~\ref{fig:logit}.
    
    \subsection{\gr{Multivariate illustrative} examples}
    \label{subsec:multivariate_toy_examples}
    
    Similarly to Section~\ref{subsec:toy_examples}, we illustrate our theoretical study in the multivariate case with simulations based on toy distributions for the chains.
    Especially, we consider multivariate normal distributions, and focus on the case where all the margins are the same (typically distributed according to a standard normal distribution).
    This leads to
    \begin{equation*}
        \boldsymbol{\theta}^{(i,j)} \sim \mathcal{N}(\textbf{0}, \boldsymbol{\Sigma}_j),
    \end{equation*}
    $i \in \{1, \ldots, n\}$ and $j \in \{1, \ldots, m\}$, where $\boldsymbol{\Sigma}_j$ is the covariance matrix of the chain $j$, with diagonal elements equal to one to keep standard Gaussian margins.
    
\begin{figure}
        \begin{tabular}{cc}
            \textbf{Example~\hyperref[fig:bivariate]{4}:} Bivariate normal &  Brooks--Gelman $\hat{R}$, $\hat{R}_\infty$, and $R_\infty$   wrt $\rho_m$ \\
            \vspace{-.3cm}
            \includegraphics[trim=0cm 0cm 0cm 0cm, clip, width = .48\textwidth]{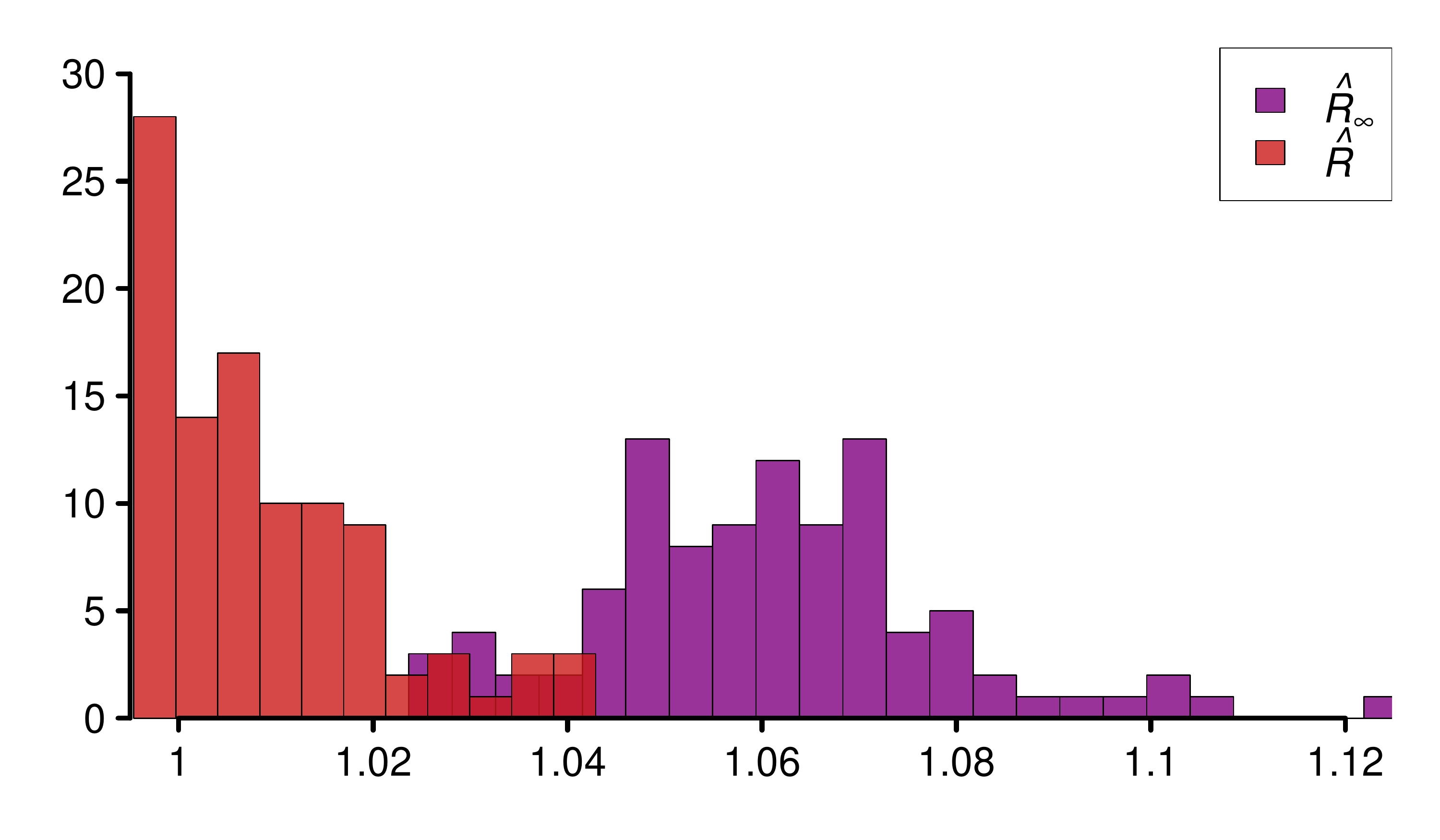} & 
            \includegraphics[trim=0cm 0cm 0cm 0cm, clip, width = .48\textwidth]{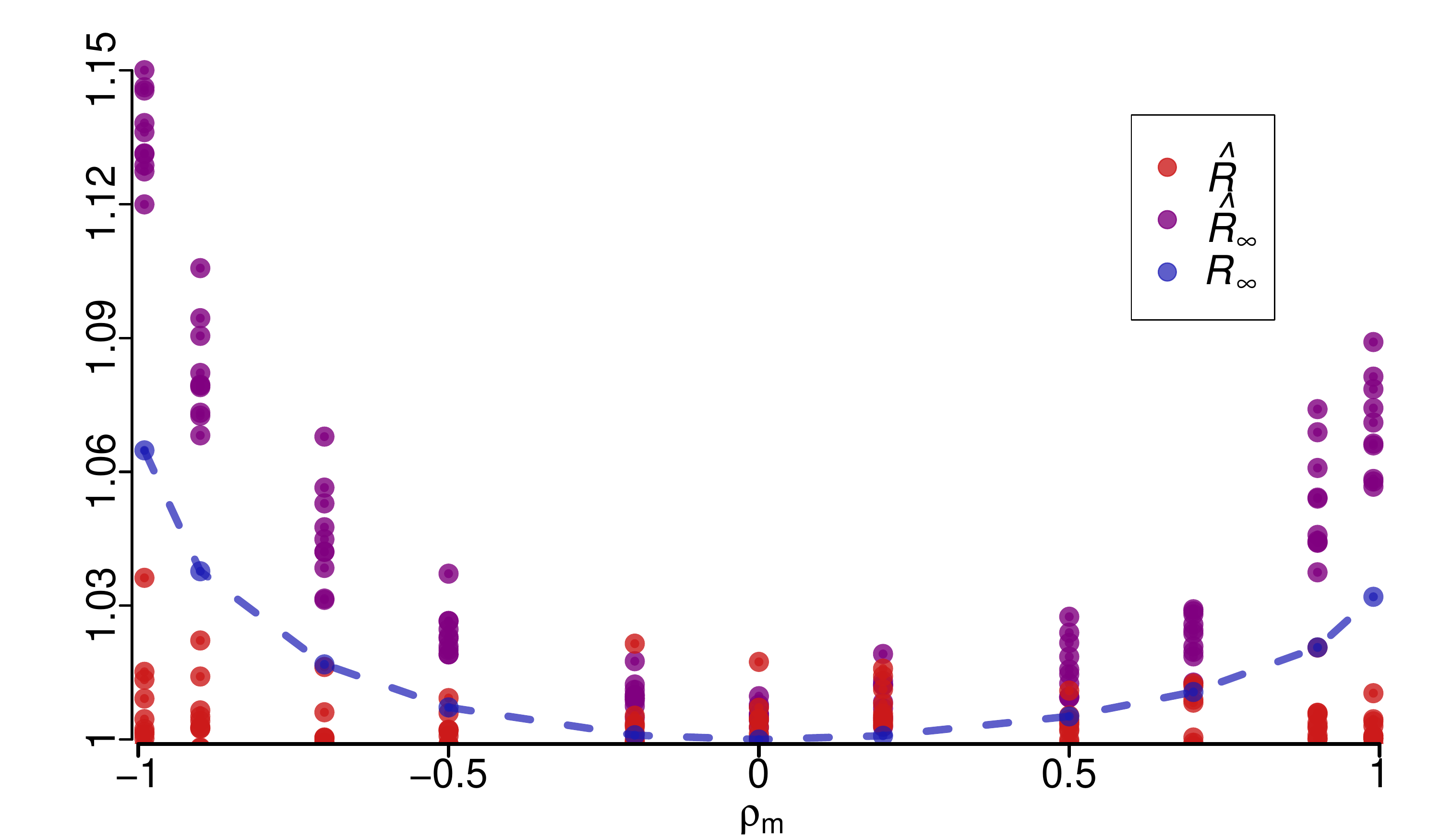}\\
            \vspace{-.3cm}
            \hspace{-1.5cm}\footnotesize{\color{orange2}$\hat{R}$} \hspace{2cm} \footnotesize{\color{violet2}$\hat{R}_\infty$} & 
        \end{tabular}
        \caption{Behaviour of Brooks--Gelman {\color{orange2}$\hat{R}$} (in orange) and multivariate {\color{violet2}$\hat{R}_\infty$} (in violet) in the case of chains with bivariate normal distributions, with different off-diagonal elements in the covariance matrix.
        On the left: Histograms with $100$ replications with one standard normal chain and one with $\rho_m = 0.9$.
        On the right: The same experiment with $10$ replications for different values of $\rho_m$, plotted as a function of $\rho_m$, and the corresponding population {\color{blue}$R_\infty$} in blue.
        }
        \label{fig:bivariate}
    \end{figure}

    \paragraph{Example~\hyperref[fig:bivariate]{4}: Bivariate normal distributions with different correlation terms.}
    
    In the bivariate case, the dependence structure is driven by only one value, which is the off-diagonal element $\rho_j \in (-1,1)$ of  $\boldsymbol{\Sigma}_j$. 
    Similarly to other examples, we suppose that we have $m-1$ converging chains with identity covariance matrix ($\rho_1 = \cdots =  \rho_{m-1} = 0$) while $\rho_m \in (-1,1)$ for the last one.
    
    Results are shown in Figure~\ref{fig:bivariate}, with a comparison of $\hat{R}_\infty$ with the multivariate $\hat{R}$ of \cite{brooks1998general}. The histogram on the left represents the values of the two diagnostics for $100$ replications with $m=2$, $n=200$ and $\rho_m = 0.9$.
    Despite a large difference on the covariance term between the chains, we can see that Brooks--Gelman $\hat{R}$ fails to correctly diagnose this difference, as most of the values are between $1$ and $1.01$, contrary to $\hat{R}_\infty$.
    Due to the i.i.d nature of the example, the recent \gr{proposal} of \cite{vats2018revisiting} for a multivariate $\hat{R}$ does not detect any convergence issue \gr{\cancel{either}} as the diagnostic is not based on a comparison between chains.
    This difference of behaviour is confirmed on the right panel of Figure~\ref{fig:bivariate}, which illustrates 10 replications of both diagnostics as a function of $\rho_m$.
    For instance, if $\rho_m = 0$ then the four chains are identically distributed and no convergence issue should be raised.
    Conversely, the value of $\hat{R}$ should increase when $|\rho_m| \to 1$, as the difference between the last chain and the other ones increases.
    For the Brooks--Gelman version, we can see that the value of $\hat{R}$ is almost constant and thus insensitive to $\rho_m$, which is not satisfactory, contrary to $\hat{R}_\infty$ which has a parabolic shape.
    
    As discussed in Section~\ref{sub-direction}, the behaviour of $\hat{R}_\infty$ is not symmetric when $\rho_m \to -1$ and $\rho_m \to 1$: the upper bound corresponding to positive dependence diverges with the dimension (Corollary~\ref{cor:PLOD_bound} for PLOD copulas) whereas the one for negative dependence is bounded by approximately $1.14$ (Corollary~\ref{cor:NLOD_bound} for NLOD copulas).
    This leads to the intuition that the convergence diagnostic is more sensitive in the PLOD case than in the NLOD, but this observation is asymptotic and when $d=2$, the two bounds are respectively equal to $1.08$ and $1.04$, so the statement is reversed.
    This asymmetry is illustrated in Figure~\ref{fig:bivariate} on theoretical $R_\infty$ (in blue) and  estimations $\hat{R}_\infty$ (in purple).

\begin{figure}
        \begin{center}
        \textbf{Example~\hyperref[fig:hist_dim]{5}:} Multivariate normal with $d \in \{2, \ldots, 6\}$\\
        \includegraphics[width = \textwidth]{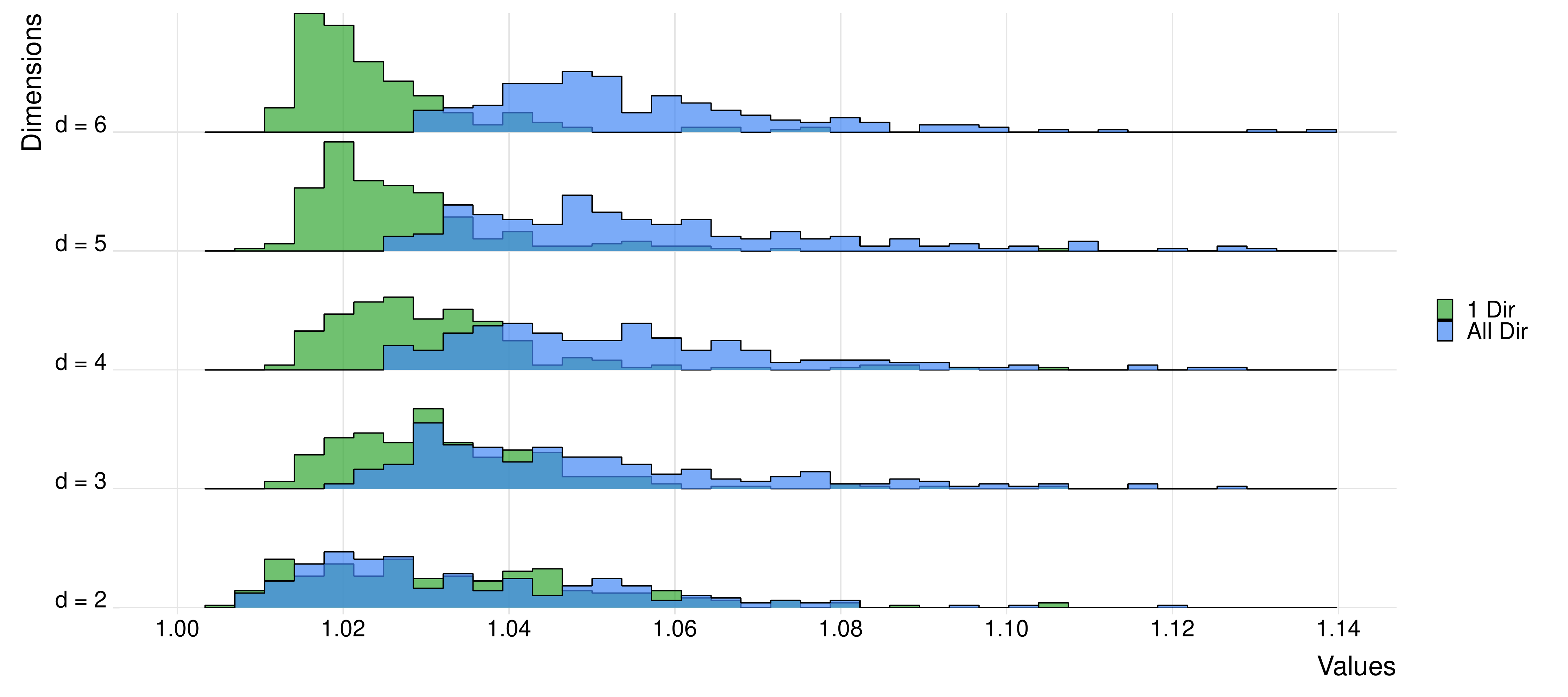}\\
        \vspace{-.8cm}
        \end{center}
        \hspace{-6cm}
        \footnotesize{{\color{forestgreen}$\hat{R}_\infty$}} \hspace{.3cm} \footnotesize{\color{blue}$\hat{R}^{(\text{max})}_\infty$}
        \caption{Comparison between {\color{forestgreen}$\hat{R}_\infty$} computed on one direction (in {\color{forestgreen}green}), and {\color{blue}$\hat{R}^{(\text{max})}_\infty$}, the maximum of $\hat{R}_\infty$ computed on all possible indicator functions (in {\color{blue}blue}). For each $d\in \{2, \ldots, 6\}$, 200 replications are done where a new covariance matrix is generated for the normal distribution, which leads to different directions of dependence among the replications.
        }
        \label{fig:hist_dim}
    \end{figure}

    \paragraph{Example~\hyperref[fig:hist_dim]{5}: Evolution of the behaviour when the dimension increases.}
    
    In the general case of dimensionality $d > 2$, we still compare $m-1$ chains that follow a multivariate standard normal distribution with one that has a given covariance matrix $\boldsymbol{\Sigma}_m$.
    To obtain $\boldsymbol{\Sigma}_m$, we generate a matrix $\boldsymbol{S}$ according to Wishart distribution with $d$ degrees of freedom, and we transform $\boldsymbol{S}$ in order to have one on the diagonal to keep the same margins for all chains (while remaining semi definite positive):
    \begin{equation*}
        \boldsymbol{\Sigma}_m = \boldsymbol{D}^{-1/2} \boldsymbol{S} \boldsymbol{D}^{-1/2}, 
        \quad \text{ with } \quad 
        \boldsymbol{D} = \text{diag}\left(s_{1,1}, \ldots, s_{d,d}\right).
    \end{equation*}
    To illustrate the influence of the dependence direction (Section~\ref{sub-direction}), a new matrix $\boldsymbol{\Sigma}_m$ is generated for each simulation, in order to have varying directions across replications.
    Then, we compare $\hat{R}_\infty$ with $\hat{R}^{(\text{max})}_\infty$, the maximum of $\hat{R}_\infty$ over all $2^{d-1}$ possible directions for the indicator functions.
    
    Results are shown in Figure~\ref{fig:hist_dim}, where $200$ replications are shown for $\hat{R}_\infty$ and $\hat{R}^{(\text{max})}_\infty$ for $d \in \{2, \ldots, 6\}$.
    As $\hat{R}^{(\text{max})}_\infty$ requires the computation of $2^{d-1}$ different $\hat{R}_\infty$, obtaining these histograms quickly becomes infeasible for larger dimensions.
    When $d=2$, we can see that there is no significant difference between $\hat{R}_\infty$ and $\hat{R}^{(\text{max})}_\infty$, but as the dimension increases the values of $\hat{R}_\infty$ become more concentrated and closer to one.
    Indeed, as the number of possible directions increases exponentially, it is more and more rare to obtain the one to which $\hat{R}_\infty$ is sensitive. 
    On the contrary, $\hat{R}^{(\text{max})}_\infty$ seems to stay robust with respect to this curse of dimensionality in terms of sensitivity, as the histograms look invariant when $d$ increases.

\section{Empirical results}
\label{sec:simulations}

    In Section~\ref{subsec:toy_examples} and Section~\ref{subsec:multivariate_toy_examples}, we considered toy examples where the distribution of the chains is known in order to control the value of the population $R_\infty$ and illustrate the robustness when other versions of $\hat{R}$ fail.
    Here we extend to other models in a more practical case for Bayesian inference. 
    We adopt a baseline similar to the one used by \cite{Vehtari} and \cite{lambert2021r} to illustrate the behaviour of $\hat{R}_\infty$ on Bayesian models, and add a multivariate example studied in \cite{vats2019multi}.
    For all examples in this section, we choose $4$ chains \gr{and therefore a threshold $R_{\infty,\text{lim}} = 1.02$ in the univariate case (according to Section~\ref{subsec:convergence-rhat}), and $(R_{\infty,\text{lim}}^{(C)}, R_{\infty,\text{lim}}^{(M)})= (1.03,1.03)$ in the multivariate one (according to Section~\ref{subsec:multivariate_def}).
    For each univariate study}, we plot an example of $\hat{R}(x)$ as a function of $x$, and we recommend this illustration to users who want to analyse more carefully a given value of $\hat{R}_\infty$.
    Together with this figure, we also show histograms of replications to check the behaviour of the different $\hat{R}$ more rigorously. 
    All experiments are done on R using \texttt{rstan} library \citep{rstan} and the package \texttt{localrhat} that we propose with this paper \citep{moinsLocalrhat}. 
    Additional experiments have also been conducted on Python using OpenTURNS \citep{baudin2015open}.
    All the code concerning these experiments and the additional ones are available in the online appendix (link in the Introduction).

\begin{figure}
        \begin{tabular}{rl}
            \textbf{Example \hyperref[fig:autoreg]{6}:} & Autoregressive model \\
            \vspace{-.3cm}
        \includegraphics[trim=0cm 0cm 0cm 0cm, clip, width = .49\textwidth]{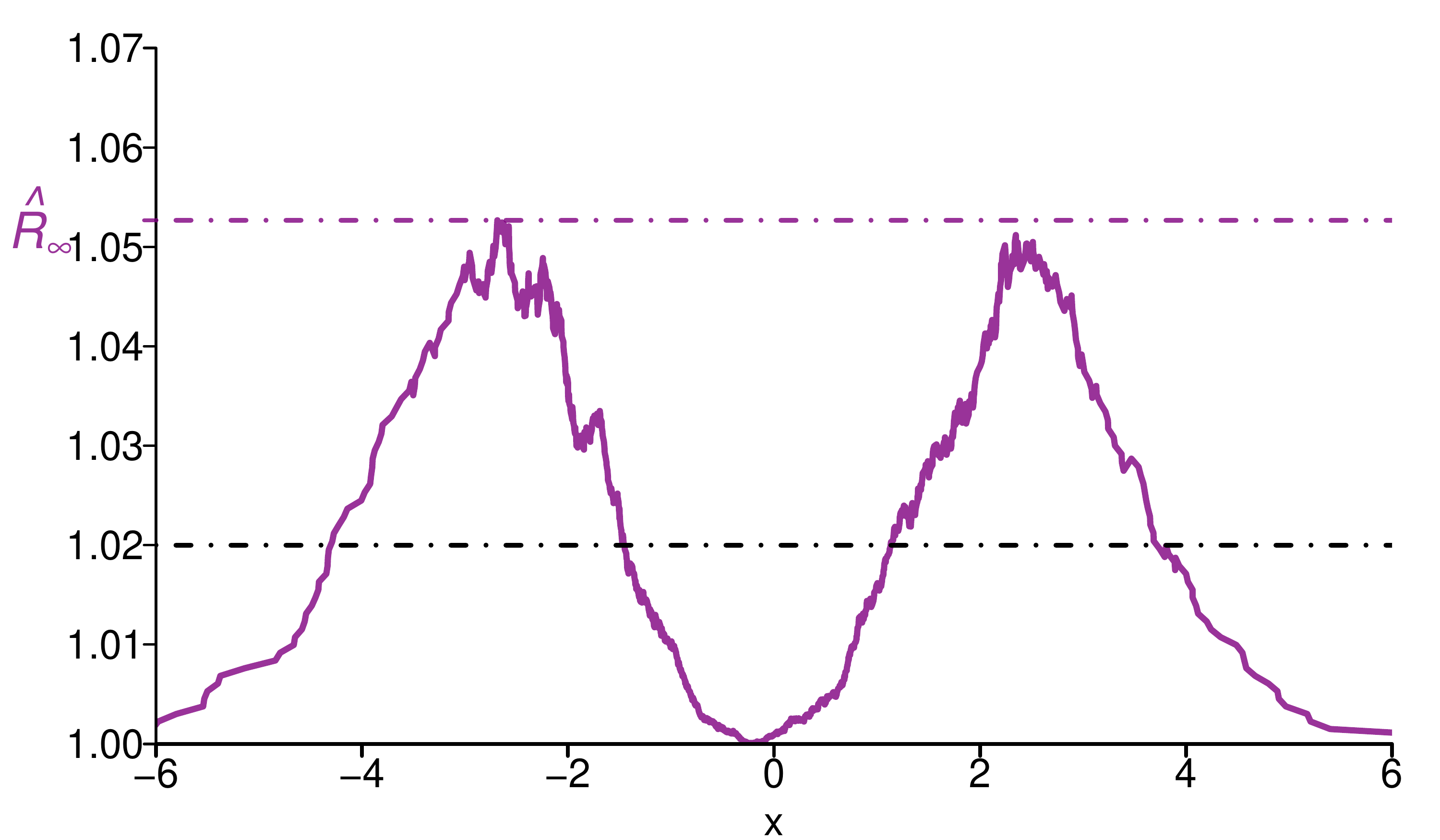} & 
            \includegraphics[trim=0cm 0.5cm 0cm 0cm, clip, width = .49\textwidth]{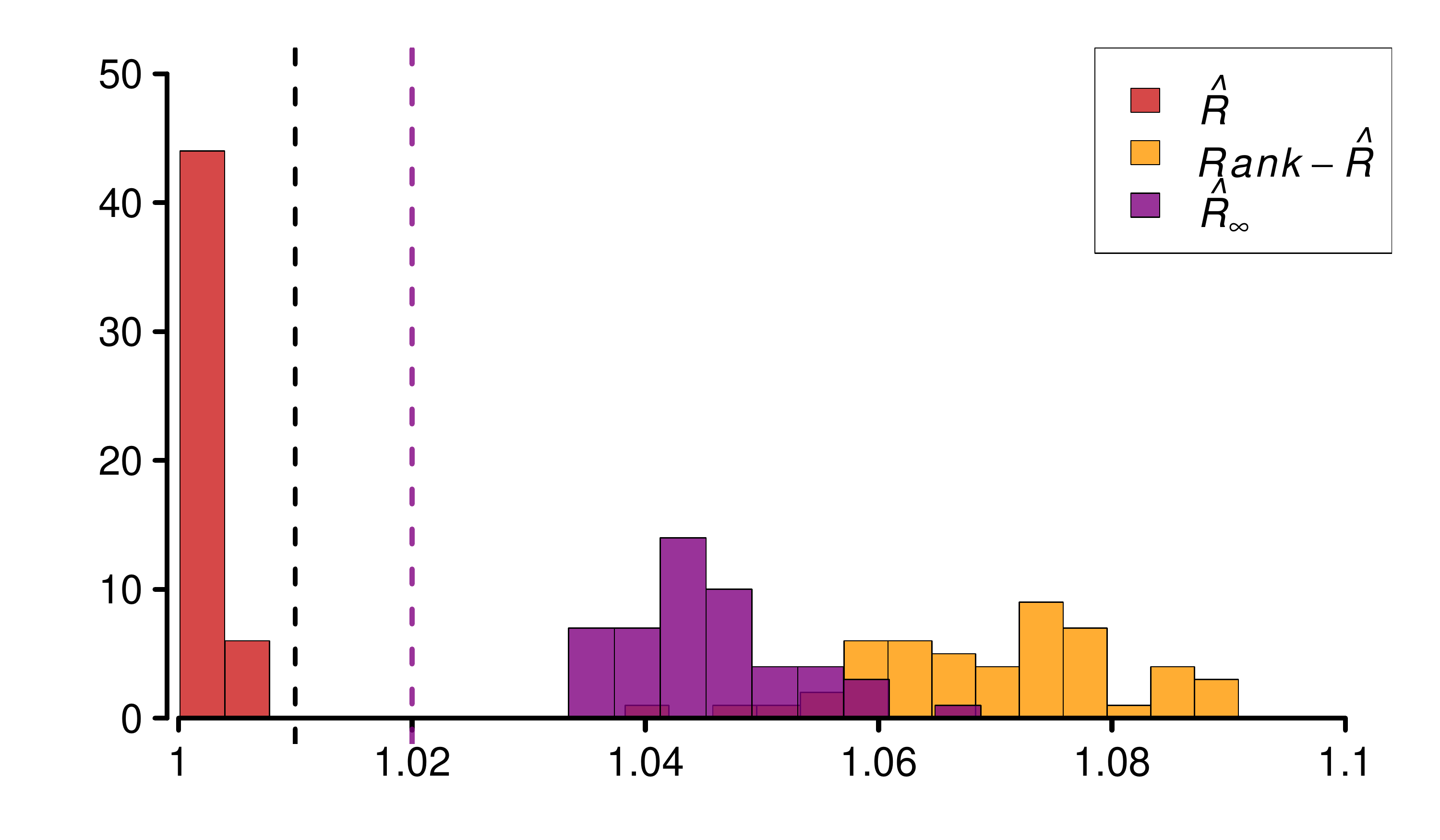}\\
        \vspace{-.3cm}
        & \hspace{1cm}\footnotesize{\color{orange2}$\hat{R}$}  \hspace{1.9cm} \footnotesize{\color{violet2}$\hat{R}_\infty$} \hspace{.7cm} \footnotesize{\color{red2}rank-$\hat{R}$}
        \end{tabular}
        \caption{Behaviour of $\hat{R}_\infty$ on the autoregressive example described in Section~\ref{subsec:autoreg}, with $m=4$ chains of size $n=500$ and $(\sigma, \sigma_m, \rho) = (1, 2, 1/2)$. 
        On the left: $\hat{R}(x)$ as a function of $x$ for one replication.
        On the right: Histograms of $50$ replications of {\color{orange2}$\hat{R}$}, {\color{red2}rank-$\hat{R}$} and {\color{violet2}$\hat{R}_\infty$}. The dashed lines correspond to thresholds of $1.01$ and $1.02$.
        }
        \label{fig:autoreg}
    \end{figure}

    \paragraph{Example \hyperref[fig:autoreg]{6}: Autoregressive model with different variances.}
    \label{subsec:autoreg}

    The first example is a basic autoregressive model to study the behavior of $\hat{R}_\infty$ in the case of Markov chains with different variances: 
    we consider $m$ chains of size $n$ such that for $i\in \{1,\ldots, n-1\}$ and $j\in \{1, \ldots, m\}$,
    \begin{equation*}
        \theta^{(i+1,j)} = \rho \theta^{(i,j)} + \epsilon_{i,j}, 
        \quad \text{with} \quad
        \epsilon_{i,j} \sim \mathcal{N}(0,\sigma^2_j),
    \end{equation*}
    where $\rho \in (0,1)$ and $\sigma_j > 0$.
    In particular, assume that the first $m-1$ chains are generated using the same process: $\sigma_1 = \cdots = \sigma_{m-1} = \sigma$, while for the last chain $\sigma_m \neq \sigma$.
    
    Results are illustrated in Figure~\ref{fig:autoreg} with $m=4$, $\sigma = 1$, $\sigma_m = 2$ and $\rho = 1/2$ on $50$ replications, and an example of $\hat{R}(x)$ as a function of $x$ on the left panel.
    Similarly to the rank-$\hat{R}$ replications, the $\hat{R}_\infty$ values remain far from the threshold of $1.02$ which confirms the sensitivity to this convergence defect. 
    This corroborates in a more practical case the results of Example~\hyperref[fig:example_known_dist]{1} in Section~\ref{subsec:toy_examples}, on the sensitivity of $\hat{R}_\infty$ on chains with same mean and different variances.
    Note that the value $R(0) = 1$ is due to the fact that all the chains share the same median equal to zero.

\begin{figure}[h!]
        \begin{tabular}{rl}
            \textbf{Example~\hyperref[fig:cauchy]{7}.a:} & HMC on nominal Cauchy\\
            \vspace{-.3cm}
            \includegraphics[trim=0cm 0cm 0cm 0cm, clip, width = .49\textwidth]{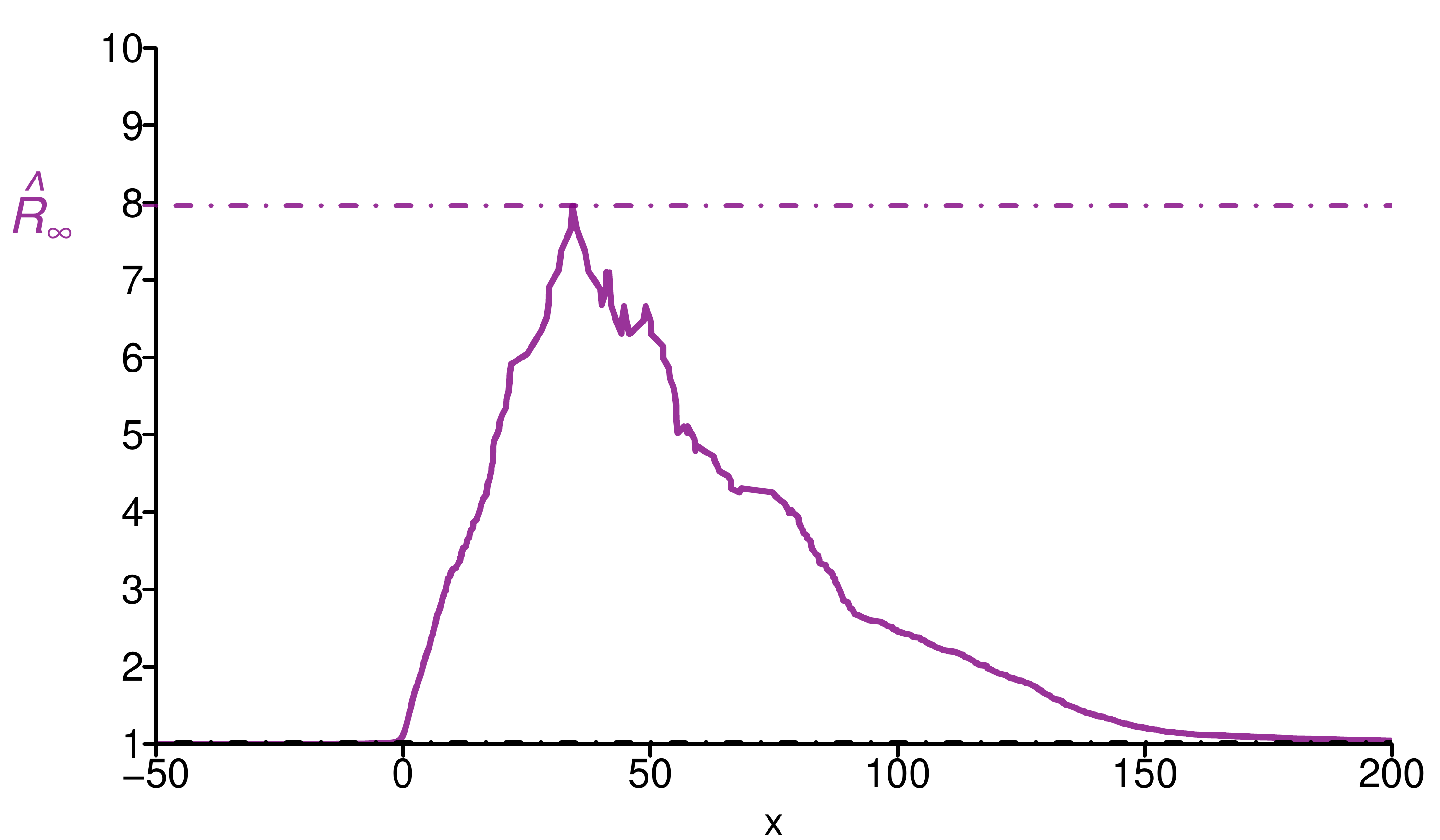} & 
            \includegraphics[trim=0cm 0.5cm 0cm 0cm, clip, width = .49\textwidth]{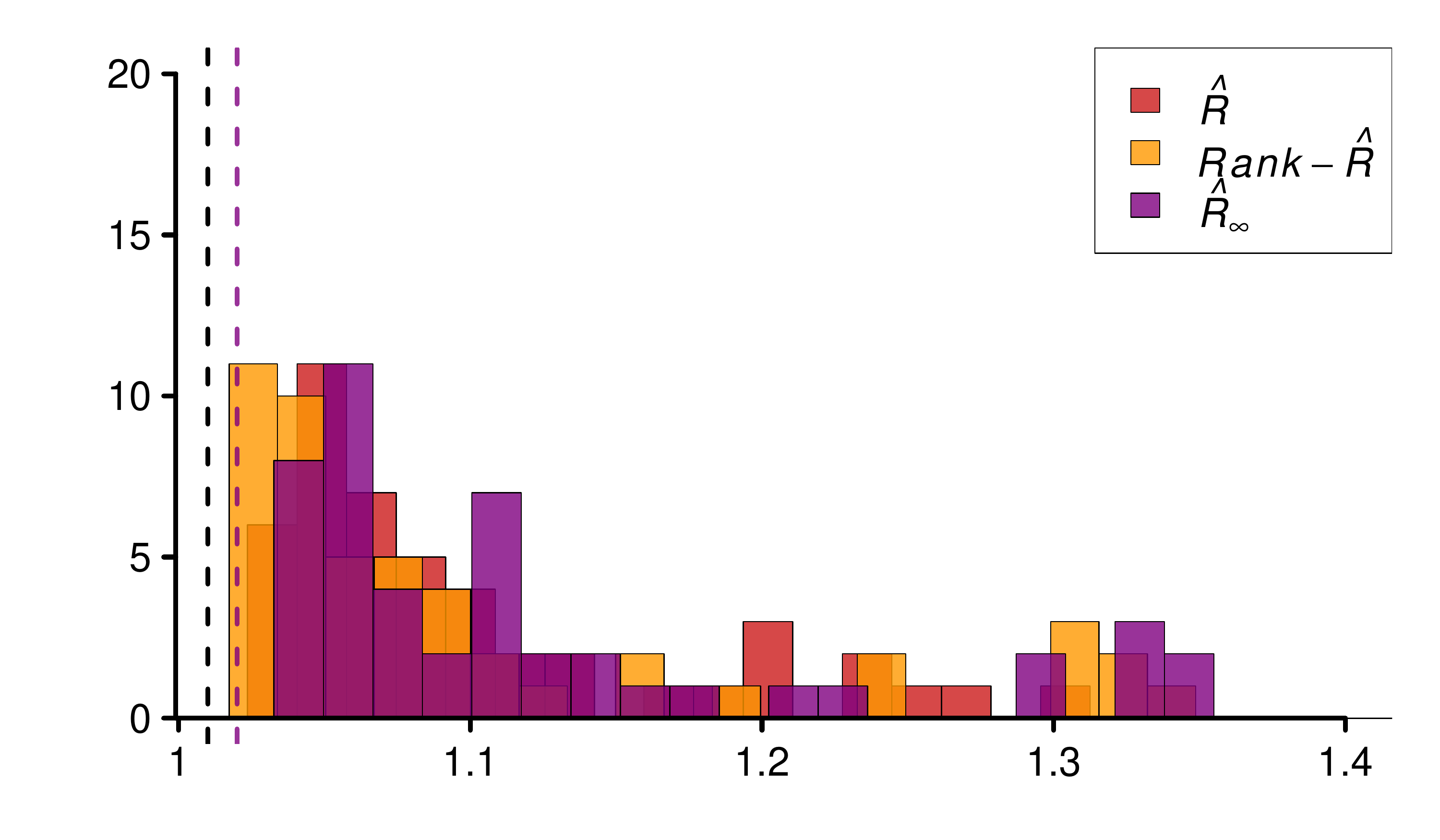}\\
            \vspace{0.1cm}
            & \hspace{0.4cm}\footnotesize{\color{red2}rank-$\hat{R}$}  
            \hspace{0cm} 
            \footnotesize{\color{orange2}$\hat{R}$} 
            \hspace{0cm} 
            \footnotesize{\color{violet2}$\hat{R}_\infty$}\\
            \textbf{Example~\hyperref[fig:cauchy]{7}.b:} & HMC on alternative Cauchy \\
            \vspace{-.3cm}
            \includegraphics[trim=0cm 0cm 0cm 0cm, clip, width = .49\textwidth]{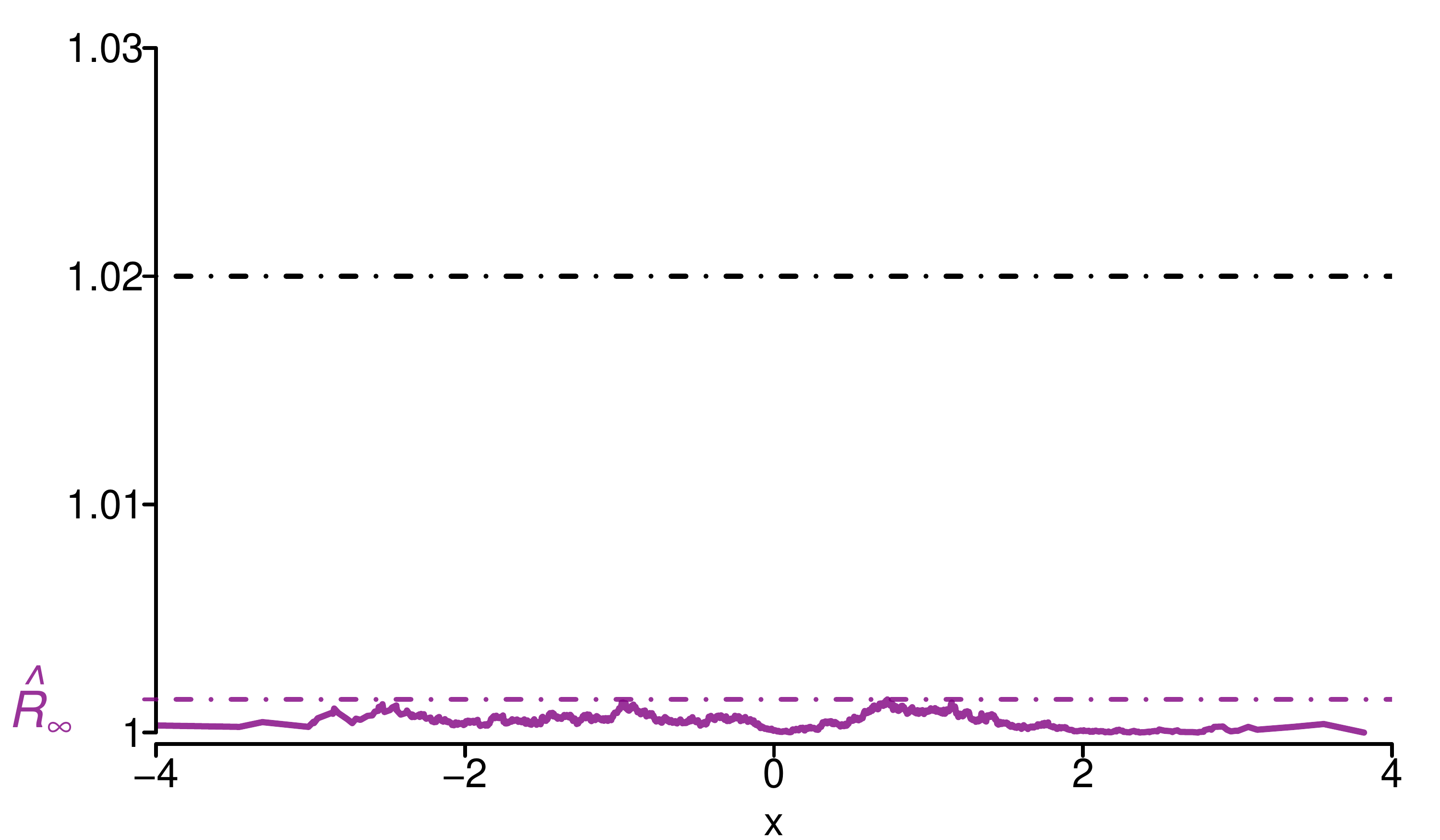} & 
            \includegraphics[trim=0cm 0.5cm 0cm 0cm, clip, width = .49\textwidth]{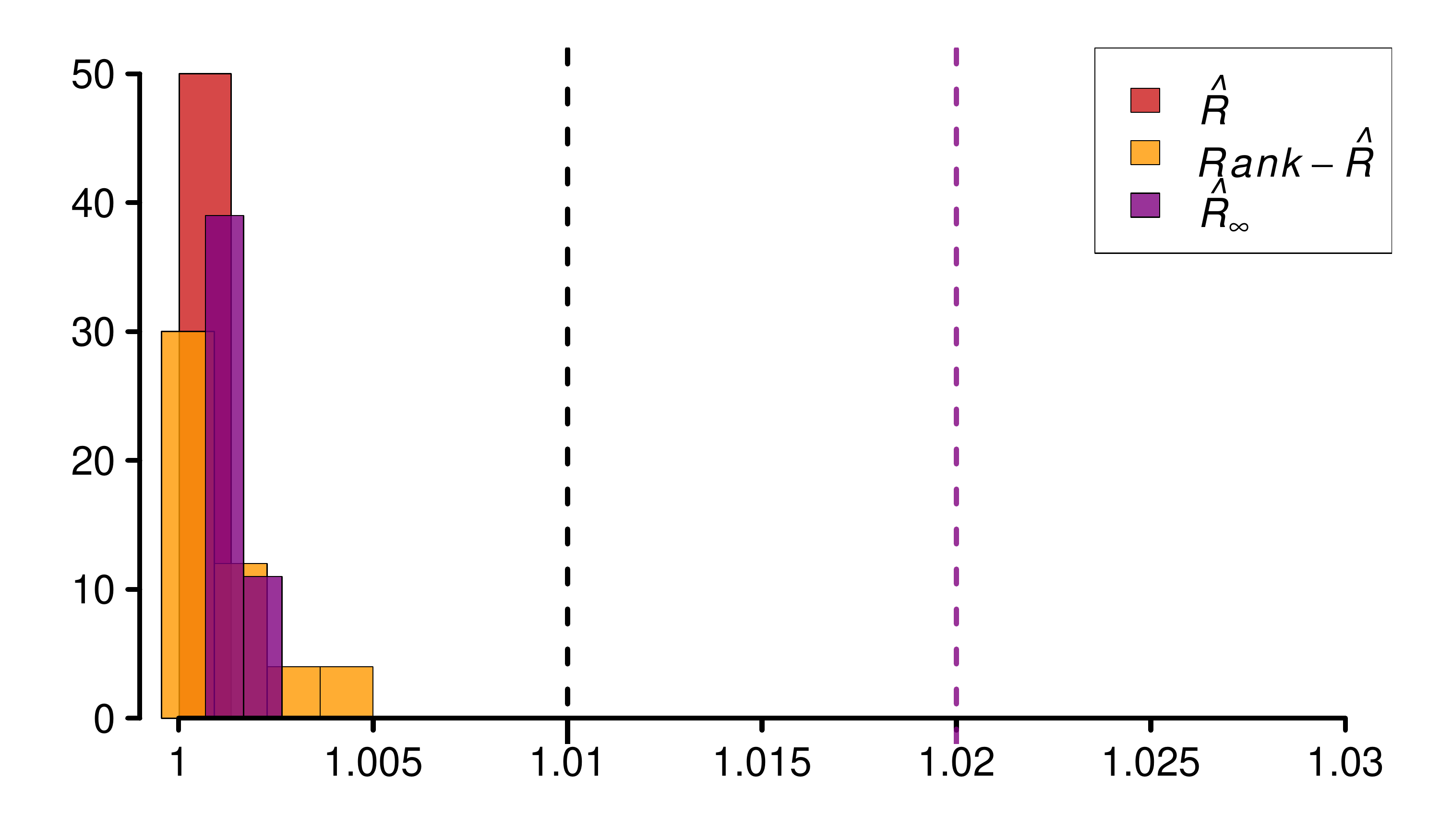}\\
        \vspace{-.3cm}
        & \hspace{-0.1cm}\footnotesize{\color{red2}rank-$\hat{R}$}  
        \hspace{0cm} 
        \footnotesize{\color{orange2}$\hat{R}$} 
        \hspace{0cm} 
        \footnotesize{\color{violet2}$\hat{R}_\infty$}\\
        \end{tabular}
        \caption{Behaviour of $\hat{R}_\infty$ on the Cauchy example described in Section~\ref{subsec:cauchy} for the two parameterisations. 
        On the left: $\hat{R}(x)$ as a function of $x$ for one replication.
        On the right: Histograms of $50$ replications of {\color{orange2}$\hat{R}$}, {\color{red2}rank-$\hat{R}$} and {\color{violet2}$\hat{R}_\infty$}. The dashed lines correspond to thresholds of $1.01$ and $1.02$.
        }
        \label{fig:cauchy}
    \end{figure}

    \paragraph{Example~\hyperref[fig:cauchy]{7}: HMC on Cauchy distribution.}
    \label{subsec:cauchy}
    As an extension of Example~\hyperref[fig:example_known_dist]{2} in Section~\ref{subsec:toy_examples}, we analyze the behaviour of $\hat{R}_\infty$ in the case of heavy-tailed distributions.
    We run Hamiltonian Monte Carlo (HMC) \citep{neal2011mcmc} using Stan on Cauchy distributions for $50$ variables. We consider the one with the most important mixing issue diagnosed with $\hat{R}_\infty$.
    Due to the tail heaviness of Cauchy distributions, the HMC iterations on a given chain can \gr{get trapped} in a tail, which causes mixing issues.
    One solution to avoid this is to use an alternative parameterisation that avoids sampling from a heavy-tailed distribution:
    
    \begin{tabular}{ll}
        \textbf{Example~\hyperref[fig:cauchy]{7}.a:} & Nominal parameterisation \\[0.1cm]
        & $x_j \sim \text{Cauchy}(0,1), \quad j\in \{1, \ldots, 50\}$.\\[0.2cm]
        \textbf{Example~\hyperref[fig:cauchy]{7}.b:} & Alternative parameterisation \\[0.1cm]
         & $x_j = a_j / \sqrt{b_j}$,\quad $a_j \sim \mathcal{N}(0,1),\quad b_j \sim \chi^2_1$.\\[0.2cm]
    \end{tabular}

    \noindent One would expect convergence issues with the nominal parameterisation and not with the alternative one.
    For both, the process of selecting the worst parameters among the $50$ ones is iterated for the generation of replications, and results are shown in Figure~\ref{fig:cauchy}.
    Histograms on the top right confirm the risk of diverging chains with the nominal parameterisation, as all the values are above $1.02$ for all the versions of $\hat{R}$.
    This means that it is very likely to have at least one chain out of the $50$ with a convergence issue in this experiment.
    This divergence can be really extreme, as it is shown on the top left panel where the value of $\hat{R}_\infty$ is over seven, due to a mixing issue in the right tail of the distribution. 
    The opposite occurs with the other parameterisation, as all the convergence diagnostics indicate no mixing issues (see bottom row of Figure~\ref{fig:cauchy}), which means no counter-indications that the chains for the 50 variables have converged.
    Looking at $\hat{R}(x)$ function on one replication in the bottom left panel, the curve seems to be very noisy and close to 1 compared to $1.02$ (even sometimes less than 1) so the difference with 1 seems only due to Monte Carlo noise.
    
\begin{figure}
        \begin{tabular}{rl}
            \textbf{Example~\hyperref[fig:hierarchical]{8}.a:} & Centered eight schools\\
            \vspace{-0.3cm}
            \includegraphics[trim=0cm 0cm 0cm 0cm, clip, width = .49\textwidth]{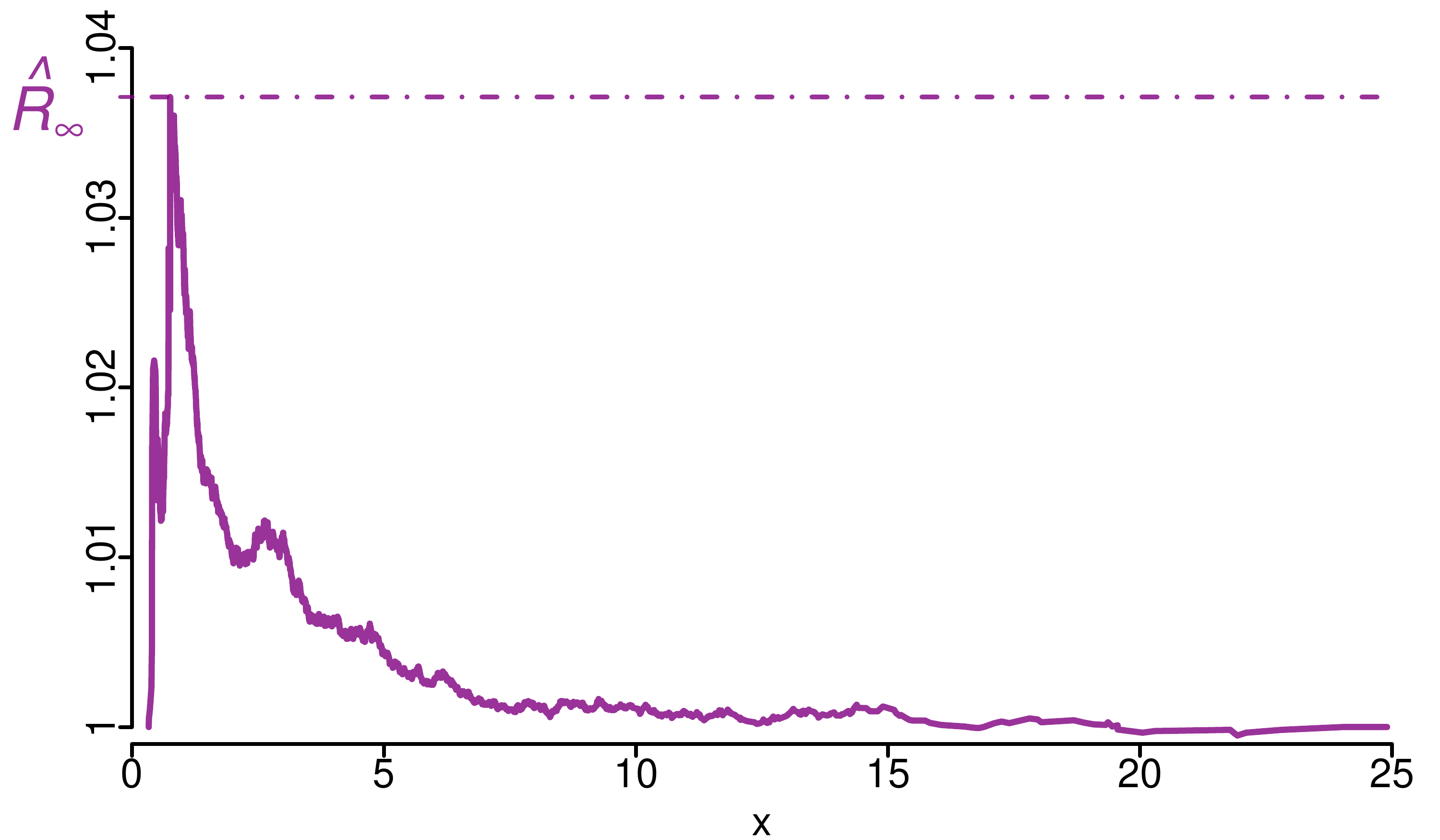} & 
            \includegraphics[trim=0cm 0.5cm 0cm 0cm, clip, width = .49\textwidth]{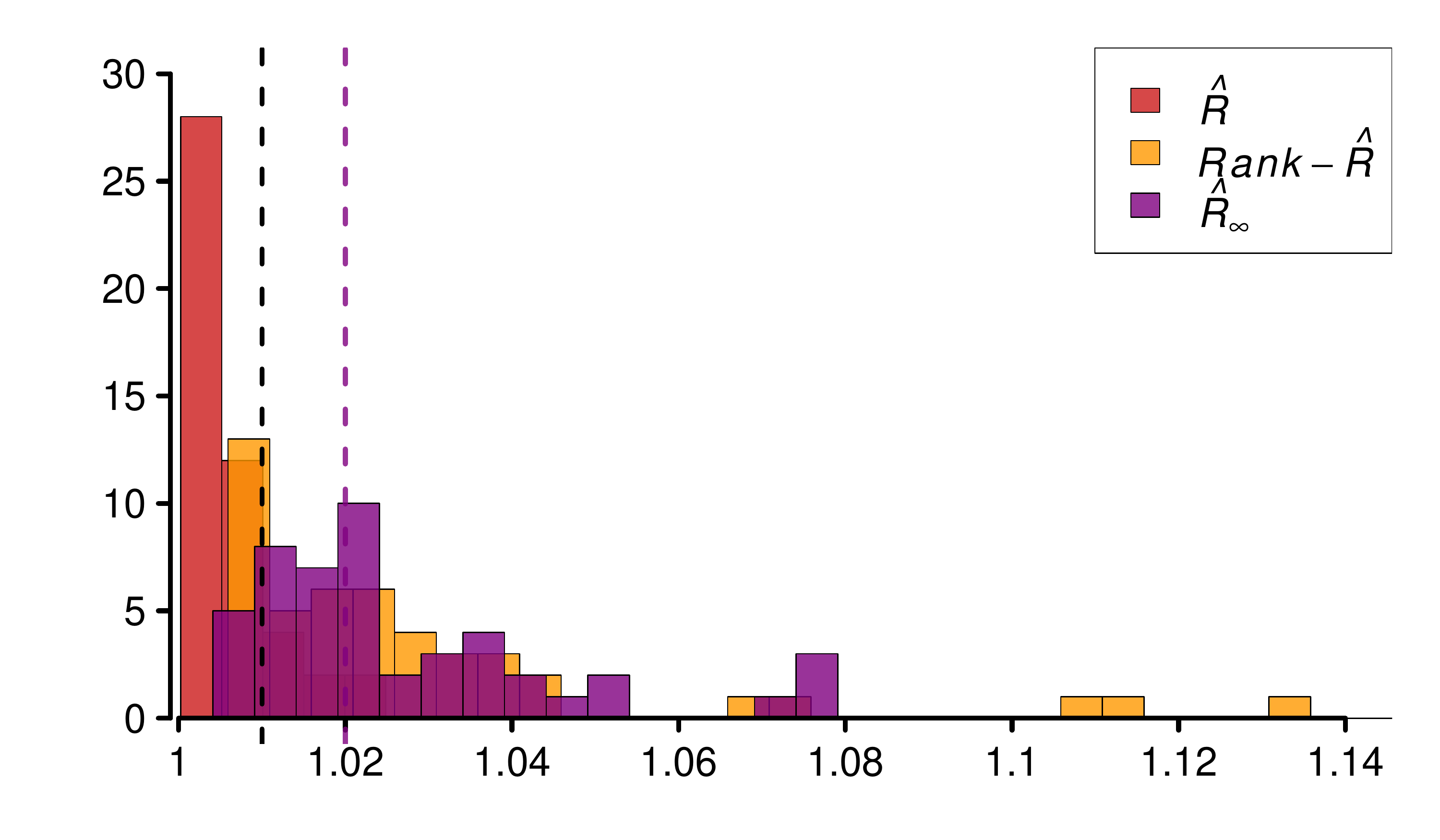}\\
            \vspace{0.1cm}
            & 
            \hspace{0.4cm}
            \footnotesize{\color{orange2}$\hat{R}$}  
            \hspace{-0.1cm} 
            \footnotesize{\color{red2}rank-$\hat{R}$} 
            \hspace{-0.1cm} 
            \footnotesize{\color{violet2}$\hat{R}_\infty$}\\
            \textbf{Example~\hyperref[fig:hierarchical]{8}.b:} & Non-centered eight schools\\
            \vspace{-0.3cm}
            \includegraphics[trim=0cm 0cm 0cm 0cm, clip, width = .49\textwidth]{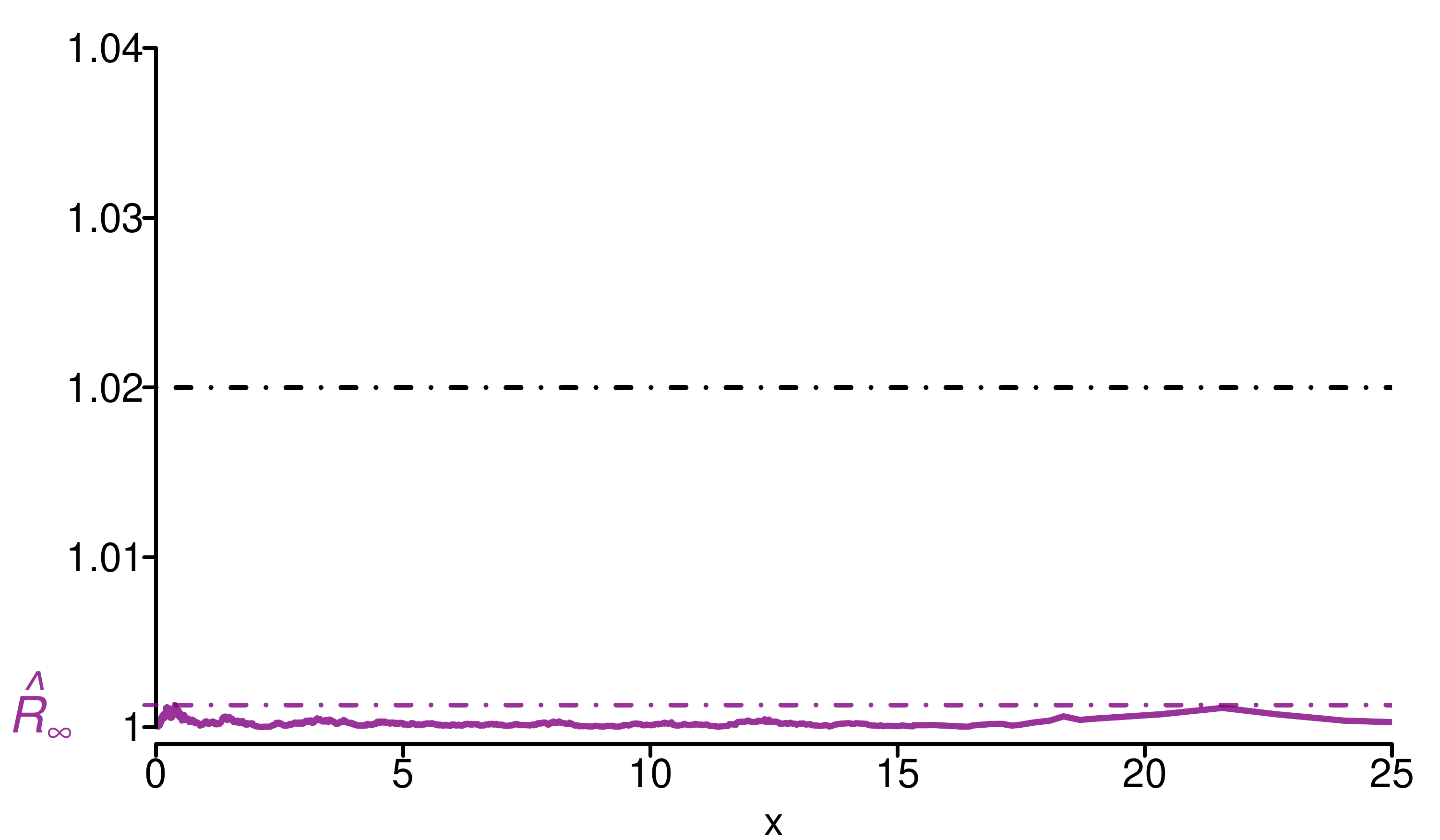} & 
            \includegraphics[trim=0cm 0.5cm 0cm 0cm, clip, width = .49\textwidth]{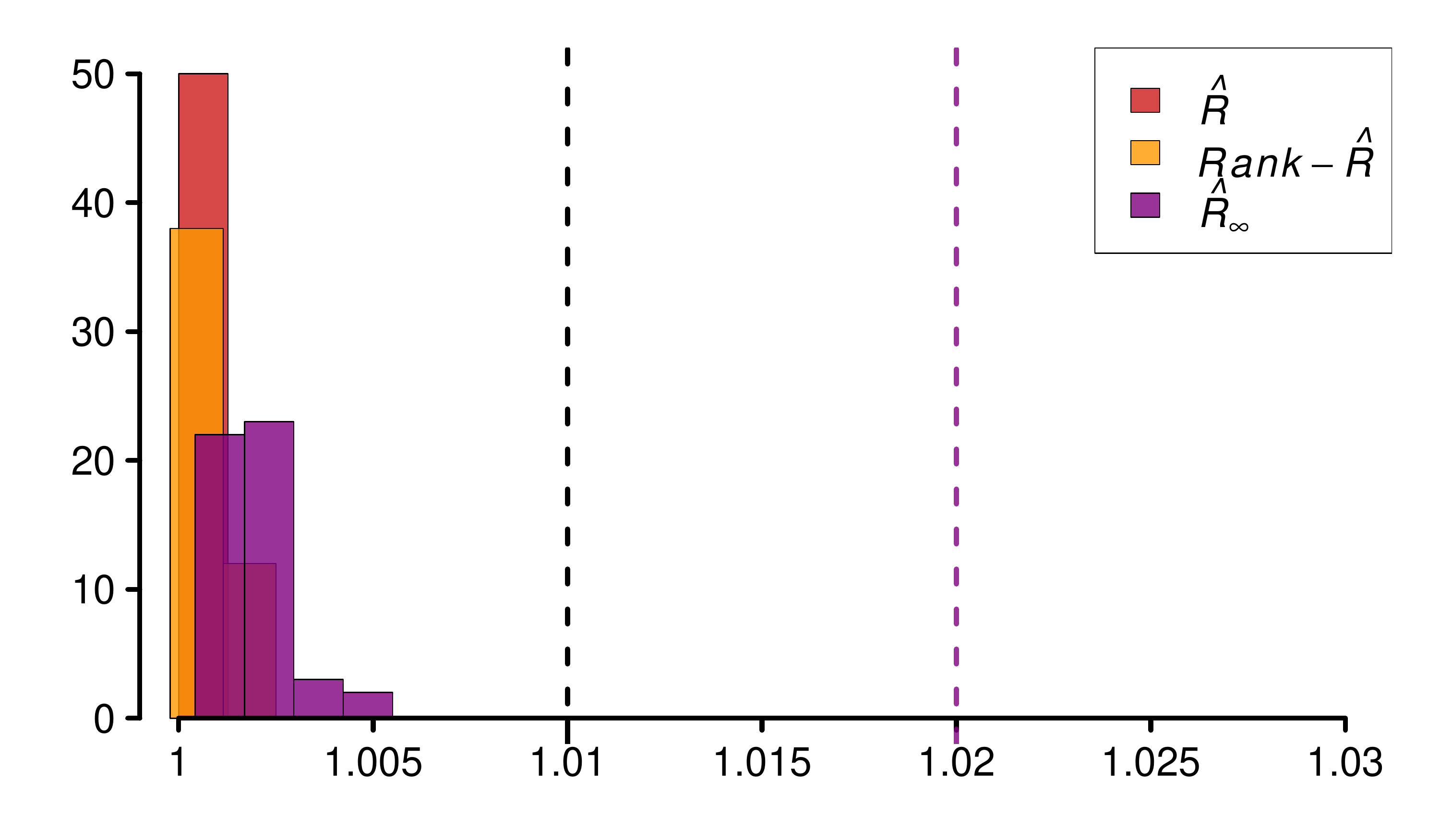}\\
            \vspace{-0.3cm}
            & 
            \hspace{-0.1cm}
            \footnotesize{\color{red2}rank-$\hat{R}$}  
            \hspace{-0.1cm} 
            \footnotesize{\color{orange2}$\hat{R}$} 
            \hspace{-0.1cm} 
            \footnotesize{\color{violet2}$\hat{R}_\infty$}
        \end{tabular}
        \caption{Behaviour of $\hat{R}_\infty$ on the hierarchical example for $\tau$ described in Section~\ref{subsec:hierarchical} for the centered and non-centered version. 
        On the left: $\hat{R}(x)$ as a function of $x$ for one replication.
        On the right: Histograms of $50$ replications of {\color{orange2}$\hat{R}$}, {\color{red2}rank-$\hat{R}$} and {\color{violet2}$\hat{R}_\infty$}. The dashed lines correspond to thresholds of $1.01$ and $1.02$.
        }
        \label{fig:hierarchical}
    \end{figure}

    \paragraph{Example~\hyperref[fig:hierarchical]{8}: Hierarchical Bayesian model on two parameterisations.}
    \label{subsec:hierarchical}
    
    As a classical Bayesian example, we consider using HMC on a hierarchical Bayesian model and in particular the eight-school \citep[Section~5.5]{gelman2013bayesian}, where two parameterisations are possible to model the problem:\\
    
    \begin{tabular}{ll}
        \textbf{Example~\hyperref[fig:hierarchical]{8}.a:} & Centered parameterisation (CP) \\[0.1cm]
        & $\theta_j \sim \mathcal{N}(\mu,\tau), \quad 
        y_j \sim \mathcal{N}(\theta_j,\sigma^2_j)$.\\[0.2cm]
        \textbf{Example~\hyperref[fig:hierarchical]{8}.b:} & Non-centered parameterisation (NCP) \\[0.1cm]
        & $\bar{\theta}_j \sim \mathcal{N}(0,1), \quad 
        \theta_j = \mu + \tau \bar{\theta}_j, \quad y_j \sim \mathcal{N}(\theta_j,\sigma^2_j)$.\\[0.2cm]
    \end{tabular}
    
    \noindent In the first parameterisation (CP), a prior dependence is between $(\mu, \tau)$ and the population parameters $\theta_j$, whereas in the other case (NCP), $\bar{\theta}_j$ is a priori independent of $(\mu, \tau)$, and $\theta_j$ is just a function of $\bar{\theta}_j$ and  $(\mu, \tau)$ \citep[see for example][for more details]{papaspiliopoulos2003}.
    \cite{Vehtari} argue in favour of the NCP for the eight-school example, by analysing the convergence of the chains associated with the parameter $\tau$.
    We also focus on computing $\hat{R}_\infty$ for $\tau$: results for $\hat{R}_\infty$ and comparison with other versions of $\hat{R}$ are shown in Figure~\ref{fig:hierarchical}.
    In the first row, we can see that the $\hat{R}_\infty$ diagnostic confirms the one of rank-$\hat{R}$, as the two corresponding histograms are similar in the top right panel and conclude for a lack of convergence in most of the cases. 
    However, for both diagnostics, a significant number of cases are also below $1.02$ (respectively $1.01$ for rank-$\hat{R}$), and most of them are close to this threshold, which is represented on the top left panel.
    In spite of this, the bottom row of Figure~\ref{fig:hierarchical} shows a clear difference and NCP seems to help for chain convergence.
    
\begin{figure}
        \begin{tabular}{rl}
            \textbf{Example~\hyperref[fig:logit]{9}:} & Bayesian logistic regression \\
            \vspace{-0.4cm}
            \includegraphics[trim=0cm 0cm 0cm 0cm, clip, width = .49\textwidth]{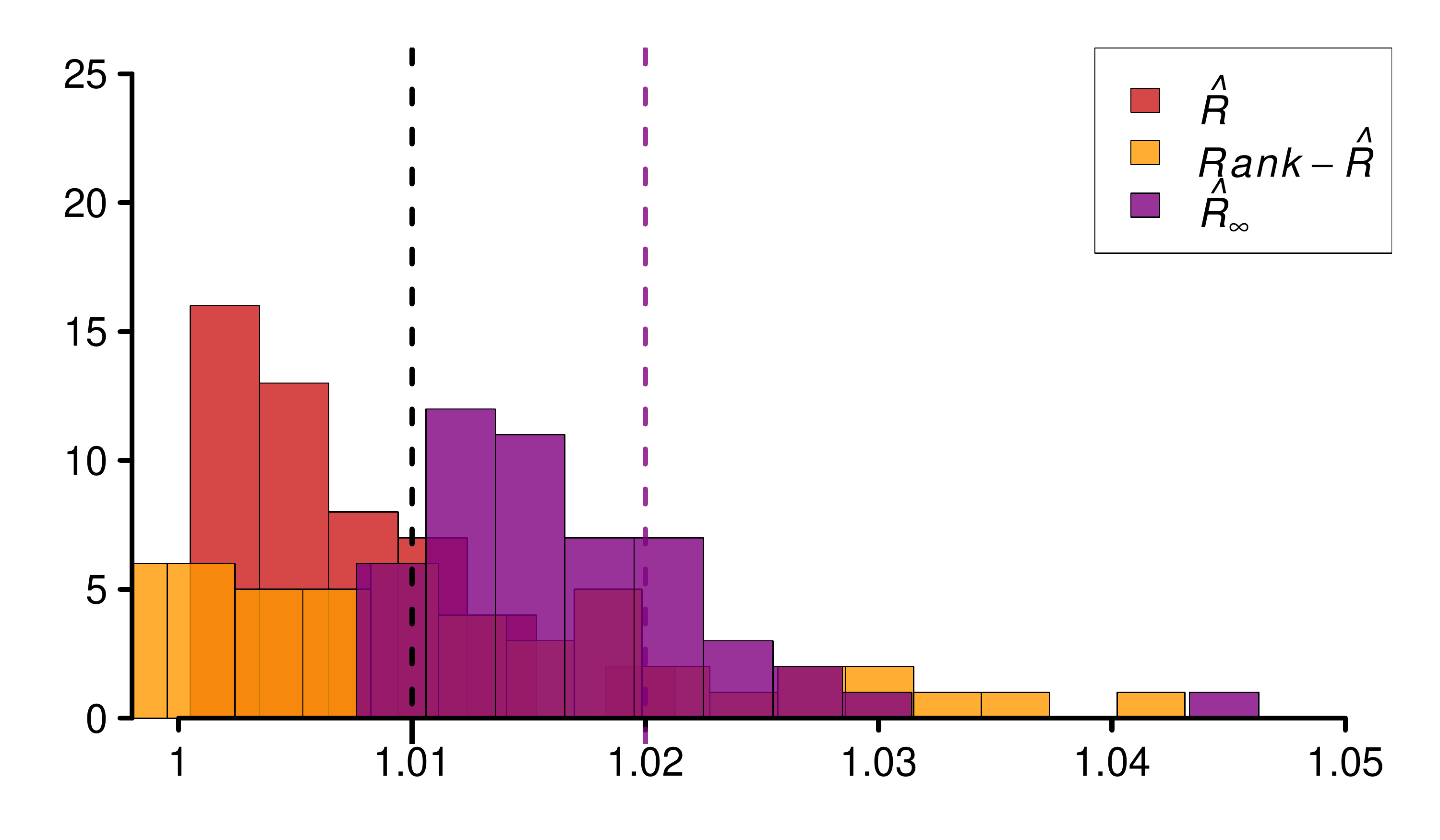} & 
            \includegraphics[trim=0cm 0cm 0cm 0cm, clip, width = .49\textwidth]{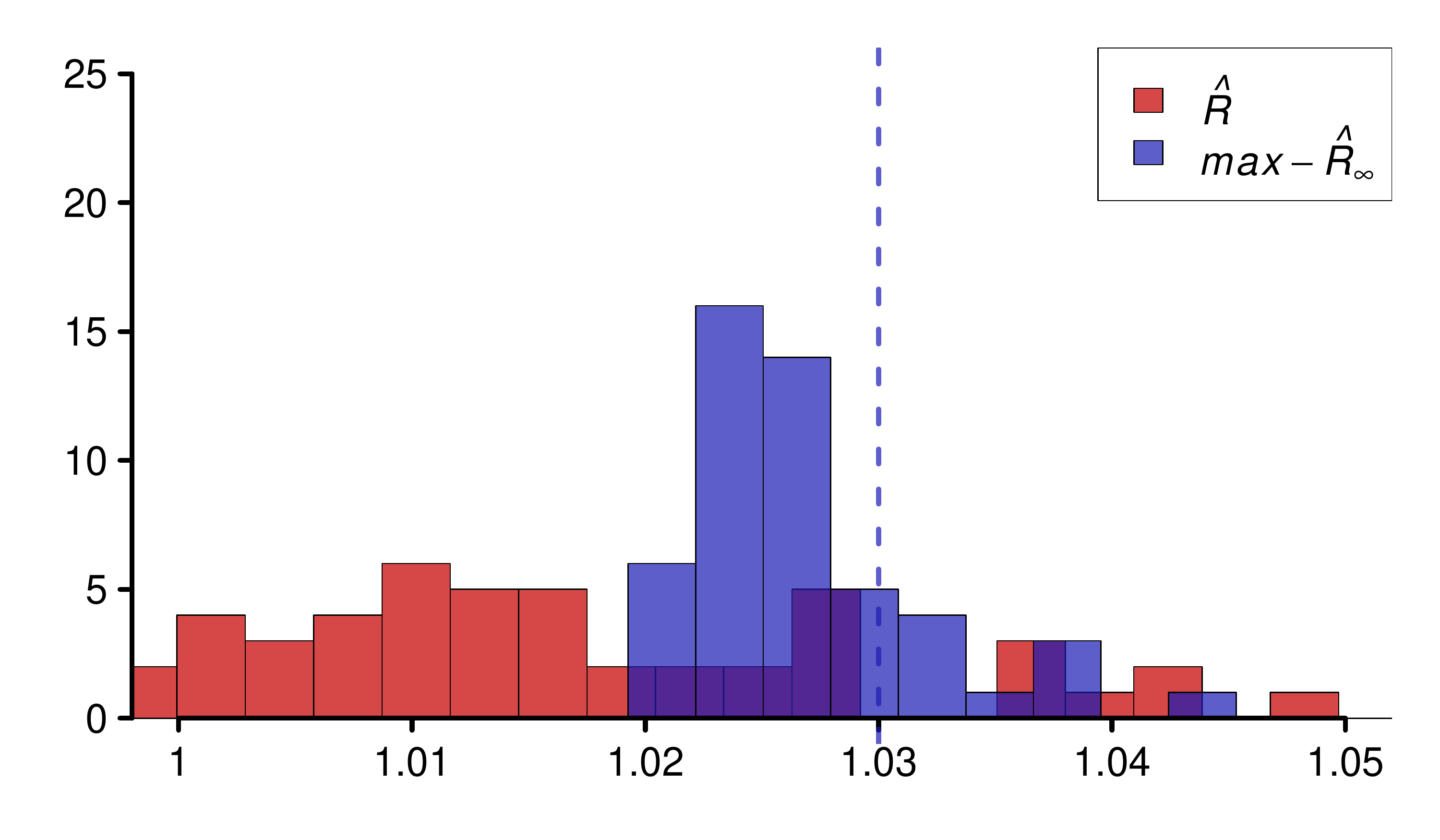}\\
            \vspace{-0.3cm}
            \footnotesize{\color{red2}rank-$\hat{R}$} \hspace{.1cm} 
            \footnotesize{\color{orange2}$\hat{R}$} \hspace{.7cm} 
            \footnotesize{\color{violet2}$\hat{R}_\infty$}
            \hspace{3.7cm} 
            &
            \hspace{1.6cm} 
            \footnotesize{\color{orange2}$\hat{R}$} \hspace{1cm} 
            \footnotesize{\color{blue}$\hat{R}^{(\text{max})}_\infty$}
        \end{tabular}
        \caption{Behaviour of multivariate and univariate $\hat{R}_\infty$ on the Bayesian logistic regression example, with $m=4$ chains of size $n=200$. 
        On the left: Histograms of $50$ replications of $\hat{R}$, rank-$\hat{R}$ and univariate $\hat{R}_\infty$ all applied on the log-posterior.
        On the right: Histograms of $50$ replications of Brooks--Gelman $\hat{R}$ and $\hat{R}^{(\text{max})}_\infty$. The dashed line corresponds to different thresholds: on the left, $1.01$ in black for $\hat{R}$ and rank-$\hat{R}$, $1.02$ in violet for $\hat{R}_\infty$, and on the right $1.03$ in blue for $\hat{R}^{(\text{max})}_\infty$.
        }
        \label{fig:logit}
    \end{figure}

    \paragraph{Example~\hyperref[fig:logit]{9}: Bayesian logistic regression.}
    This example is related to the extension of $\hat{R}_\infty$ in the multivariate case as proposed in Section~\ref{sec:multivariate}.
    As a multivariate Bayesian example, we run Stan on a basic hierarchical logistic model using the dataset \texttt{logit} available in the R package \texttt{mcmc}:
    \begin{equation*}
        \boldsymbol{\beta} \sim \mathcal{N}(0, 0.35^2 \boldsymbol{I}_4), \quad
        y_j \sim \text{Bernoulli}\left(\frac{1}{1+e^{-\boldsymbol{x}_j^\top\boldsymbol{\beta}}}\right).
    \end{equation*}
    Here the posterior is intractable and \cite{vats2019multi} showed that the posterior coefficients $\boldsymbol{\beta}$ could be significantly correlated, encouraging a multivariate diagnostic to check the convergence of the dependence structure.
    We run $m=4$ chains each of size $n=200$ after a burn-in of $100$.
    In this configuration, despite a low number of iterations, all the different univariate $\hat{R}_\infty$ are mostly below $1.02$ when replicated, and the rank-$\hat{R}$ are below $1.01$.

    When applied to the log posterior, the diagnostic is less clear and results are shown in the left panel of Figure~\ref{fig:logit}: a significant part of the histogram for $\hat{R}_\infty$ is below the threshold, meaning that the number of iterations is almost sufficient but is not yet.
    Looking at the right plot of Figure~\ref{fig:logit}, we notice in this example that the sensitivity of $\hat{R}^{(\text{max})}_\infty$ is approximately the same as the univariate version on the left, as the proportion of values over the threshold is similar
    \gr{(the choice of $R_{\infty,\lim}^{(C)} = 1.03$ is made according to Table~\ref{tab:rhat_max_tab})}.
    Although the computation of $\hat{R}^{(\text{max})}_\infty$ is possible here as the number of dimensions is small, computing a univariate $\hat{R}_\infty$ on the log posterior instead seems satisfactory here.

\section{Discussion}
\label{sec:conclusion}

In this paper we propose a new version of the Gelman--Rubin diagnostic called $\hat{R}_\infty$, which improves MCMC convergence diagnostics on several aspects.
Firstly, it uses a localized version $\hat{R}(x)$ which assesses convergence at a given quantile $x$ of the target distribution.
Moreover, it is also based on a theoretical study of what $\hat{R}(x)$ is actually estimating:
assuming stationarity to focus only on the mixing property, the population version can be seen as a distance measure between the distributions of the chains.
This allows us to obtain convergence properties of $\hat{R}(x)$ and to tune the usual threshold of $1.01$ (Section~\ref{subsec:convergence-rhat}) based on a given confidence level and on the number of chains.
We show theoretically (Section~\ref{subsec:toy_examples}) and using experiments (Section~\ref{sec:simulations}) that our version is efficient to diagnose convergence.
Finally, we suggest a two-step algorithm for a multivariate diagnosis (Section~\ref{subsec:multivariate_def}), and reinforce the second step to consider all the directions of the space, as we show that the natural extension cannot be used directly (Section~\ref{sub-direction}).
Therefore, in the high-dimensional case where this computation is likely to be too expensive, we suggest to replace it by a univariate calculation on the log-likelihood or the log-posterior.
Diagnosing convergence in the multivariate case remains an open problem, and this is our hope that the local approach advocated here will trigger more research in this direction in the future.

\section*{Acknowledgement}

\gr{We would like to thank a Reviewer and an Editor for providing us with valuable comments that helped us improving the manuscript. Specifically, comments from an Editor allowed us to deal with multi-stage testing in an appealing and satisfactory way. S.~Girard acknowledges the support of the Chair Stress Test, Risk Management and Financial Steering, led by the  \'Ecole polytechnique and its Foundation and sponsored by BNP Paribas.
J. Arbel acknowledges the support of the French National Research Agency (ANR-21-JSTM-0001).
}

\begin{appendices}

\renewcommand\thefigure{\thesection.\arabic{figure}}
\renewcommand\thetable{\thesection.\arabic{table}}

\section{Construction of rank-$\hat{R}$ false negatives}
\label{app:counter_example}

In this section we detail different ways to construct false negative distributions on the chains that fool rank-$\hat{R}$, in the sense that they return rank-$\hat{R} \approx 1$ with non-identical distributions.

\subsection{Two-parameter distributions}

\begin{figure}
    \begin{tabular}{rl}
        \textbf{Example \ref{fig:counter_example_app}:} & Laplace $\mathcal{L}(0,\sigma)$ and $\mathcal{U}(-2\sigma, 2\sigma)$\\
        \vspace{-.3cm}
        \includegraphics[trim=0cm 0cm 0cm 0cm, clip, width = .49\textwidth]{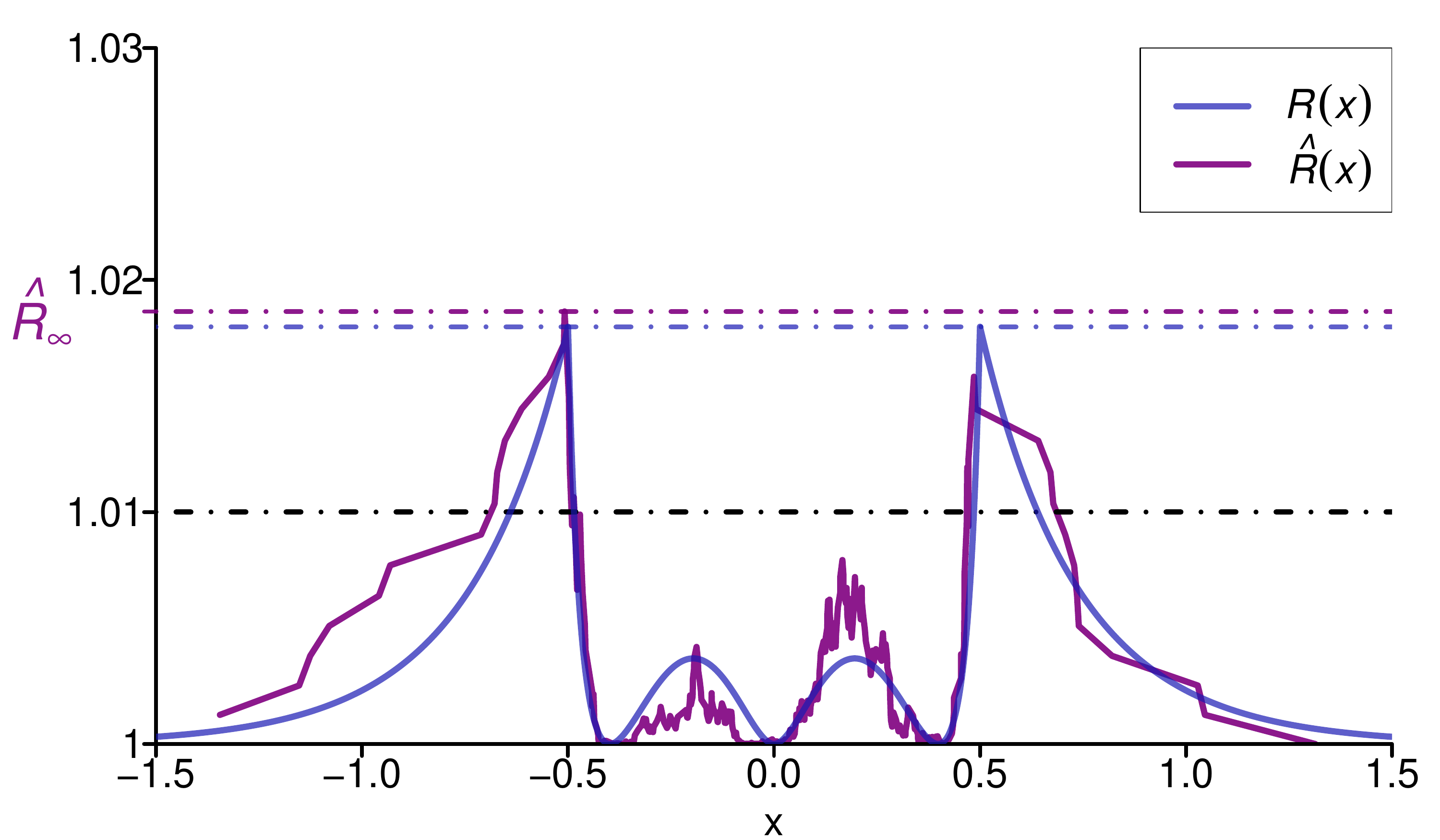} & 
        \includegraphics[trim=0cm 0.5cm 0cm 0cm, clip, width = .49\textwidth]{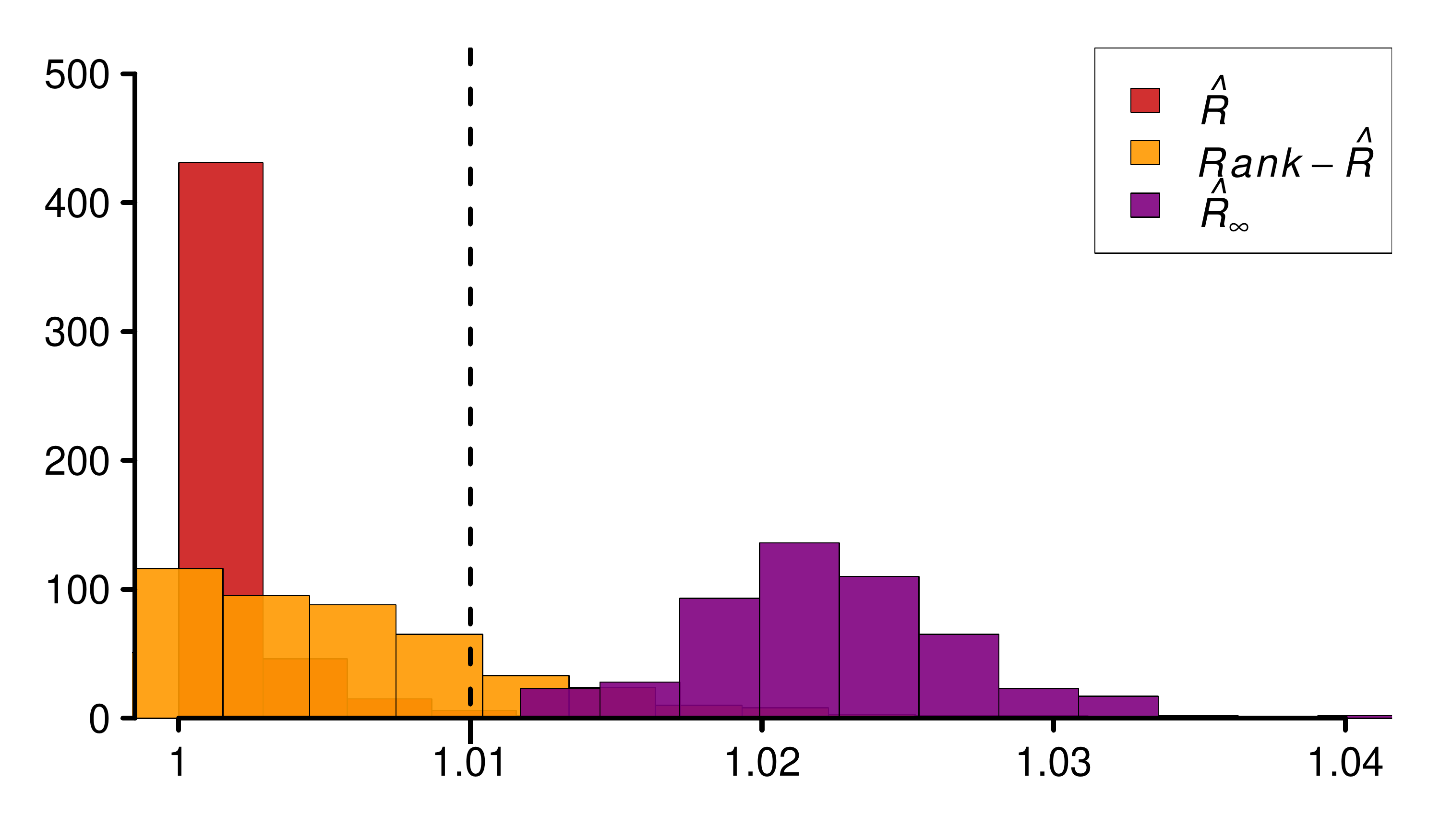}\\
        \vspace{-.3cm}
        & \hspace{-0.1cm}
        \footnotesize{\color{red2}rank-$\hat{R}$} 
        \hspace{-0.1cm} 
        \footnotesize{\color{orange2}$\hat{R}$} 
        \hspace{2.1cm} 
        \footnotesize{\color{violet2}$\hat{R}_\infty$}
    \end{tabular}
    \caption{Behaviour of $\hat{R}_\infty$ on the case of $m=2$ chains, $n=500$ each where $F_1$ is a Laplace distribution $\mathcal{L}(0,1/4)$ and $F_2$ is $\mathcal{U}(-1/2; 1/2)$.
    On the left: $\hat{R}(x)$ as a function of $x$ for one simulation.
    On the right: Histograms of $500$ replications of {\color{red2}$\hat{R}$}, {\color{orange2}rank-$\hat{R}$} of and {\color{violet2}$\hat{R}_\infty$}. 
    The dashed line corresponds to the threshold of $1.01$ which holds here for all versions of $\hat{R}$.
    }
    \label{fig:counter_example_app}
\end{figure}

As mentioned in Section~\ref{subsec:intro_comparison_rhat}, one way to fool rank-$\hat{R}$ is to fix two constraints on the distributions, which are identical mean and mean over the median on all chains. 
Many two-parameter distributions can be tuned to respect these constraints.
Consider for example a uniform chain $\mathcal{U}(-2\sigma, 2\sigma)$ and another one from a Laplace distribution $\mathcal{L}(0, \sigma)$, where $\sigma>0$.
The resulting $R_\infty$ does not depend on the scale parameter $\sigma$ and is given by
\begin{equation*}
    R_\infty = 
    \sqrt{1+\frac{1}{2(2e^2-1)}} \approx 1.018,
\end{equation*}
see Lemma~\ref{lem:laplace_uniform_proof}.
Note that computations are done in the case of $m=2$ chains, and so the traditional threshold of $1.01$ holds for $\hat{R}_\infty$ (see Table~\ref{tab:rhat_inf_tab}).
Therefore we expect $\hat{R}_\infty$ to diagnose convergence, and not $\hat{R}$ nor rank-$\hat{R}$.
Results are provided in Figure~\ref{fig:counter_example_app} and confirm this hypothesis.

\subsection{General framework with Generalized Pareto Distribution}

To find a general way to construct counter-examples, let us consider the Generalised Pareto Distribution (GPD), parametrised by $(\mu, \sigma, \xi)$: 
\begin{equation*}
    F_{\text{GPD}}(x) = 
    \begin{cases}
        1 - \left\{1 + \xi \left(\frac{x-\mu}{\sigma}\right) \right\}_+^{-\frac{1}{\xi}}
        &\quad \text{if $\xi \neq 0$} \,,\\
        1 - \exp\left(-\frac{x-\mu}{\sigma}\right)
        &\quad \text{if $\xi = 0$},
    \end{cases}
\end{equation*}
with $\{.\}_+ = \max \{0, .\}$. The support of this distribution depends on the parameters: $[\mu; +\infty)$ if $\xi \geq 0$, $[\mu; \mu - \sigma/\xi)$ otherwise. The expectation exists only if $\xi < 1$ and is equal to $\mu + \frac{\sigma}{1-\xi}$, and the median is given by $x_{\text{med}} = \mu + \sigma \frac{2^\xi - 1}{\xi}$.

An interesting property here is that conditioned on exceeding a value, the distribution is still a GPD distribution:
If $X \sim \text{GPD}(\mu, \sigma, \xi)$ then 
$X \mid X > u \sim \text{GPD}(u, \Tilde{\sigma}, \xi)$ with $\Tilde{\sigma} = \sigma + \xi(u-\mu)$.
Therefore, $X\mid X > x_{\text{med}} \sim \text{GPD}(u, \Tilde{\sigma}, \xi)$, and
\begin{align*}
    \mathbb{E}(X\mid X > x_{\text{med}}) &= x_{\text{med}} + \frac{\sigma + \xi(x_{\text{med}}-\mu)}{1-\xi},\\
    &= \mu + \frac{\sigma}{1-\xi}\left( 1 + \frac{2^\xi - 1}{\xi} \right),\\
    &= \mathbb{E}(X) + \sigma \frac{(2^\xi - 1)}{\xi(1-\xi)}.
\end{align*}
Then, by considering $(m-1)$ chains that follow a $\text{GPD}(\mu_1, \sigma_1, \xi_1)$, and one that follows a $\text{GPD}(\mu_2, \sigma_2, \xi_2)$, we can solve the system of two equations that links the two means and the two means over the medians, to obtain a range of possible parameters:
\begin{equation}
    \begin{cases}
        \mu_1 + \frac{\sigma_1}{1-\xi_1} = \mu_2 + \frac{\sigma_2}{1-\xi_2},\\
        \sigma_1 \frac{(2^{\xi_1}-1)}{\xi_1(1-\xi_1)} = \sigma_2 \frac{(2^{\xi_2}-1)}{\xi_2(1-\xi_2)}.
    \end{cases}
    \label{eq:system}
\end{equation}
In order to obtain different values of parameters, we should choose $\xi_1 \neq \xi_2$ and $\sigma_1 \neq \sigma_2$.
One way to characterize the set of solutions is as follows:
\begin{enumerate}
    \item Fix $\xi_1$ and $\xi_2$ such that $\xi_1 \neq \xi_2$, and define $\lambda = \frac{f(\xi_1)}{f(\xi_2)}$ with $f(\xi) = \frac{(2^\xi - 1)}{\xi(1-\xi)}$.
    \item Fix $\sigma_1$, and using the second equation of (\ref{eq:system}) define $\sigma_2 = \lambda\sigma_1$.
    \item Finally, choose $\mu_1$ and $\mu_2$ such that $\mu_1 - \mu_2 = \sigma_1 (\frac{\lambda}{1-\xi_2} - \frac{1}{1-\xi_1})$.
\end{enumerate}
An example of a solution is as follows:
\begin{enumerate}
    \item Choose $(\mu_1, \sigma_1, \xi_1) = (0,1,0)$, the standard exponential distribution Exp$(1)$, and $\xi_2 = -1$, a uniform distribution.
    \item Following the method described before, we obtain $\lambda = \frac{f(\xi_1)}{f(\xi_2)} = 4\log(2)$, so $\sigma_2 = \lambda \sigma_1 = 4\log(2)$.
    \item Following the last point, we set $\mu_2 = \mu_1 - \sigma_1 (\frac{\lambda}{1-\xi_2} - \frac{1}{1-\xi_1}) = 1 - 2\log(2)$.
\end{enumerate}
To conclude, if some chains follow an Exp$(1)$ distribution and the other ones a $\mathcal{U}(1 - 2\log(2); 1 + 2\log(2))$ distribution, the difference between the chains should not be detected by the rank-$\hat{R}$, which is illustrated in the last column of Figure~\ref{fig:example_known_dist}.
Similarly to the other examples in Section~\ref{subsec:toy_examples}, $\hat{R}_\infty$ manages to diagnose the convergence issue contrary to other versions.

\section{Proofs}

\subsection{Proofs in the univariate case}

\paragraph{Proof of Proposition~\ref{prop-calcul-R}.}
    
    Let$ x\in\mathbb{R}$. The within- and between-variances are given by:
    \begin{align}
    W(x) &= \mathbb{E}\left[\text{Var}[I_x \mid Z ]\right] \nonumber\\
        \label{eq:theoretical_W_int}
        &= \frac{1}{m} \sum_{j=1}^m \left(\mathbb{E}\left[I_x^2 \mid Z = j \right] - \mathbb{E}^2\left[I_x \mid Z = j \right]\right)\\
        &= \frac{1}{m} \sum_{j=1}^m \left(\mathbb{P}\left(\theta \leq x \mid Z = j \right) - \mathbb{P}^2\left(\theta \leq x \mid Z = j \right) \right) \nonumber \\
        \label{eq:theoretical_W}
        &= \frac{1}{m} \sum_{j=1}^m \left(F_j(x) - F^2_j(x)\right), \\
    B(x) &= \text{Var}\left[\mathbb{E}[I_x \mid Z ]\right]
        = \mathbb{E}\left[F_Z(x)^2\right] - \mathbb{E}[F_Z(x)]^2 \nonumber\\
        \label{eq:theoretical_B_int}
        &= \frac{1}{m} \sum_{j=1}^m F^2_j(x) - \left(\frac{1}{m} \sum_{j=1}^m F_j(x)\right)^2\\
        &= \frac{m-1}{m^2} \sum_{j=1}^m F^2_j(x) - \frac{2}{m^2}\sum_{j<k}F_j(x)F_{k}(x)\nonumber\\
        &= \frac{1}{m^2} \sum_{j<k} (F^2_j(x) + F^2_{k}(x) - 2F_j(x)F_{k}(x))\nonumber\\
                \label{eq:theoretical_B}
        &= \frac{1}{m^2}\sum_{j<k} \left(F_j(x)-F_{k}(x)\right)^2,
    \end{align}
  where (\ref{eq:theoretical_W_int}) and (\ref{eq:theoretical_B_int}) are a consequence of $\mathbb{P}(Z = j) = 1/m$ for all $j \in \{1, \cdots, m\}$.
  The conclusion follows.

  \paragraph{Proof of Proposition~\ref{prop-first}.} Proofs of (i) and (ii) are straightforward, let us focus on (iii) and (iv).
    
    (iii): Without loss of generality, we assume that all $F_j$'s are defined on $\mathbb{R}$ and we denote $\thickbar{F}_j = 1 - F_j$ the associated survival function.
    From Proposition~\ref{prop-calcul-R}, we can write, for any $x\in\mathbb{R}$:
    \begin{equation*}
        R^2(x)
        = \frac{\left(\frac{1}{m}\sum_{j=1}^m F_j(x)\right) \left(\frac{1}{m}\sum_{j=1}^m \thickbar{F}_j(x)\right)}{\frac{1}{m}\sum_{j=1}^m \thickbar{F}_j(x)F_j(x)}.
    \end{equation*}
    Denote by $f(x) \simeq g(x)$ two functions $f$ and $g$ that are asymptotically equivalent when $x \to a$, $a \in \mathbb{R} \cup \{-\infty, +\infty\}$. Consider first the case $x\to-\infty$.
    Clearly, for all $j\in\{1,\dots,m\}$, $F_j(x) \to 0$ and $\thickbar{F}_j(x) \to 1$ as $x\to-\infty$, so that
    \begin{align*}
        \left(\frac{1}{m}\sum_{j=1}^m F_j(x)\right) \left(\frac{1}{m}\sum_{j=1}^m \thickbar{F}_j(x)\right) &\simeq \frac{1}{m}\sum_{j=1}^m {F}_j(x), \\
        \text{and} \quad \frac{1}{m}\sum_{j=1}^m \thickbar{F}_j(x)F_j(x) &\simeq \frac{1}{m}\sum_{j=1}^m {F}_j(x),
    \end{align*}
    leading to $R(x) \to 1$ as $x\to-\infty$. 
    Second, remarking that $R(x)$ is symmetric with respect to $(F_j(x), \thickbar{F}_j(x))$, we similarly have $R(x) \to 1$ as $x\to\infty$.\\
    
    (iv): It is clear that $B(x)$ and $W(x)$ are continuous since the $F_j$'s are continuous, so the only thing to prove is that the denominator $W(x)$ never vanishes, except if extending by continuity is possible. 
    Remarking that for all $x$, we have $F_j(x)(1-F_j(x)) \geq 0$, $W(x) = 0$ if and only if $F_j(x) = 0$ or $F_j(x) = 1$ for $j \in \{1, \ldots, m\}$.
    Almost all combinations are avoided by the assumption that the supports must overlap, except $F_1(x) = \cdots = F_m(x) = 0$ and $F_1(x) = \cdots = F_m(x) = 1$. In these two latter cases, an extension by continuity is possible as $R(x) = 1$ from (iii).

    \paragraph{Proof of Proposition~\ref{prop:asymptotic_R-hat}.} 
    (i): Assume for the sake of simplicity that $F$ is continuous and strictly increasing.
    Let us show that the corresponding $\hat{R}_\infty$ has the same distribution as in the standard uniform case. In view of~(\ref{eq:R_theorique}), $\hat{R}_\infty$ can be written as
    \begin{equation*}
        \hat{R}_\infty = \sup_{x\in {\mathbb R}} \hat{R}(x) 
        = \sup_{x\in {\mathbb R}}  \sqrt{1 + \frac{\sum_{j=1}^m\sum_{k=j+1}^m \left(\hat{F}_j(x)-\hat{F}_k(x)\right)^2}{m\sum_{j=1}^m \hat{F}_j(x)(1-\hat{F}_j(x))}},
    \end{equation*}
    where for any $j\in \{1,\dots,m\}$, 
    \begin{equation*}
        \hat{F}_j(x) 
        = \frac{1}{n}\sum_{i=1}^n \mathbb{I}\{\theta^{(i,j)} \leq x \}
        = \frac{1}{n}\sum_{i=1}^n \mathbb{I}\{F(\theta^{(i,j)}) \leq F(x) \},
    \end{equation*}
    and where the random variables $F(\theta^{(i,j)})$ are standard uniformly distributed for any $i\in\{1,\dots,n\}$ and $j\in\{1,\dots,m\}$. 
    Finally, observing that 
    \begin{equation*}
        \sup_{x\in {\mathbb R}} \hat{R}(x) 
        = \sup_{y\in [0,1]} \hat{R}(F^{-1}(y)) 
    \end{equation*}
    concludes the proof.
    
    (ii): Let $x \in \mathbb{R}$ and remark that~(\ref{eq:CLT_MC_j}) implies
    that $\hat F_j(x)\toP F(x)$ as $n\to\infty$ for all $j\in\{1,\dots,m\}$.
    Then it follows from~(\ref{eq:theoretical_W}) in the proof of Proposition~\ref{prop-calcul-R} that
    \begin{equation}
        \label{eq-part1}
        \hat{W}(x)  = \frac{1}{m} \sum_{j=1}^m \hat{F}_j(x)(1-\hat{F}_j(x)) = F(x) (1-F(x)) + o_p(1).
    \end{equation}       
    Besides, from~(\ref{eq:CLT_MC_j}) and~(\ref{eq-def-local-ESS}), for $j = 1, \ldots, m$, $\hat{F}_j(x)$ can be expanded as
    $$
    \hat{F}_j(x) = F(x) + \sqrt{\frac{F(x)(1-F(x))}{\text{ESS}(x)/m}} \xi_{j,n},
    $$
    with $\xi_{j,n}\tod\mathcal{N}(0,1)$ as $n\to\infty$.
    Thus (\ref{eq:theoretical_B}) in the proof of Proposition~\ref{prop-calcul-R} entails
\begin{equation}
    \label{eq-part2}
\hat{B}(x) 
        = \frac{1}{m^2}\sum_{j<k} \left(\hat{F}_j(x) - \hat{F}_k(x)\right)^2
        = \frac{F(x)(1-F(x))}{\text{ESS}(x)}\times \frac{1}{m} \sum_{j<k} \left(\xi_{j,n} - \xi_{k,n}\right)^2.
\end{equation}
Introducing the random vector $\boldsymbol{\xi}_n = (\xi_{1,n}, \ldots, \xi_{m,n})^\top \tod \mathcal{N}(0,\boldsymbol{I}_m)$ where $\boldsymbol{I}_m$
is the $m\times m$ identity matrix, one has
    \begin{align*}
     \frac{1}{m}   \sum_{j<k} \left(\xi_{j,n} - \xi_{k,n}\right)^2
        = \frac{m-1}{m} \sum_{j=1}^m \xi_{j,n}^2 - \frac{1}{m}\sum_{j\neq k} \xi_{j,n} \xi_{k,n} 
        = \boldsymbol{\xi}_n^\top \boldsymbol{A} \boldsymbol{\xi}_n,
    \end{align*}
    with $\boldsymbol{A} =  \boldsymbol{I}_m - \boldsymbol{J}_m/m$, and where $\boldsymbol{J}_m$ is the $m\times m$ matrix filled with ones.
    The symmetric matrix 
    $\boldsymbol{A}$ can be eigen-decomposed as
    $\boldsymbol{A} = \boldsymbol{Q}^\top \boldsymbol{\Lambda} \boldsymbol{Q}$ with $\boldsymbol{Q}$ an orthogonal matrix and
        $\boldsymbol{\Lambda} = \text{diag}(\lambda_1, \cdots, \lambda_d) = \text{diag}(1, \ldots, 1, 0)$.
    Remark that $\boldsymbol{U}_n \coloneqq  \boldsymbol{Q}\boldsymbol{\xi}_n \tod \mathcal{N}(0, \boldsymbol{I}_m)$ so that
\begin{equation}
    \label{eq-part3}
    \frac{1}{m}  \sum_{j<k} \left(\xi_{j,n} - \xi_{k,n}\right)^2
        = (\boldsymbol{Q}\boldsymbol{\xi}_n)^\top \boldsymbol{\Lambda} (\boldsymbol{Q}\boldsymbol{\xi}_n) 
        = \sum_{j=1}^m \lambda_{j} U_{j,n}^2
        =  \sum_{j=1}^{m-1} U_{j,n}^2 \tod \chi^2_{m-1},
    \end{equation}
as $n\to\infty$.    Collecting~(\ref{eq-part1}),~(\ref{eq-part2}) and~(\ref{eq-part3}) yields
    \begin{equation*}
  \text{ESS}(x)(\hat{R}^2(x) - 1) =    \text{ESS}(x)   \frac{\hat{B}(x)}{\hat{W}(x)} \tod \chi^2_{m-1},
    \end{equation*}
  as $n\to\infty$ and the result is proved.

\subsection{Proofs in the multivariate case}
    \label{sec:proof_multivariate}
    \begin{LemApp}[Standardization of margins]
        \label{lem-invariance}
        Assume the assumptions of Proposition~\ref{prop-calcul-R} hold.
        If the margins of $F_1, \ldots, F_m$ coincide, then the multivariate $R_\infty$ is the same as the one calculated on the associated copulas $C_1,\dots, C_m$.
    \end{LemApp}

\paragraph{Proof.} Denote by $\phi_1, \ldots, \phi_d$ the common margins of the cdf's $F_1, \ldots, F_m$, 
so that $F_j(\boldsymbol{x})=C_j(\phi_1(x_1),\dots,\phi_d(x_d))$ for any $j\in\{1,\dots,m\}$.
Letting $\boldsymbol{y} = (\phi_1(x_1), \ldots, \phi_d(x_d))$, we have
\begin{equation*}
    \sup_{\boldsymbol{x}\in {\mathbb R}^d} R(\boldsymbol{x}) = \sup_{\boldsymbol{y}\in[0,1]^d} R(\phi^{-1}_1(y_1), \ldots, \phi^{-1}_d(y_d)).
\end{equation*}
Besides, in view of~(\ref{eq:R_theorique_multi}),  $m(R^2(\phi^{-1}_1(y_1), \ldots, \phi^{-1}_d(y_d)) - 1)$ can be written as
\begin{multline*}
   \frac{\sum_{j=1}^m\sum_{k=j+1}^m \left(F_j(\phi^{-1}_1(y_1), \ldots, \phi^{-1}_d(y_d))-F_k(\phi^{-1}_1(y_1), \ldots, \phi^{-1}_d(y_d))\right)^2}{\sum_{j=1}^m F_j(\phi^{-1}_1(y_1), \ldots, \phi^{-1}_d(y_d))(1-F_j(\phi^{-1}_1(y_1), \ldots, \phi^{-1}_d(y_d)))} \\
     =\frac{\sum_{j=1}^m\sum_{k=j+1}^m \left(C_j(y_1, \ldots, y_d)-C_k(y_1, \ldots, y_d)\right)^2}{\sum_{j=1}^m C_j(y_1, \ldots, y_d)(1-C_j(y_1, \ldots, y_d))},
\end{multline*}
which proves the result.

\paragraph{Proof of Lemma~\ref{lem:copule_bound}.} 
Let function $\Tilde{R}$ be defined on $[0,1]^2$ by 
\begin{equation}
    \Tilde{R}(c_1, c_2) = \sqrt{1+\frac{1}{2}\frac{(c_1-c_2)^2}{c_1(1-c_1) + c_2(1-c_2)}}.        
\end{equation}
    For any $c_1 \in [0,1]$, $c_2 \mapsto \Tilde{R}(c_1,c_2)$ is decreasing on $[0,c_1]$, and increasing on $[c_1,1]$.
    Let $(c_{-}, c_{+}) \in [0,1]^2$ and $(c_1, c_2) \in [c_-,c_+]^2$.
    Without loss of generality, assume that $c_1 \leq c_2$.
    Let $c\in[0,c_2]$. Since $c_2 \leq c_{+}$ and $\Tilde{R}(c,\cdot)$ is increasing on $[c,1]$, we have $\Tilde{R}(c, c_2) \leq \Tilde{R}(c, c_{+})$.
    In particular, for $c = c_{-}$:
    \begin{equation}
        \Tilde{R}(c_{-}, c_2) \leq \Tilde{R}(c_{-}, c_{+}).
        \label{eq:proof_lemma_c1}
    \end{equation}
    Moreover, $\Tilde{R}(c_2, \cdot)$ is decreasing on $[0,c_2]$ and $c_{-} \leq c_1 \leq c_2$, we also have
    \begin{equation}
        \Tilde{R}(c_{-}, c_2) =  \Tilde{R}(c_2, c_{-}) \geq \Tilde{R}(c_2, c_1).
        \label{eq:proof_lemma_c2}
    \end{equation}
    Combining (\ref{eq:proof_lemma_c1}) and (\ref{eq:proof_lemma_c2}), we finally obtain $\Tilde{R}(c_2, c_1) \leq \Tilde{R}(c_{-}, c_{+})$, which concludes the proof.

\paragraph{Proof of Proposition~\ref{prop:global_bound}.}
By definition of lower and upper Fréchet--Hoeffding bounds \citep{Nelsen2006}, $W_d(\boldsymbol{u})\leq C(\boldsymbol{u}) \leq M_d(\boldsymbol{u})$ for any copula $C$ and $\boldsymbol{u}\in[0,1]^d$.
Thus in view of Lemma~\ref{lem:copule_bound} it only remains to prove that $R_\infty(W_d, M_d) = \sqrt{(d+1)/{2}}$.
To this end, let $p$ be the index such that $u_p = \min\left\{ u_1, \ldots, u_d \right\}$. Two cases arise:

(i) If $1 - d + \sum_{i=1}^d u_i \leq 0$, then 
\begin{equation}
\label{eq-cas1}
    f(\boldsymbol{u}) \coloneqq 2(R^2(\boldsymbol{u}) - 1) = \frac{u^2_p}{u_p(1-u_p)} = \frac{1}{{1}/{u_p}-1}.
\end{equation}
As a consequence, $R^2(\boldsymbol{u})$ is maximum when $u_p$ is maximum under the constraints 
$$\begin{cases}
    u_p \leq u_i \quad \forall i \neq p,\\
    u_p \leq d-1-\sum_{i\neq p} u_i.
\end{cases}$$
It is easily seen that the maximum occurs in the equality case $u_1 = \cdots = u_d = (d-1)/{d}$ and thus $  2(R^2(\boldsymbol{u}) - 1)=d-1$.

(ii) Conversely, if $1 - d + \sum_{i=1}^d u_i \geq 0$, then
\begin{equation*}
    f(\boldsymbol{u}) = \frac{\left(1 - d + \sum_{i\neq p} u_i\right)^2}{u_p(1-u_p) + (1 - d + \sum_{i = 1}^d u_i)(d-\sum_{i=1}^d u_i)}.
\end{equation*}
The problem therefore amounts to maximising $f(\boldsymbol{u})$ under the constraints
\begin{equation}
\label{eq:constraints1}
\begin{cases}
    \boldsymbol{u} \in [0,1]^d,\\
    u_p \leq u_i \quad \forall i \neq p,\\
    1-d + \sum_{i=1}^d u_i \geq 0.
\end{cases}
\end{equation}
The associated Lagrangian can be written with $\boldsymbol{\lambda}_1, \boldsymbol{\lambda}_2, \boldsymbol{\lambda}_3 \in \mathbb{R}^d$ and $\lambda_4 \in \mathbb{R}$ as:
\begin{equation*}
    \label{eq:lagrangian1}
    \mathfrak{L}(\boldsymbol{u}, \boldsymbol{\lambda}_1, \boldsymbol{\lambda}_2, \boldsymbol{\lambda}_3, \lambda_4) = f(\boldsymbol{u}) 
    + \sum_{i=1}^d \lambda_{1,i} u_i
    + \sum_{i=1}^d \lambda_{2,i} (1-u_i) 
    + \sum_{i\neq p} \lambda_{3,i} (u_i-u_p) 
    + \lambda_{4} \left(1-d+ \sum_{i=1}^d u_i \right).
\end{equation*}
The first-order conditions are given by (\ref{eq:constraints1}) and $\nabla_{\boldsymbol{u}} \mathfrak{L}(\boldsymbol{u}, \boldsymbol{\lambda}_1, \boldsymbol{\lambda}_2, \boldsymbol{\lambda}_3, \lambda_4) = 0$,
and the Karush--Kuhn--Tucker conditions are
\begin{equation}
\label{eq:kkt_conditions1}
\begin{cases}
    \lambda_{1,i} \geq 0 \quad &\text{with} \quad \lambda_{1,i} u_i = 0, \quad \forall i\in\{1,\dots,d\},\\
    \lambda_{2,i} \geq 0 \quad &\text{with} \quad \lambda_{2,i}(1-u_i) = 0, \quad \forall i\in\{1,\dots,d\},\\
    \lambda_{3,i} \geq 0 \quad &\text{with} \quad \lambda_{3,i}(u_i - u_p) = 0, \quad \forall i\in\{1,\dots,d\} \quad \text{s.t.} \quad i\neq p,\\
    \lambda_4 \geq 0 \quad &\text{with} \quad \lambda_4 \left(1-d + \sum_{i=1}^d u_i\right) = 0.
\end{cases}
\end{equation}
We distinguish different cases: 
\begin{itemize}
    \item If $\lambda_{1,i_0} \neq 0$ for some $i_0 \in \{1, \ldots, d\}$, then $u_{i_0} = 0$ and combined with (\ref{eq:constraints1}) we obtain necessarily $i_0 = p$ and $1-d + \sum_{i=1}^d u_i = 0$, leading to a non-optimal solution for $f(\boldsymbol{u})$.
    \item If $\lambda_{2,i_0} \neq 0$ for some $i_0 \in \{1, \ldots, d\}$, then $u_{i_0} = 1$, and the problem is exactly the same written in dimension $d-1$, and a recurrence proves that the maximum is equal to $\sqrt{(d+1)/2}$ in dimension $d$. So the maximum is increasing with the dimension and therefore is not reached in this case.
\end{itemize}
One can thus assume $\boldsymbol{\lambda}_1 = 0$ and $\boldsymbol{\lambda}_2 = 0$.
Moreover, for all $(i,j)$ such that $i \neq p$ and $j \neq p$, we have $\frac{\partial f}{\partial u_i} = \frac{\partial f}{\partial u_j}$. 
Combined with $\nabla_{\boldsymbol{u}} \mathfrak{L} = 0$, we obtain $\lambda_{3,i} = \lambda_{3,j}$, which leads to considering the simplified Lagrangian:
\begin{equation*}
    \label{eq:lagrangian2}
    \mathfrak{L}(\boldsymbol{u}, \lambda_3, \lambda_4) = f(\boldsymbol{u})
    + \lambda_{3} \left(\sum_{i=1}^d u_i - d u_p\right) 
    + \lambda_{4} \left(1-d+ \sum_{i=1}^d u_i \right).
\end{equation*}
In that form, the function $f$ and the constraints of the Lagrangian can be written only as a function of $(1-d+\sum_{i=1}^d u_i, u_p, \lambda_3, \lambda_4)$. 
Let $\Tilde{\mathfrak{L}}(x = 1-d+\sum_{i=1}^d u_i, y = u_p, \lambda_3, \lambda_4) = \mathfrak{L}(\boldsymbol{u}, \lambda_3, \lambda_4)$.
Since for all $i$
\begin{equation*}
    \frac{\partial \mathfrak{L}}{\partial u_i} 
    = \frac{\partial \Tilde{\mathfrak{L}}}{\partial x}\frac{\partial x}{\partial u_i} + \frac{\partial \Tilde{\mathfrak{L}}}{\partial y}\frac{\partial y}{\partial u_i}
    = \frac{\partial \Tilde{\mathfrak{L}}}{\partial x} + \frac{\partial \Tilde{\mathfrak{L}}}{\partial y}\mathbb{I}\{i=p\},
\end{equation*}
solving $\nabla_{\boldsymbol{u}}\mathfrak{L} = 0$ is equivalent to solving $\nabla_{(x,y)}\mathfrak{\Tilde{L}} = 0$.
The corresponding problem is therefore 
\begin{equation}
    \max \frac{(x-y)^2}{x(1-x) + y(1-y)}, \quad \text{under the constraints} \quad
    \begin{cases}
        x \geq \max\{0, 1-d+dy\},\\
        x \leq y \leq 1.
    \end{cases} 
    \label{eq:opti_prob_2D}
\end{equation}
Combining the constraints with the study of function in the proof of Lemma~\ref{lem:copule_bound} leads to the solution 
$(x,y) = (0, \frac{d-1}{d})$. 
So 
$$
\sum_{i=1}^d u_i = d-1 \quad \text{ and } \min\{u_1, \ldots, u_d\} = \frac{d-1}{d},
$$ 
so necessarily $u_1 = \cdots = u_d = \frac{d-1}{d}$, which concludes the proof.

\paragraph{Proof of Corollary~\ref{cor:global_bound_d}} 
For all $\boldsymbol{u}\in[0,1]^d$, we have 
\begin{align*}
    m\left(R^2(\boldsymbol{u}) -1 \right)
    &= \sum_{j=1}^m \sum_{k=j+1}^m \frac{\left(C_j(\boldsymbol{u}) - C_k(\boldsymbol{u})\right)^2}{\sum_{\ell=1}^m C_\ell(\boldsymbol{u})(1-C_\ell(\boldsymbol{u}))},\\
    &\leq \sum_{j=1}^m \sum_{k=j+1}^m \frac{\left(C_j(\boldsymbol{u}) - C_k(\boldsymbol{u})\right)^2}{C_j(\boldsymbol{u})(1-C_j(\boldsymbol{u})) + C_k(\boldsymbol{u})(1-C_k(\boldsymbol{u}))},\\
    &= \sum_{j=1}^m \sum_{k=j+1}^m 2\left(R^2(W_d(\boldsymbol{u}), M_d(\boldsymbol{u})) -1 \right),
\end{align*}
so that an upper bound on the multivariate $R_\infty$ can be expressed thanks to bivariate as follows:
\begin{align*}
    m\left(R_\infty^2(C_1, \ldots, C_m) -1 \right)
    &\leq \sum_{j=1}^m \sum_{k=j+1}^m 2\left(R_\infty^2(C_j, C_k) -1 \right),\\
    &\leq \sum_{j=1}^m \sum_{k=j+1}^m 2\left(R_\infty^2(W_d, M_d) -1 \right),\\
    &= \frac{m(m-1)}{2} (d-1),
\end{align*}
using Lemma~\ref{lem:copule_bound} and Proposition~\ref{prop:global_bound}. The result is thus proved.

\paragraph{Proof of Corollary~\ref{cor:PLOD_bound}.} One can 
prove that the maximum of $\boldsymbol{u}\in[0,1]^d\mapsto R(\Pi_d(\boldsymbol{u}), M_d(\boldsymbol{u}))$ is reached at $u_1 = \cdots = u_d \coloneqq u \in [0,1]$,
which leads to studying the maximum of 
$$
f_d(u) \coloneqq  2(R^2(u,\ldots,u)-1)=\frac{(u^d-u)^2}{u^d(1-u^d) + u(1-u)}.
$$
The first derivative of $f_d$ is proportional to 
    $$g_d(u) = -2(d-1)u^{2d-1} + d u^{2d-2} - 2(d-1)u^{d} + 3(d-1)u^{d-1} -1.$$
Routine calculations show the existence of a unique root in $[0,1]$, but finding the explicit value does not seem possible when $d>2$
since  $g_d$ is a polynomial of order $2d-1$. We thus restrict ourselves to an asymptotic analysis when $d\to\infty$.
First, Lemma~\ref{lem:copule_bound} and Proposition~\ref{prop:global_bound} entail that $$\max\limits_{u\in[0,1]}  f_d(u) \leq 2(R^2_\infty(W_d,M_d) - 1) = d-1.$$
Second, a lower bound can be obtained by letting $u_d = 1 - (\log d)/{d}$. Indeed, $u_d\to 1$ and
$u_d^d =  \exp\left(-(\log d)(1+o(1))\right) \to 0$  so that the numerator of $f_d(u)$ satisfies 
$(u_d^d - u_d)^2 \to 1$ as $d\to\infty$. Moreover, the denominator satisfies 
$$
u_d^d(1-u_d^d) + u_d(1-u_d)
= \exp\left(-\log(d)(1+o(1))\right)(1+o(1)) + \frac{\log d}{d}(1+o(1)) =\frac{\log d}{d}(1+o(1)),
$$
as $d\to\infty$. As a consequence, 
$f_d(u_d) = \frac{d}{\log d}(1+o(1))$ and we have proved that 
$$        
\frac{d}{\log d}(1+o(1)) \leq \max\limits_{u\in[0,1]}\ f_d(u) \leq d-1.
$$     
The result follows.

\paragraph{Proof of Corollary~\ref{cor:NLOD_bound}.} 
From Proposition~\ref{lem:copule_bound} and the definition of a NLOD copula, the proof reduces to calculating $R_\infty(W_d,\Pi_d)$. 
Two cases are considered.

(i) First, if $1 - d + \sum_{i=1}^d u_i \leq 0$, then $W_d(\boldsymbol{u}) = 0$ and 
\begin{equation*}
    2(R^2(\boldsymbol{u}) - 1) = \frac{1}{\frac{1}{\prod_{i=1}^d u_i}-1}.
\end{equation*}
The maximisation of $2(R^2(\boldsymbol{u}) - 1)$ is then equivalent to solving:
\begin{equation}
    \max \prod_{i=1}^d u_i, \quad \text{under the constraints} \quad
    \begin{cases}
        \boldsymbol{u} \in [0,1]^d,\\
        1-d + \sum_{i=1}^d u_i \leq 0.
    \end{cases} 
    \label{eq:opti_prob}
\end{equation}
Since the constraints are linear and the objective function is convex, the above optimization problem is convex.
The Lagrangian associated with~(\ref{eq:opti_prob}) can be written with $\lambda_1 \in \mathbb{R}$, $\boldsymbol{\lambda}_2, \boldsymbol{\lambda}_3 \in \mathbb{R}^d$, as:
\begin{equation*}
    \mathfrak{L}(\boldsymbol{u}, \lambda_1, \boldsymbol{\lambda}_2, \boldsymbol{\lambda}_3) = \prod_{i=1}^d u_i - \lambda_1 \left(1-d + \sum_{i=1}^d u_i\right)
    - \sum_{i=1}^{d} \lambda_{2,i} (u_i-1)
    + \sum_{i=1}^{d} \lambda_{3,i} u_i,
\end{equation*}
The first-order conditions are
\begin{equation}
\label{eq:first_order_cond}
\begin{cases}
    \nabla_{\boldsymbol{u}} \mathfrak{L}(\boldsymbol{u}, \lambda_1, \boldsymbol{\lambda}_2, \boldsymbol{\lambda}_3) = 0,\\
    0 \leq u_i \leq 1 \quad \forall i\in\{1,\dots,d\},\\
    1-d + \sum_{i=1}^d u_i \leq 0,
\end{cases}
\end{equation}
and the Karush--Kuhn--Tucker conditions are
\begin{equation}
\label{eq:kkt_conditions}
\begin{cases}
    \lambda_1 \geq 0 \quad &\text{with} \quad \lambda_1 (1-d + \sum_{i=1}^d u_i) = 0,\\
    \lambda_{2,i} \geq 0 \quad &\text{with} \quad \lambda_{2,i}(u_i-1) = 0, \quad \forall i\in\{1,\dots,d\},\\
    \lambda_{3,i} \geq 0 \quad &\text{with} \quad \lambda_{3,i} u_i = 0, \quad \forall i\in\{1,\dots,d\}.
\end{cases}
\end{equation}
If there exists $i_0$ such that $u_{i_0} = 0$, then $\prod_{i=1}^d u_i = 0$, which is clearly non-optimal. 
Thus,~(\ref{eq:kkt_conditions}) implies $\lambda_{3,i} = 0$ for all $i \in \{1, \ldots, d\}$.
Moreover,  for all $i \in \{1, \ldots, d\}$,
\begin{equation*}
    \frac{\partial \mathcal{L}}{\partial u_i} = 0 \implies \lambda_1 + \lambda_{2,i} = \prod_{j\neq i} u_j, 
\end{equation*}
so that, for all $(j,k) \in \{1, \ldots, d\}^2$, 
\begin{equation}
    \label{eq:lambda_constant}
    (\lambda_1 + \lambda_{2,j})u_j = (\lambda_1 + \lambda_{2,k})u_k.
\end{equation}
If there exists $i_0$ such that $\lambda_{2,i_0} \neq 0$, then $u_{i_0}=1$ from~(\ref{eq:kkt_conditions}) and consequently, for all $k \neq i_0$,
\begin{equation*}
    u_k = \frac{\lambda_1 + \lambda_{2,i_0}}{\lambda_1 + \lambda_{2,k}}.
\end{equation*}
 Taking account of~(\ref{eq:first_order_cond}) yields $u_k \leq 1$ which implies in turn $\lambda_{2,k} \geq \lambda_{2,i_0} \neq 0$ and $u_k = 1$ for all $k \in \{1, \ldots, p\}$ from~(\ref{eq:kkt_conditions}). The resulting $\boldsymbol{u}$ does not fulfil the third constraint in (\ref{eq:first_order_cond}).
Then, necessarily, $\lambda_{2,i} = 0$ for all $i\in \{1, \ldots, d\}$, and combining with~(\ref{eq:lambda_constant}), it follows that the optimum is reached when $u_1 = \cdots = u_d = u$. Replacing in~(\ref{eq:opti_prob}) yields the optimization problem
\begin{equation}
    \max u^d, \quad \text{under the constraints} \quad
    \begin{cases}
        0 \leq u \leq 1, \\
        1 + d(u-1) \leq 0,
    \end{cases} 
    \label{eq:opti_prob_1D}
\end{equation}
whose solution is $u = (d-1)/{d}$.

(ii) Second, if $1 - d + \sum_{i=1}^d u_i \geq 0$, the problem to solve is
\begin{equation*}
    \max f_d(\boldsymbol{u}) \coloneqq \frac{\left(\prod_{i=1}^d u_i - 1 + d - \sum_{i=1}^d u_i\right)^2}{\left(\prod_{i=1}^d u_i\right)\left(1-\prod_{i=1}^d u_i\right) + \left(1 - d + \sum_{i=1}^d u_i\right)\left(d - \sum_{i=1}^d u_i\right)},
\end{equation*}
\begin{equation*}
    \text{under the constraints} \quad
    \begin{cases}
        \boldsymbol{u} \in [0,1]^d, \\
        1 - d + \sum_{i=1}^d u_i \geq 0.
    \end{cases} 
\end{equation*}
The Lagrangian can be written, with $\lambda_1 \in \mathbb{R}$, $\boldsymbol{\lambda}_2, \boldsymbol{\lambda}_3 \in \mathbb{R}^d$:
\begin{equation*}
    \mathfrak{L}(\boldsymbol{u}, \lambda_1, \boldsymbol{\lambda}_2, \boldsymbol{\lambda}_3) = f_d(\boldsymbol{u})
    + \lambda_1 \left(1-d + \sum_{i=1}^d u_i\right)
    + \sum_{i=1}^{d} \lambda_{2,i} (1-u_i)
    + \sum_{i=1}^{d} \lambda_{3,i} u_i.
\end{equation*}
In the same way as in the proof of Proposition~\ref{prop:global_bound}, we can focus on the solution such that $\boldsymbol{\lambda}_2 = 0$ and $\boldsymbol{\lambda}_3 = 0$.
The first order condition $\nabla_{\boldsymbol{u}}\mathfrak{L} = 0$ leads to the solution
\begin{equation*}
    u_j = \prod_{i=1}^d u_i \times \frac{2r(\boldsymbol{u}) - \left(\prod_{i=1}^d u_i - 1 + d - \sum_{i=1}^d u_i \right)\left(1-2\prod_{i=1}^d u_i\right)}{2r(\boldsymbol{u}) + \left(\prod_{i=1}^d u_i - 1 + d - \sum_{i=1}^d u_i \right)\left(1-2(1 - d + \sum_{i=1}^d u_i)\right) - \lambda_1},
\end{equation*}
for $j \in \{1, \ldots, d\}$ with $r(\boldsymbol{u}) = \prod_{i=1}^d u_i(1-\prod_{i=1}^d u_i) + (1 - d + \sum_{i=1}^d u_i)(d - \sum_{i=1}^d u_i)$.
Note that from this expression, the maximum verifies $u_1 = \cdots = u_d$, and so the initial $d$-dimensional optimization problem amounts to the one-dimensional problem:
\begin{equation}
    \label{eq:1D_opti_prob}
    \max f_d(u)\coloneqq \frac{\left(u^d - 1 + d -d u\right)^2}{u^d(1-u^d) + (1 - d + d u)(d - d u)},
    \quad \text{under the constraints} \quad
    \begin{cases}
        u \in [0,1]^d, \\
        u \geq \frac{d-1}{d}.
    \end{cases}
\end{equation}
Iterated derivative computations allow to show that $f_d$ is a decreasing function on $[\frac{d-1}{d}, 1]$, so the maximum is reached at $u = \frac{d-1}{d}$, which concludes the proof.

\paragraph{Proof of Corollary~\ref{cor:changement_bound}.}
Let us denote by $R({\mathbb I}\{\theta_1^{(\cdot)} \leq x_1, \theta_2^{(\cdot)} \leq x_2\})$ and $R({\mathbb I}\{\theta_1^{(\cdot)} \leq x_1, \theta_2^{(\cdot)} \geq x_2\})$ the two considered versions of $R(\boldsymbol{x})$ in the bivariate case, with standard uniform margins. We clearly have:
    \begin{align*}
    R_\infty^- 
    &=  \max\limits_{(x_1,x_2)\in [0,1]^2} R\left({\mathbb I}\left\{\theta_1^{(\cdot)} \leq x_1, \theta_2^{(\cdot)} \geq x_2\right\}\right)\\
    &=  \max\limits_{(x_1,x_2)\in [0,1]^2} R\left({\mathbb I}\left\{\theta_1^{(\cdot)} \leq x_1, 1-\theta_2^{(\cdot)} \leq 1-x_2\right\}\right)\\
    &=  \max\limits_{(x_1,x_2)\in [0,1]^2} R\left({\mathbb I}\left\{\theta_1^{(\cdot)} \leq x_1, 1-\theta_2^{(\cdot)} \leq x_2\right\}\right).
\end{align*}
Then, remarking that 
\begin{align*}
    (\theta_1^{(\cdot)}, \theta_2^{(\cdot)}) \sim M_2 &\implies (\theta_1^{(\cdot)}, 1-\theta_2^{(\cdot)}) \sim W_2,\\
            (\theta_1^{(\cdot)}, \theta_2^{(\cdot)}) \sim W_2 &\implies (\theta_1^{(\cdot)}, 1-\theta_2^{(\cdot)}) \sim M_2,\\
    (\theta_1^{(\cdot)}, \theta_2^{(\cdot)}) \sim \Pi_2 &\implies (\theta_1^{(\cdot)}, 1-\theta_2^{(\cdot)}) \sim \Pi_2
\end{align*}
proves the result.

\section{Threshold estimation for $\hat{R}_\infty$}
\label{app:r_inf_threshold}

This section details the computation of the empirical quantiles of $\hat{R}_\infty$ that is done to obtain Table~\ref{tab:rhat_inf_tab}.

\subsection{Univariate case}

\begin{figure}
    \hspace{-0.15\textwidth}
    \setlength\tabcolsep{2pt}
    \begin{tabular}{m{.1\textwidth}m{.33\textwidth}m{.33\textwidth}m{.33\textwidth}}
     & \multicolumn{1}{c}{$m = 2$} &  \multicolumn{1}{c}{$m = 4$} & \multicolumn{1}{c}{$m = 8$}\\
     & \includegraphics[trim=0cm 0cm 0cm 0cm, clip, width = .33\textwidth]{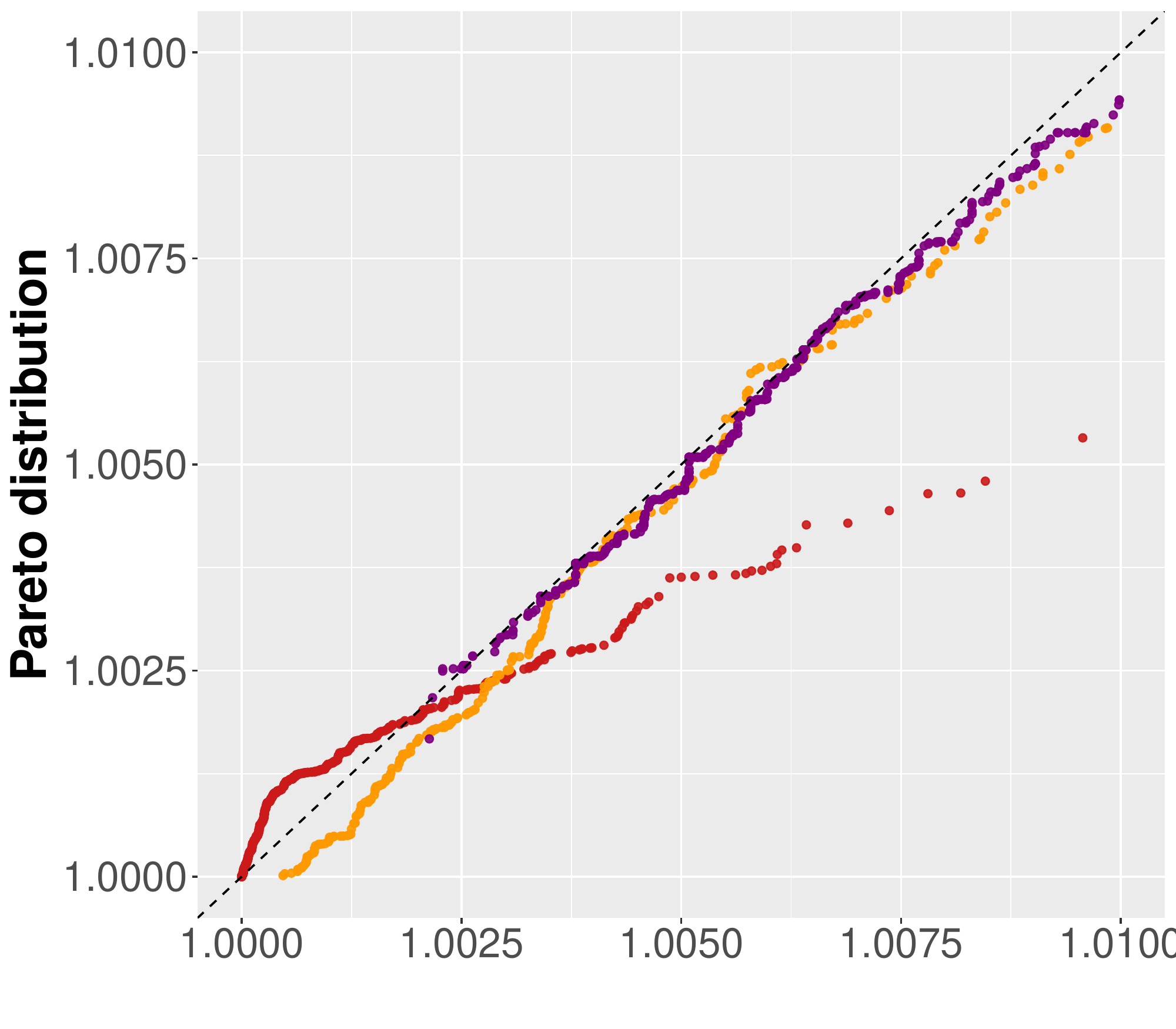} &
    \includegraphics[trim=0cm 0cm 0cm 0cm, clip, width = .33\textwidth]{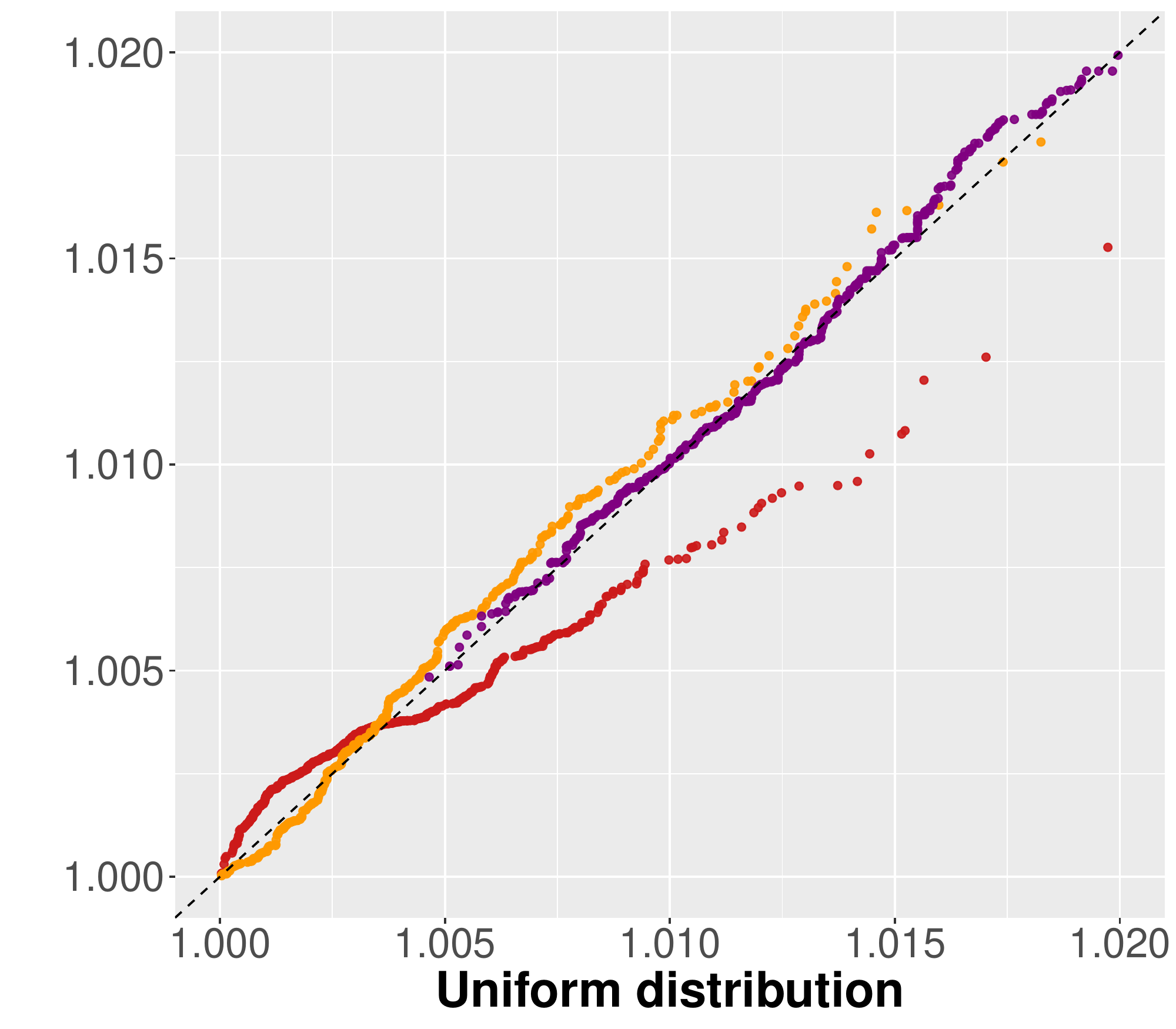} &
    \includegraphics[trim=0cm 0cm 0cm 0cm, clip, width = .33\textwidth]{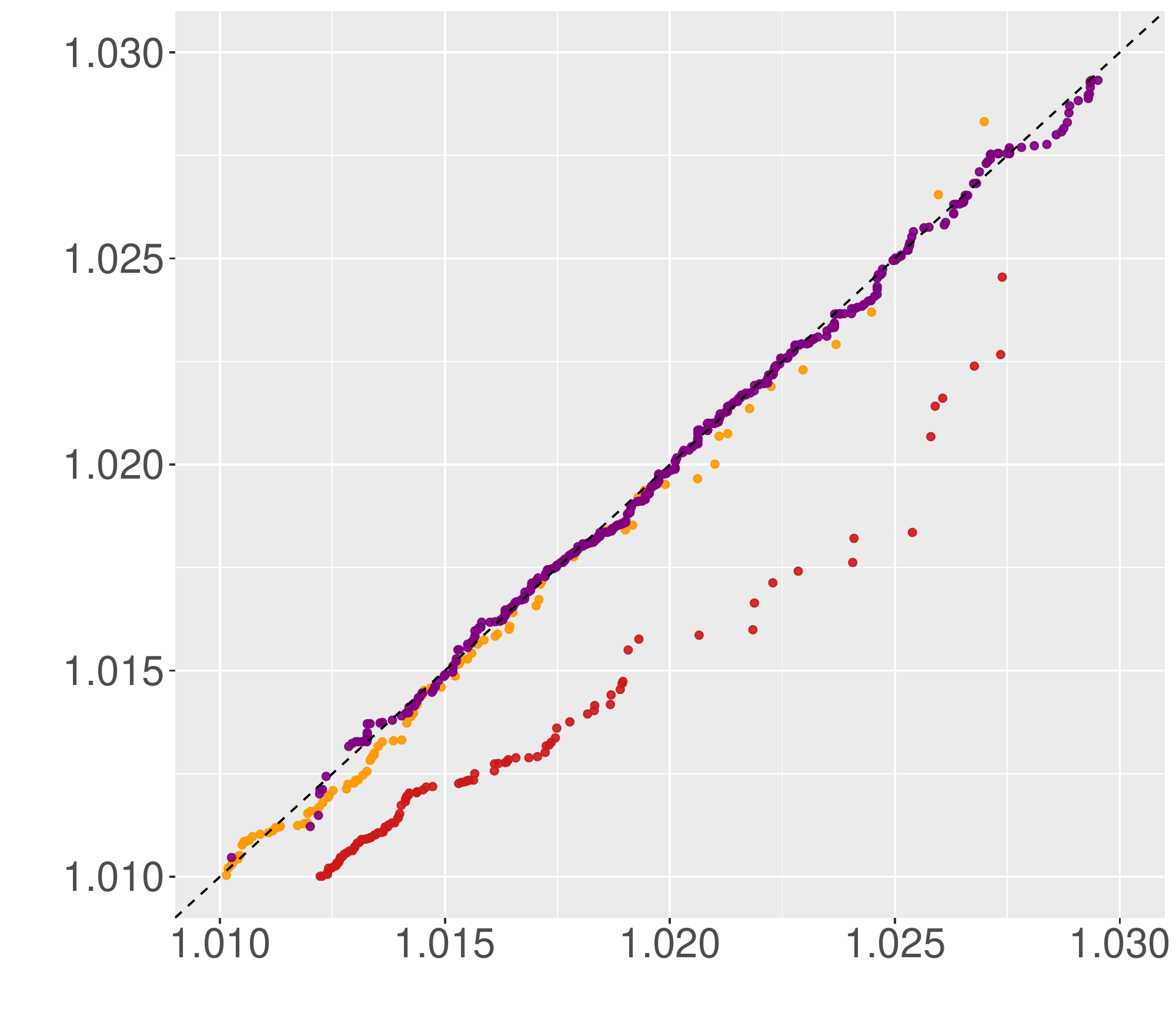} \\ 
    \multicolumn{1}{r}{$\hat{R}$} & \includegraphics[trim=0cm 0cm 0cm 0cm, clip, width = .33\textwidth]{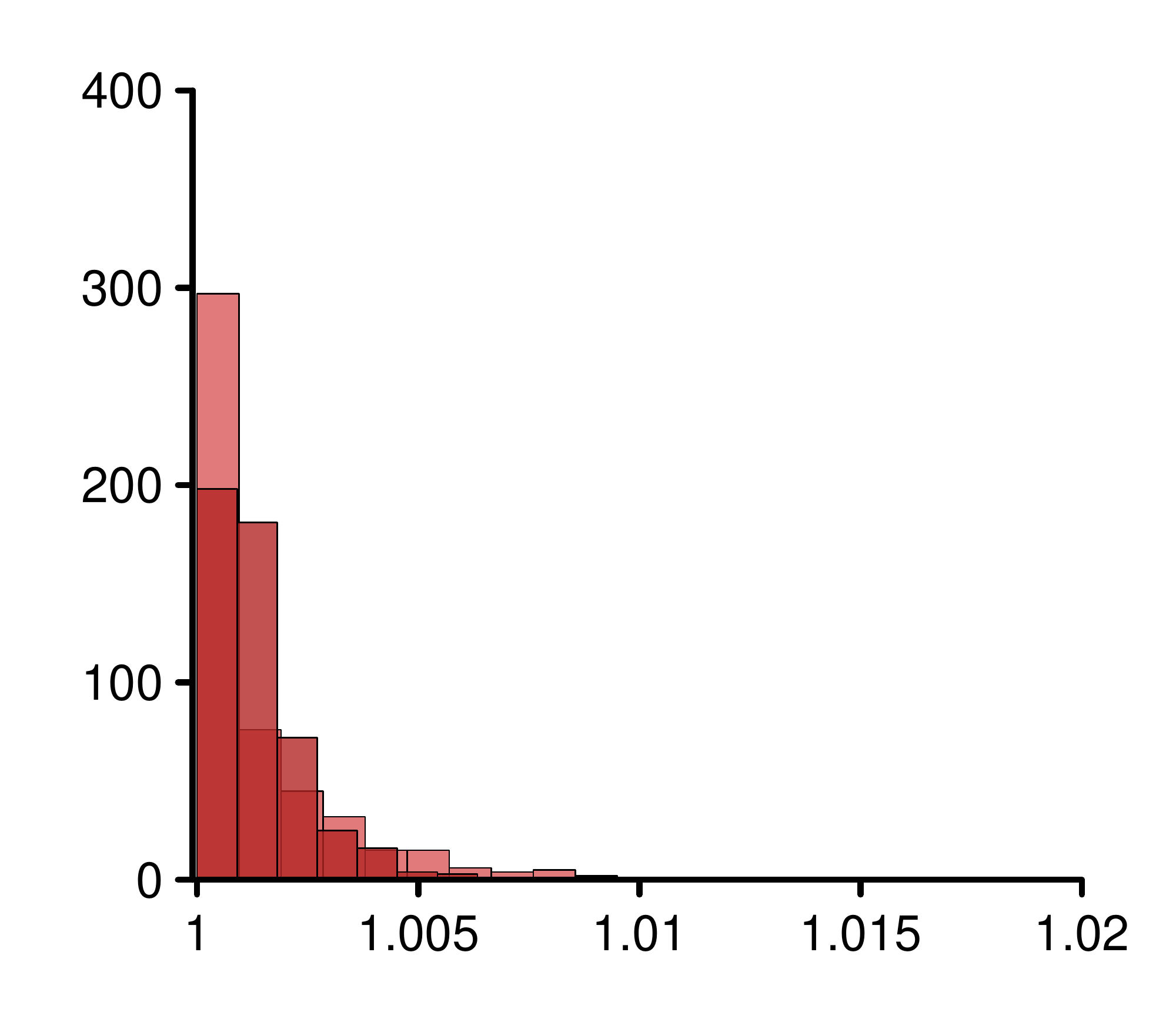} &
    \includegraphics[trim=0cm 0cm 0cm 0cm, clip, width = .33\textwidth]{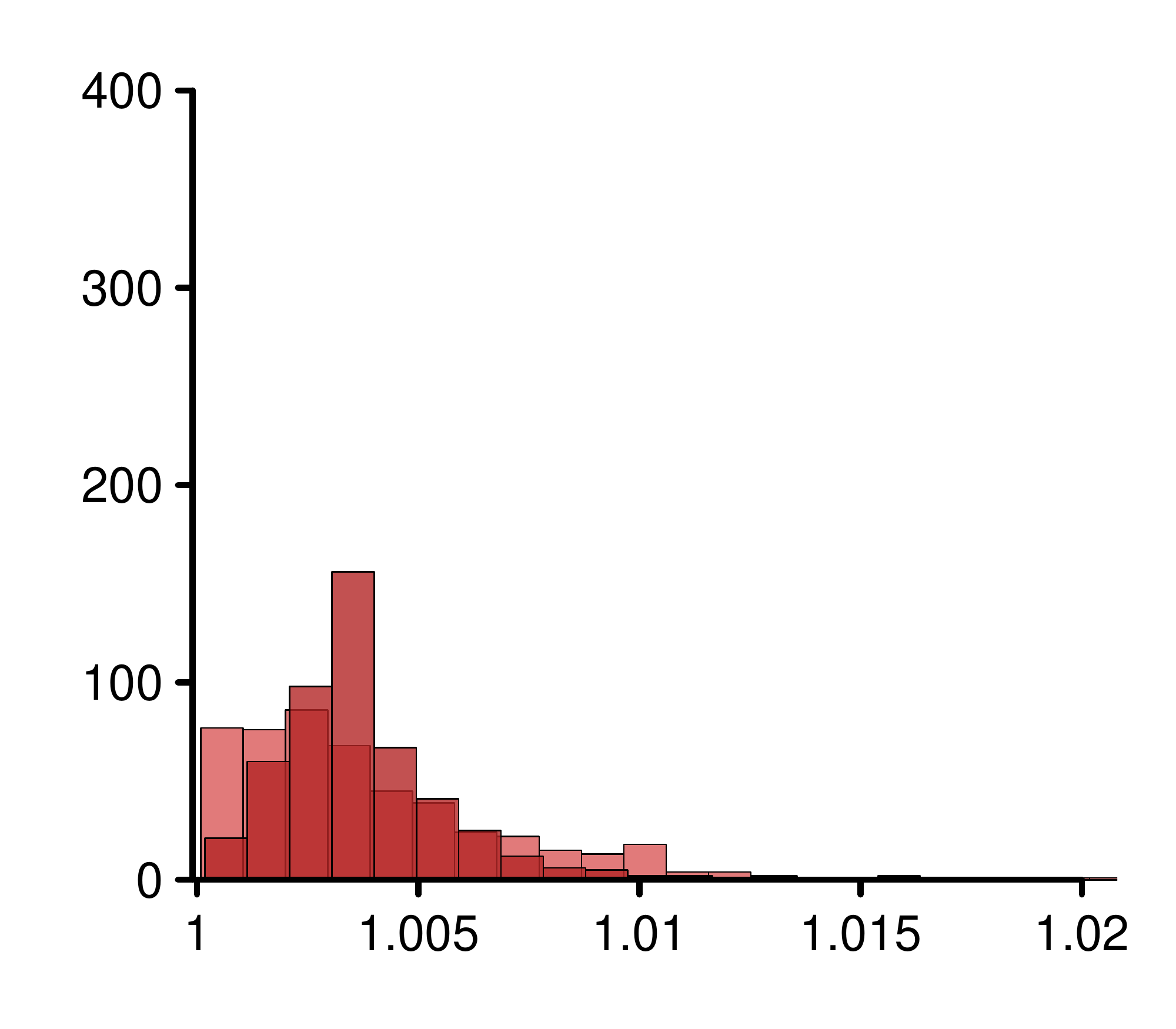} &
    \includegraphics[trim=0cm 0cm 0cm 0cm, clip, width = .33\textwidth]{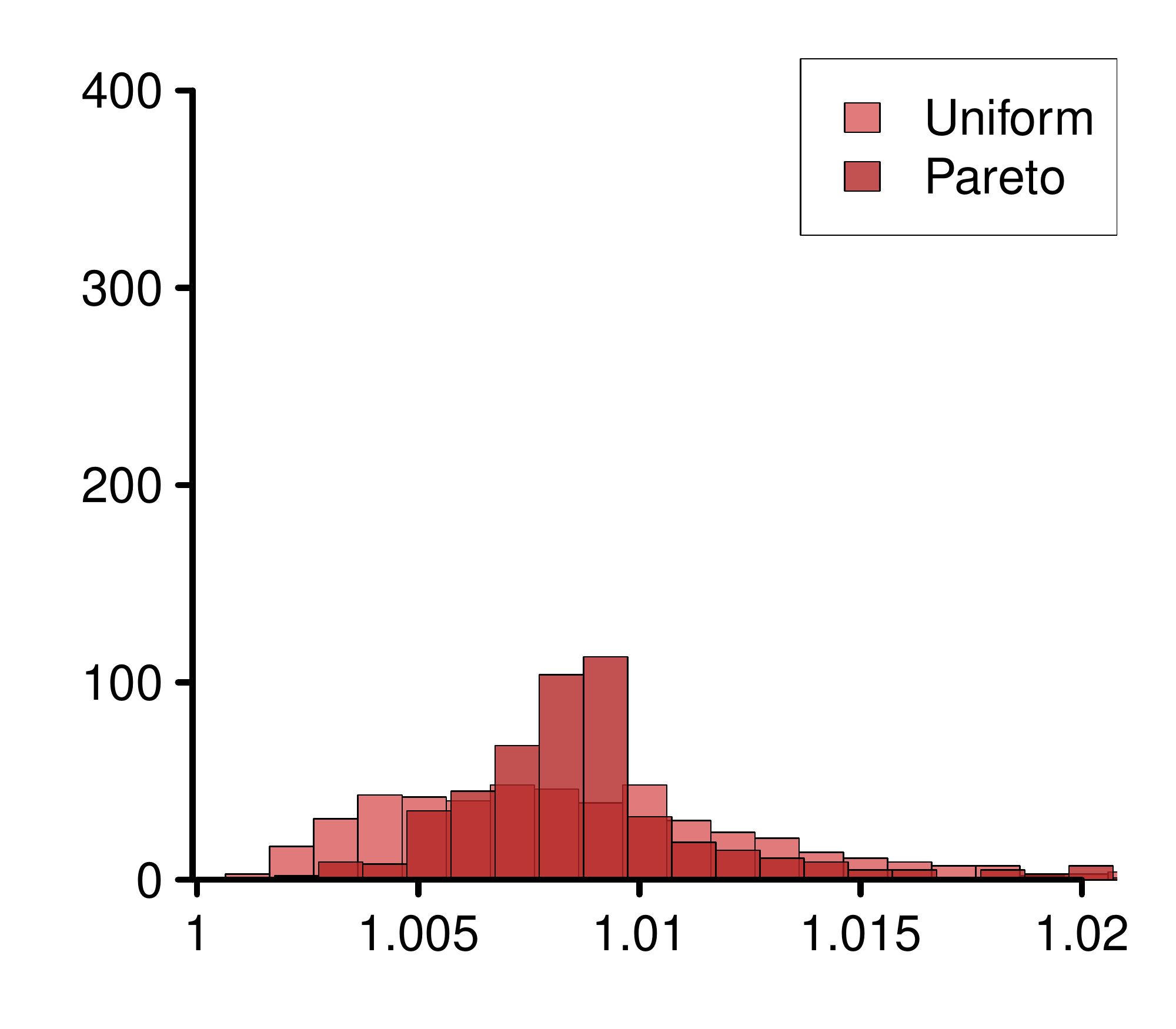} \\
    \multicolumn{1}{r}{rank-$\hat{R}$} & \includegraphics[trim=0cm 0cm 0cm 0cm, clip, width = .33\textwidth]{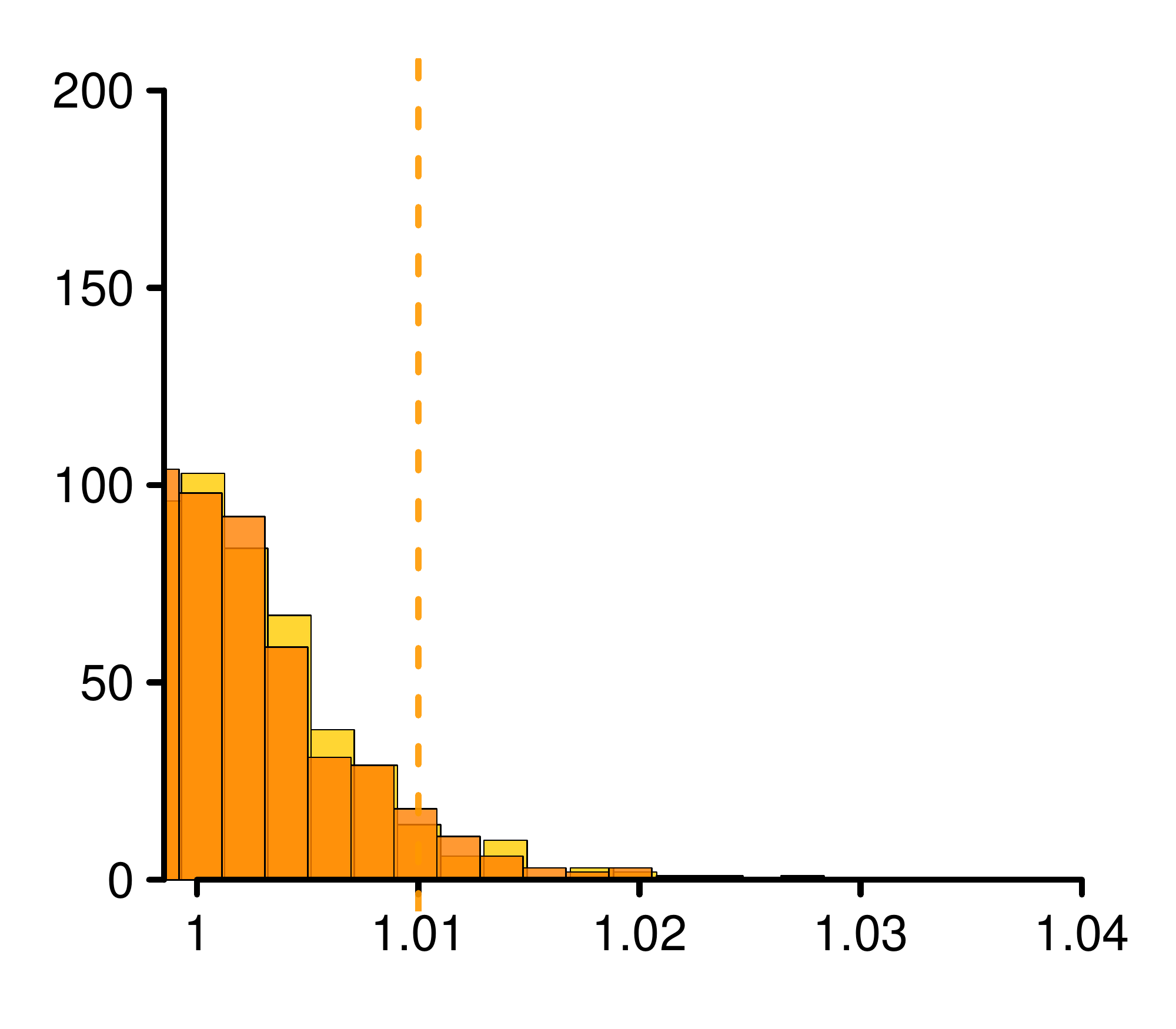} &
    \includegraphics[trim=0cm 0cm 0cm 0cm, clip, width = .33\textwidth]{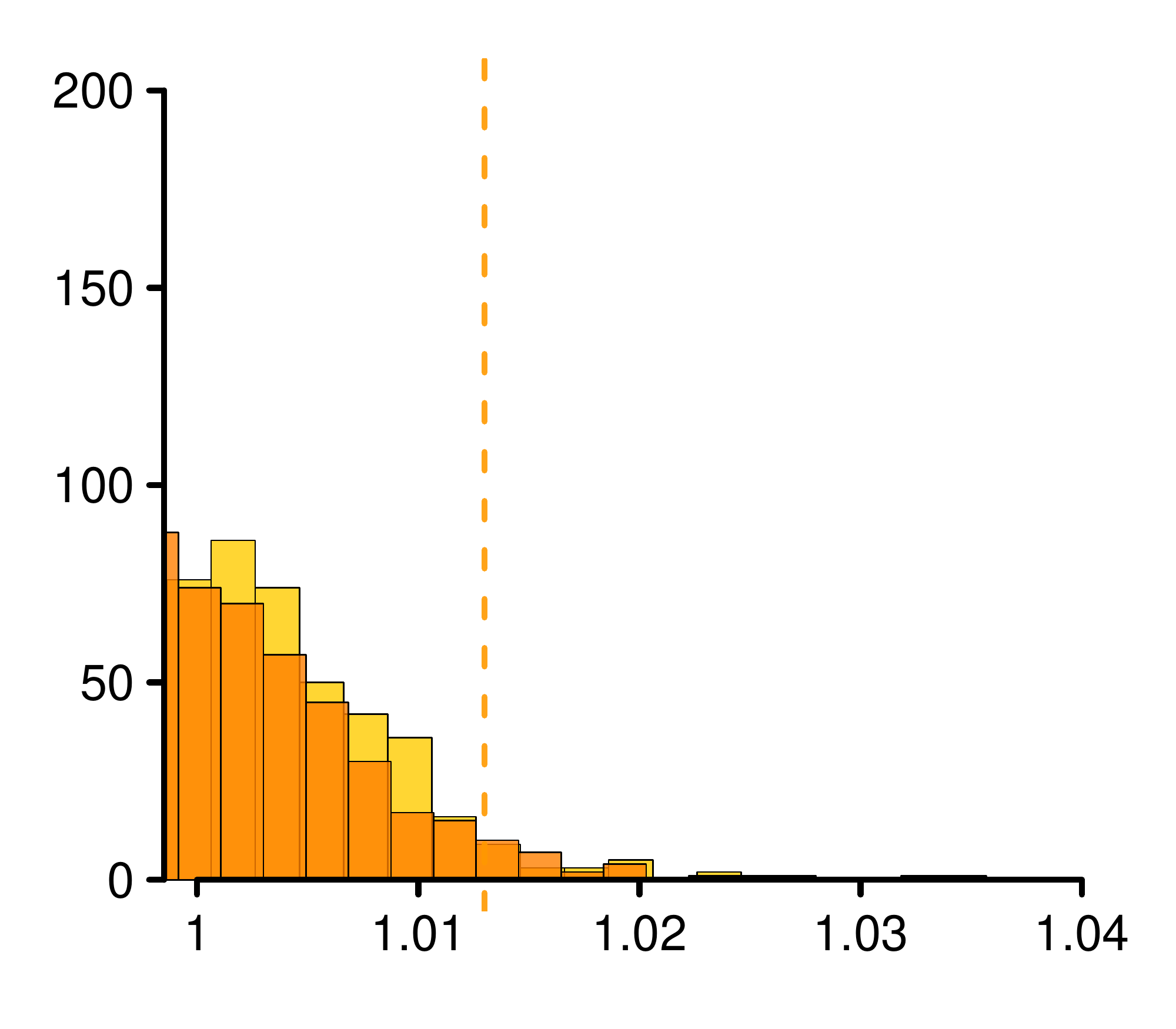} &
    \includegraphics[trim=0cm 0cm 0cm 0cm, clip, width = .33\textwidth]{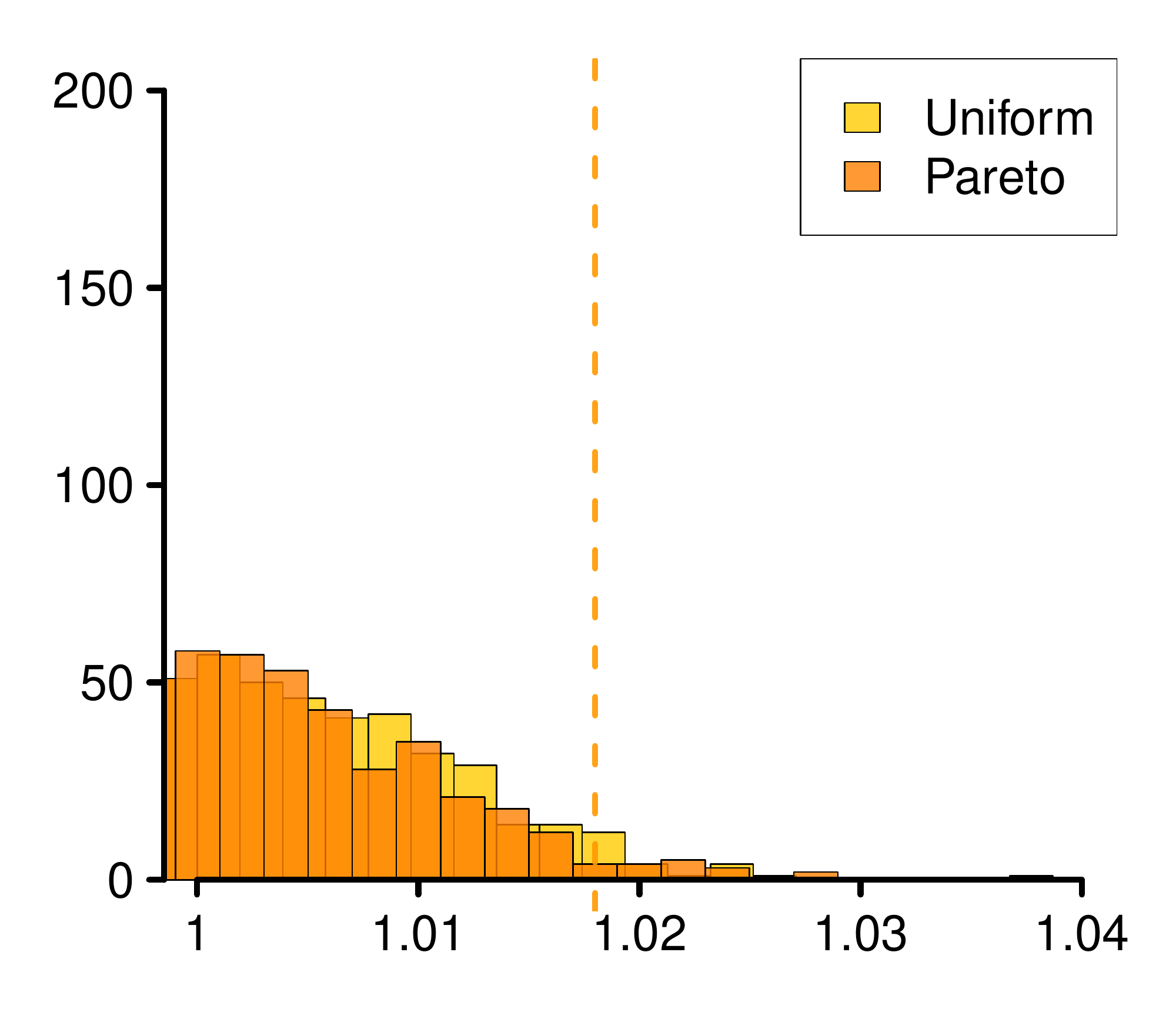} \\   
    \multicolumn{1}{r}{$\hat{R}_\infty$} & \includegraphics[trim=0cm 0cm 0cm 0cm, clip, width = .33\textwidth]{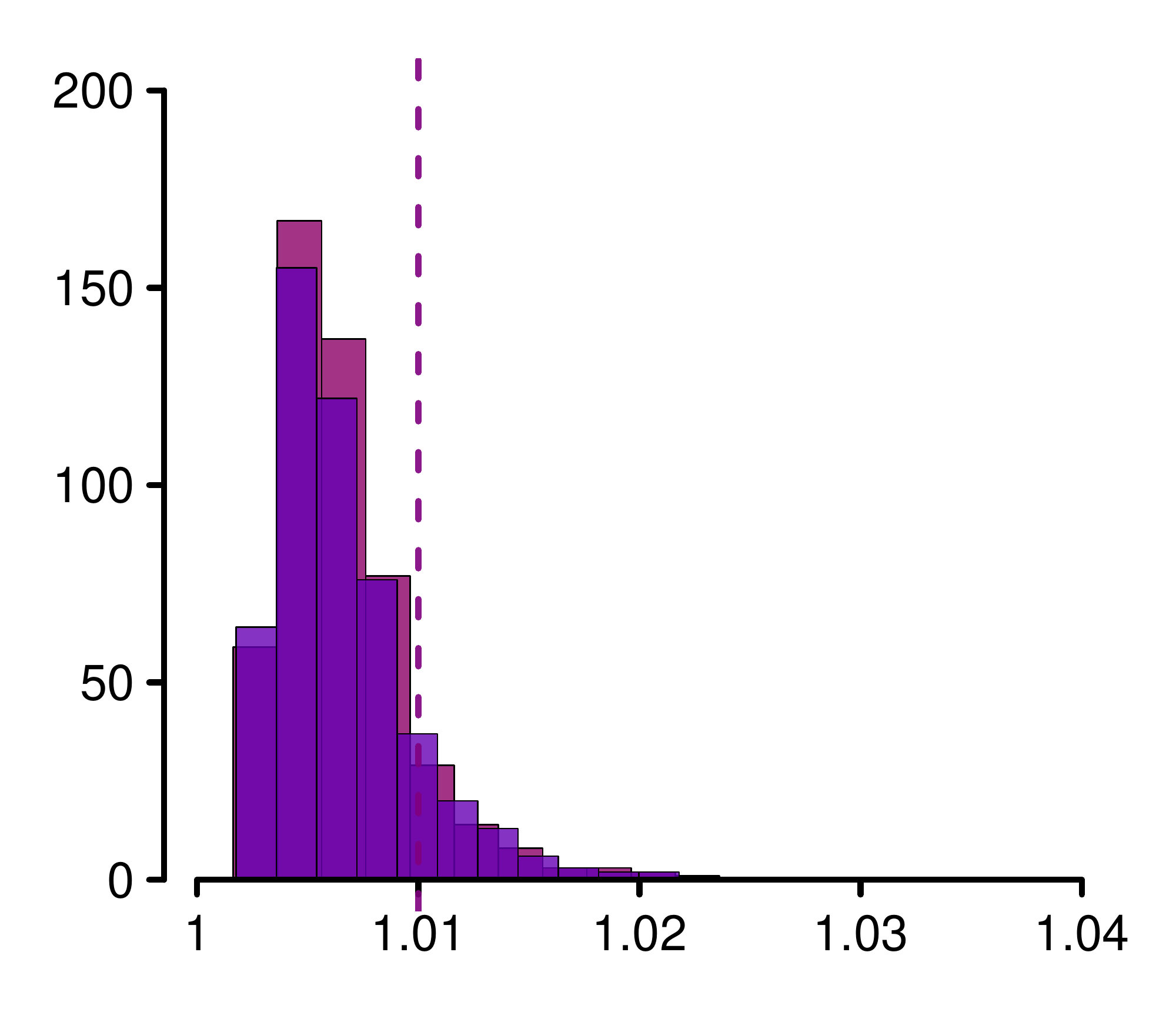} &
    \includegraphics[trim=0cm 0cm 0cm 0cm, clip, width = .33\textwidth]{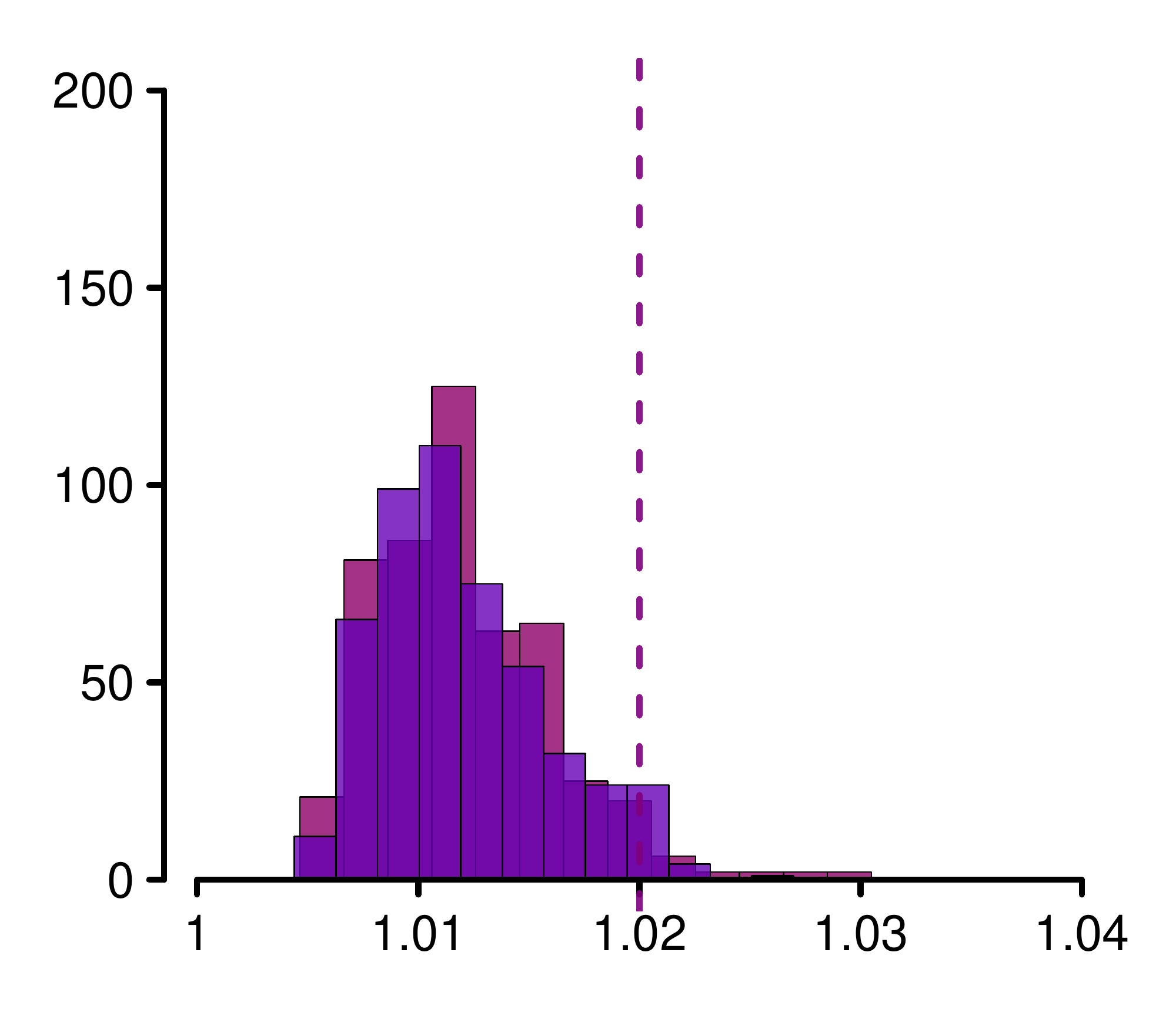} &
    \includegraphics[trim=0cm 0cm 0cm 0cm, clip, width = .33\textwidth]{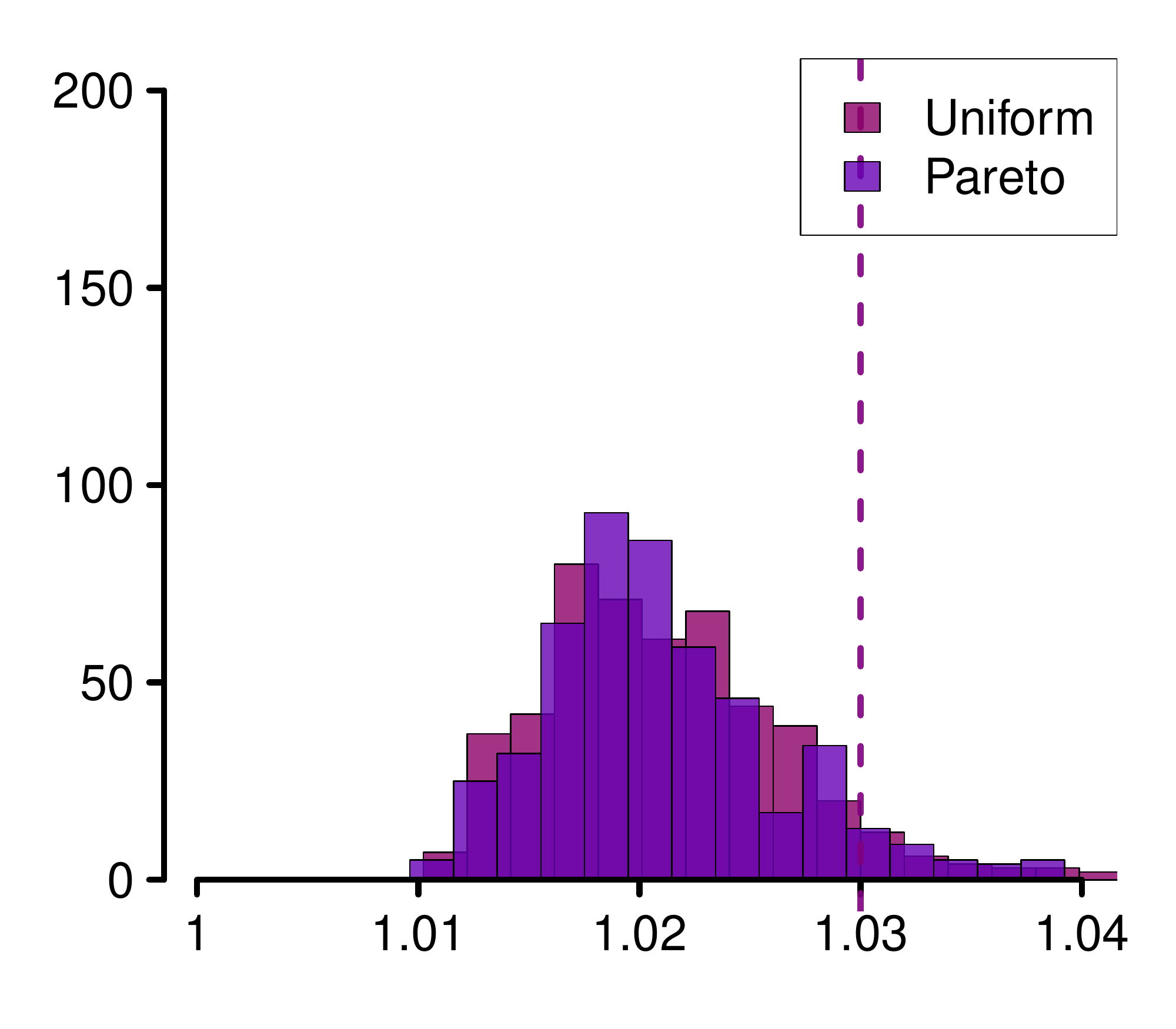} \\
    \end{tabular}
    \caption{Study of  {\color{orange2}$\hat{R}$}, {\color{red2}rank-$\hat{R}$} and {\color{violet2}$\hat{R}_\infty$} when all the distributions are the same, for a number of chains $m \in \{2, 4, 8\}$.
    First row: Q-Q plot that compares $500$ replications computed on $m$ uniform chains $\mathcal{U}(-1, 1)$ with $m$ Pareto$(\alpha = 0.8, \eta = 1)$ ones.
    Second, third and fourth rows: Histograms of $500$ replications with $n \in \{200, 100, 50\}$ for {\color{orange2}$\hat{R}$} (second row), {\color{red2}rank-$\hat{R}$} (third row), and {\color{violet2}$\hat{R}_\infty$}(fourth row) in the case of uniform and Pareto chains. 
    The dashed lines represent the suggested thresholds for {\color{red2}rank-$\hat{R}$} and {\color{violet2}$\hat{R}_\infty$}, corresponding to a confidence level of approximately $95\%$.}
    \label{fig:qqplot}
\end{figure}

\paragraph{Invariance on the underlying distribution under the null hypothesis (all chains have the same distribution).}

A primary step to compute the quantiles of $\hat{R}_\infty$ when all the chain distributions are identical is to verify that such quantiles are well-defined, in the sense that they do not depend on the choice of chain distribution.
This property is expected as using a supremum over the quantiles provides invariance to bijective transformations (see Proposition~\ref{prop:asymptotic_R-hat}(i)).
The first row of Figure~\ref{fig:qqplot} illustrates the behaviour of $\hat{R}$, rank-$\hat{R}_\infty$, and $\hat{R}_\infty$ on two cases:
\begin{itemize}
    \item all chains are uniform; 
    \item all chains are Pareto distributed.
\end{itemize}
The QQ-plots seem to confirm the invariance in distribution of rank-$\hat{R}$ and $\hat{R}_\infty$ for various choices of $m$ (see the yellow and violet dots in the first row of Figure~\ref{fig:qqplot}).
This was also expected for rank-$\hat{R}$ because of the use of a rank-normalization step.
For the traditional $\hat{R}$, the QQ-plots show a difference of distribution (see the red dots in the first row of Figure~\ref{fig:qqplot}), which makes the quantile estimation ill-defined if the chains distribution is not provided. The rest of Figure~\ref{fig:qqplot} shows histograms of replications for $\hat{R}$ (second row), rank-$\hat{R}$ (third row) and $\hat{R}_\infty$ (fourth row), when the chains are distributed according to the same uniform and Pareto distributions.
Histograms for $\hat{R}$ confirm the difference of behaviour when the chain distribution is uniform or Pareto, while the histograms overlap much more for rank-$\hat{R}$ and $\hat{R}_\infty$.

\paragraph{Threshold elicitation.}

As a direct consequence of Proposition~\ref{prop:asymptotic_R-hat}(i), the empirical quantiles of $\hat{R}_\infty$ only depend on the number of chains $m$ and the length $n$. The same holds true for rank-$\hat{R}$.
Therefore, we are able to estimate empirical quantiles with replications of $\hat{R}_\infty$ using any chain distribution (typically uniform): results are reported in Table~\ref{tab:rhat_inf_tab} for a fixed value of $mn = 400$, which is the effective sample size as we are in the i.i.d case.
The same estimation can be done for rank-$\hat{R}$, in order to associate a choice of threshold to a confidence level $1-\alpha$.
Consequently, this viewpoint implies a threshold that depends on $m$ for rank-$\hat{R}$ too, as the empirical quantiles change with $m$ (third row of Figure~\ref{fig:qqplot}). 
Although \cite{Vehtari} do not suggest any tuning with respect to the number of chains, using a threshold of $1.01$ can lead to large type I error:
for example, a threshold of $1.01$ leads to approximately $5\%$ when $m=2$, but to $21\%$ when $m=8$ (third row of Figure~\ref{fig:qqplot}).
Therefore a threshold of 1.01 seems too strong for rank-$\hat{R}$ when $m$ increases: the empirical quantile at $\alpha = 0.05$ suggests to use $1.01$ when $m=2$, $1.013$ when $m=4$ and $1.018$ when $m=8$.
The same observation holds for $\hat{R}_\infty$, and looking at Table~\ref{tab:rhat_inf_tab} and the fourth row of Figure~\ref{fig:qqplot} leads to a threshold of $1.01$ for $m=2$, $1.02$ for $m=4$ and $1.03$ for $m=8$ in order to keep a type I error of approximately $5\%$.

\subsection{Multivariate case}

\gr{In the multivariate extension (see Section~\ref{sec:multivariate}), two thresholds have to be elicited: $R_{\infty,\text{lim}}^{(M)}$ for the convergence of margins, and $R_{\infty,\text{lim}}^{(C)}$ for the convergence of the copula.
Focusing on $R_{\infty,\text{lim}}^{(M)}$, the quantiles of $\hat{R}_\infty$ for the margins are the same as in the univariate case (given in Table~\ref{tab:rhat_inf_tab}), but a Bonferroni correction is necessary to take into account the multiplicity of tests.
For the copula, $\hat{R}^{(\text{max})}_\infty$ is a maximum of multiple versions of $\hat{R}_\infty$ computed on different directions of dependence, so its quantiles are different from $\hat{R}_\infty$ ones. 
However, to reduce the calculation cost, one can compute the quantiles of $\hat{R}_\infty$ under the null hypothesis for the chains, with a Bonferroni correction with $2^{d-1}$ hypotheses.
We estimate the corresponding quantiles using replications to determine the two thresholds $R_{\infty,\text{lim}}^{(M)}$ and $R_{\infty,\text{lim}}^{(C)}$.
Here, the computation is done for several values of $d$, with $d$ relatively small.
Results are reported in Table~\ref{tab:rhat_max_tab}. 
Values in bold confirm the rule of thumb given in Section~\ref{sec:multivariate}: $(R_{\infty,\text{lim}}^{(M)}, R_{\infty,\text{lim}}^{(C)}) = (1.03, 1.03)$ for $m=4$ and $(R_{\infty,\text{lim}}^{(M)}, R_{\infty,\text{lim}}^{(C)}) = (1.04, 1.05)$ for $m=8$.}

\begin{table}
    \gr{
        \begin{center}
        \begin{tabular}{cccccc}
            \toprule
            & & \multicolumn{4}{c}{$R_{\infty,\text{lim}}^{(M)}$}\\
            \cmidrule{3-6}
            & & \multicolumn{4}{c}{$\alpha$}\\
            \cmidrule{3-6}
            $d$ & $m$ & 0.005 & {0.01} & {0.05} & 0.1 \\
            \midrule
            \multirow{4}{*}{2} 
             & 2 & 1.017 & 1.017 & 1.015 & 1.013 \\
             & 3 & 1.024 & 1.023 & 1.019 & 1.018 \\
             & 4 & 1.030 & 1.028 & \textbf{1.025} & 1.022 \\
             & 8 & 1.041 & 1.041 & \textbf{1.037} & 1.034 \\
            \midrule
            \multirow{4}{*}{3} 
             & 2 & 1.022 & 1.021 & 1.018 & 1.016 \\
             & 3 & 1.032 & 1.031 & 1.023 & 1.020 \\
             & 4 & 1.032 & 1.031 & \textbf{1.026} & 1.023 \\
             & 8 & 1.045 & 1.043 & \textbf{1.037} & 1.036 \\
            \midrule
            \multirow{4}{*}{4} 
             & 2 & 1.025 & 1.023 & 1.016 & 1.015 \\
             & 3 & 1.027 & 1.026 & 1.022 & 1.020 \\
             & 4 & 1.027 & 1.026 & \textbf{1.025} & 1.023 \\
             & 8 & 1.047 & 1.044 & \textbf{1.040} & 1.038 \\
            \midrule
            \multirow{4}{*}{5} 
             & 2 & 1.025 & 1.020 & 1.018 & 1.016 \\
             & 3 & 1.027 & 1.027 & 1.025 & 1.022 \\
             & 4 & 1.030 & 1.029 & \textbf{1.026} & 1.024 \\
             & 8 & 1.053 & 1.053 & \textbf{1.038} & 1.036 \\
            \midrule
            \multirow{4}{*}{6} 
             & 2 & 1.022 & 1.021 & 1.018 & 1.015 \\
             & 3 & 1.024 & 1.023 & 1.021 & 1.021 \\
             & 4 & 1.030 & 1.030 & \textbf{1.030} & 1.024 \\
             & 8 & 1.047 & 1.046 & \textbf{1.039} & 1.037 \\
            \bottomrule
        \end{tabular} 
        \hspace{-.3cm}
        \begin{tabular}{cccccc}
            \toprule
            & & \multicolumn{4}{c}{$R_{\infty,\text{lim}}^{(C)}$}\\
            \cmidrule{3-6}
            & & \multicolumn{4}{c}{$\alpha$}\\
            \cmidrule{3-6}
             & $m$ & 0.005 & {0.01} & {0.05} & 0.1 \\
            \midrule
            \multirow{4}{*}{} 
             & 2 & 1.026 & 1.026 & 1.019 & 1.016 \\
             & 3 & 1.029 & 1.028 & 1.024 & 1.021 \\
             & 4 & 1.033 & 1.030 & \textbf{1.026} & 1.024 \\
             & 8 & 1.052 & 1.050 & \textbf{1.040} & 1.038 \\
            \midrule
            \multirow{4}{*}{} 
             & 2 & 1.022 & 1.020 & 1.019 & 1.018 \\
             & 3 & 1.031 & 1.028 & 1.025 & 1.023 \\
             & 4 & 1.038 & 1.036 & \textbf{1.030} & 1.027 \\
             & 8 & 1.052 & 1.049 & \textbf{1.047} & 1.043 \\
            \midrule
            \multirow{4}{*}{} 
             & 2 & 1.029 & 1.026 & 1.022 & 1.021 \\
             & 3 & 1.028 & 1.028 & 1.026 & 1.024 \\
             & 4 & 1.035 & 1.034 & \textbf{1.033} & 1.030 \\
             & 8 & 1.051 & 1.050 & \textbf{1.048} & 1.047 \\
            \midrule
            \multirow{4}{*}{}
             & 2 & 1.024 & 1.024 & 1.021 & 1.021 \\
             & 3 & 1.028 & 1.028 & 1.026 & 1.024 \\
             & 4 & 1.043 & 1.043 & \textbf{1.040} & 1.036 \\
             & 8 & 1.049 & 1.049 & \textbf{1.048} & 1.048 \\
            \midrule
            \multirow{4}{*}{} 
             & 2 & 1.020 & 1.020 & 1.019 & 1.018 \\
             & 3 & 1.028 & 1.028 & 1.025 & 1.024 \\
             & 4 & 1.035 & 1.035 & \textbf{1.034} & 1.033 \\
             & 8 & 1.059 & 1.059 & \textbf{1.058} & 1.055 \\
            \bottomrule
        \end{tabular}
        \end{center}
    }
    \caption{Left: Empirical quantiles of $\hat{R}_\infty$ for the margins with a Bonferroni correction, under the null hypothesis that all chains have the same distribution, for a fixed value of $mn = 400$. 
    Right: Empirical quantiles of $\hat{R}^{(\text{max})}_\infty$ for the copula, for a fixed value of $mn = 400$. 
    We have used 500 replications for estimation. Values in bold justify the rule of thumb proposed in Section~\ref{subsec:multivariate_def} for $m\in\{4,8\}.$}
    \label{tab:rhat_max_tab}
\end{table}

\section{Examples of closed-form $R(x)$ and $R_\infty$}
\label{app:example_theorique}

We start by a lemma providing useful tools to simplify the calculation of $R(x)$ and $R_\infty$
in the univariate case and when all chains but one have converged.
We then review families of distributions for which $R(x)$ and $R_\infty$ can be computed in closed-form.

\begin{LemApp}
\label{lem-support}
Assume the assumptions of Proposition~\ref{prop-calcul-R} hold with $F\coloneqq F_1=\dots=F_{m-1}\neq F_m$.
\begin{itemize}
\item[\normalfont (i)] Then~(\ref{eq:R_theorique}) can be simplified as
\begin{equation}
    R(x) = \sqrt{1 + \frac{(m-1)(F(x)-F_m(x))^2}{m\left((m-1)F(x)(1-F(x))+F_m(x)(1-F_m(x))\right)}}.
    \label{eq:R_single_degenerated_chain}
\end{equation}
\item[\normalfont (ii)] If $F$ and $F_m$ are symmetrical distributions wrt 0, then $R$ is an even function.
\item [\normalfont (iii)] Let $a$ and $a_m \in \mathbb{R} \cup \{-\infty\}$ be the starting points of $F$ and $F_m$ respectively,
and assume $a\leq a_m$. Then, $R(\cdot)$ reaches its supremum on $[a_m,\infty)$:
$$
R_\infty \geq \sqrt{1+ \frac{F(a_m)}{m\left(1-F(a_m)\right)}}.
$$
\item [\normalfont (iv)] Let $b$ and $b_m \in \mathbb{R} \cup \{+\infty\}$ be the endpoints of $F$ and $F_m$ respectively,
with $b\leq b_m$. Then, $R(\cdot)$ reaches its supremum on $(-\infty,b]$ and
$$
R_\infty \geq \sqrt{1+ \frac{(m-1)(1-F_m(b))}{m F_m(b)}}.
$$
\end{itemize}
\end{LemApp}

\paragraph{Proof.} (i) and (ii) are straightforward. For (iii) remark that,
when $x\leq a$, $F(x)=F_m(x)=0$ so that $R(x)=1$.
Besides, for all $x\in[a,a_m]$, one has $F_m(x)=0$ and thus 
$$
R^2(x)=1+ \frac{F(x)}{m(1-F(x))}.
$$
The above defined function is increasing so that the supremum of $R^2(\cdot)$ is reached for $x\geq a_m$ and therefore $R_\infty\geq R(a_m)$. 
Similarly for (iv), when $x\geq b_m$, $F(x)=F_m(x)=1$ so that $R(x)=1$.
Besides,
for all $x\in[b,b_m]$, one has $F(x)=0$ and thus 
$$
R^2(x) = 1 + \frac{(m-1)(1-F_m(x))}{m F_m(x)}.
$$
The above defined function is decreasing so that the supremum of $R(\cdot)$ is reached for $x\leq b$ and therefore $R_\infty\geq R(b)$. 
    
\begin{LemApp}[Uniform distribution] 
\label{lem-unif}
Assume that $F_1=\dots=F_{m-1}$ are the cdf of the uniform distribution ${\mathcal U}(-\sigma,\sigma)$
while $F_m$ is  the cdf of the uniform distribution ${\mathcal U}(-\sigma_m,\sigma_m)$ with $0<\sigma\leq\sigma_m$.
Then, 
$$
R^2(x) = 
    \begin{cases}
        1 + \frac{\big(\frac{1}{\sigma} - \frac{1}{\sigma_m}\big)^2}{\frac{m^2}{(m-1)x^2} - m\big(\frac{1}{\sigma^2} + \frac{1}{(m-1)\sigma_m^2}\big)}
        &\quad \text{if $|x| \leq \sigma$}\,,\\[5pt]
        1 + \frac{m-1}{m}\left(1 - \frac{2}{1+\sigma_m/|x|}\right)
        &\quad \text{if $\sigma \leq |x| \leq \sigma_m$}\,,\\[5pt]
        1  &\quad \text{if $|x| \geq \sigma_m$}.
    \end{cases}
$$
Moreover,
$$
R_\infty =   
        R(\pm\sigma) = 
        \sqrt{1 + \frac{m-1}{m}\left(1-\frac{2}{1+\frac{\sigma_m}{\sigma}}\right)}.
$$
\end{LemApp}
\paragraph{Proof.}
   Recall that
    \begin{align*}
        F_1(x) &= \cdots = F_{m-1}(x) =  \frac{x}{2\sigma} + \frac{1}{2}, \quad \forall x \in \left[-\sigma; \sigma\right],\\
        \text{and}\quad F_m(x) &=  \frac{x}{2\sigma_m} + \frac{1}{2}, \quad \forall x \in \left[-\sigma_m; \sigma_m\right].
    \end{align*}
    The case $|x| \geq \sigma_m$ is clear, we investigate the two other ones.
    First, if $\sigma \leq x \leq \sigma_m$, then $F_j(x) = 1$ for $j = 1 ,\ldots, m-1$ and Lemma~\ref{lem-support}(i) yields
    \begin{align*}
        R^2(x) = 1 + \frac{(m-1)\left(\frac{1}{2}-\frac{x}{2\sigma_m}\right)^2}{m\left(\frac{1}{2}+\frac{x}{2\sigma_m}\right)\left(\frac{1}{2}-\frac{x}{2\sigma_m}\right)}
        = 1 + \frac{m-1}{m} \left(1 - \frac{2}{1+\frac{\sigma_m}{x}}\right).
    \end{align*}
    Using, Lemma~\ref{lem-support}(ii), allows concluding for $\sigma \leq |x| \leq \sigma_m$.
    Finally, if $|x| \leq \sigma$, then:
    \begin{align*}
        R^2(x) &= 1 + \frac{x^2\big(\frac{1}{2\sigma} - \frac{1}{2\sigma_m}\big)^2}{m\left(\frac{1}{2}+\frac{x}{2\sigma}\right)\left(\frac{1}{2}-\frac{x}{2\sigma}\right) + \frac{m}{m-1}\left(\frac{1}{2}+\frac{x}{2\sigma_m}\right)\left(\frac{1}{2}-\frac{x}{2\sigma_m}\right)}\\
        &= 1 + \frac{\big(\frac{1}{\sigma} - \frac{1}{\sigma_m}\big)^2}{\frac{m^2}{(m-1)x^2} - m\big(\frac{1}{\sigma^2} + \frac{1}{(m-1)\sigma_m^2}\big)}.
    \end{align*}
    Lemma~\ref{lem-support}(iv) entails that the maximum is reached for $x\in[-\sigma; \sigma]$. The above expression shows that the maximum is located at $x = \pm \sigma$, which gives the result.

    \begin{LemApp}[Pareto distribution]
    \label{lem-pareto}
    Assume that $F_1=\cdots=F_{m-1}$ are the cdf of the Pareto$(\alpha, \eta)$ distribution with $\alpha > 0$ the shape parameter and $\eta > 0$ the position parameter.
    Let $F_m$ be the cdf of the Pareto$(\alpha, \eta_m)$ distribution  with $0<\eta \leq \eta_m$.
    Then,
    $$
    R^2(x) = 
        \begin{cases}
            1 + \frac{1}{m} \left(\left(\frac{x}{\eta}\right)^{\alpha} - 1\right)
            &\quad \text{if $\eta \leq x \leq \eta_m$} \,,\\[5pt]
            1 + \frac{1}{m} \frac{(\eta^{\alpha}-\eta_m^{\alpha})^2}{\left(\eta^{\alpha} + \frac{\eta_m^{\alpha}}{m-1}\right)x^{\alpha} - \left(\eta^{2\alpha} + \frac{\eta_m^{2\alpha}}{m-1}\right)}
            &\quad \text{if $\eta_m \leq x$}\,,\\[5pt]
            1 &\quad \text{if $x \leq \eta$}.
        \end{cases}
    $$
    \begin{equation*}
    \text{Moreover,}\quad    R_\infty = R(\eta_m) = \sqrt{1+  \frac{1}{m}\left(\left(\frac{\eta_m}{\eta}\right)^{\alpha}-1\right)}.
    \end{equation*}
    \end{LemApp}
    
    \paragraph{Proof.} Recall that $F_1(x) = \cdots = F_{m-1}(x) = 1 - \left({x}/{\eta}\right)^{-\alpha},
        \quad \forall x \in [\eta, +\infty)$ and $F_m(x) =  1 - \left({x}/{\eta_m}\right)^{-\alpha},
        \quad \forall x \in [\eta_m, +\infty)$.
    In the case where $\eta\leq x \leq \eta_m$, $F_m(x) = 0$, and using Lemma~\ref{lem-support}(iii) entails that $R(\cdot)$ is increasing on $[\eta,\eta_m]$ and
    \begin{equation*}
        R^2(x) = 1 + \frac{1}{m} \frac{F(x)}{1-F(x)}
        = 1 + \frac{1}{m} \left(\left(\frac{x}{\eta}\right)^{\alpha} - 1\right).
    \end{equation*}
    Moreover, for $\eta_m \leq x$, replacing the Pareto cdf in~(\ref{eq:R_single_degenerated_chain}) yields
    \begin{align*}
        R^2(x) &= 1 + \frac{1}{m}
                \frac{x^{-2\alpha}\left(\eta^{\alpha} - \eta_m^{\alpha}\right)^2}{x^{-\alpha}\eta^{\alpha}(1-x^{-\alpha}\eta^{\alpha}) + \frac{1}{m-1} x^{-\alpha}\eta_m^{\alpha}(1-x^{-\alpha}\eta_m^{\alpha})}\\
                &= 1 + \frac{1}{m} \frac{(\eta^{\alpha}-\eta_m^{\alpha})^2}{\left(\eta^{\alpha} + \frac{\eta_m^{\alpha}}{m-1}\right)x^{\alpha} - \left(\eta^{2\alpha} + \frac{\eta_m^{2\alpha}}{m-1}\right)}.
    \end{align*}
    Clearly, $R^2(\cdot)$ is decreasing on $[\eta_m, +\infty)$ and is extended by continuity at $x=\eta_m$. In conclusion, $R^2(\cdot)$ is maximum at $x=\eta_m$, and the result is proved.

    \begin{LemApp}[Uniform vs Laplace distribution]
    \label{lem:laplace_uniform_proof}
    Assume that $m=2$, $F_1$ is the cdf of the uniform distribution $\mathcal{U}(-\sigma, \sigma)$ and $F_2$ is the cdf of the centred Laplace distribution $\mathcal{L}(0, \sigma/2)$ with $\sigma > 0$. 
    Then for any $x$, $R(x) = R_{1}(x/\sigma)$ with 
    \begin{equation*}
        R_{1}^2(x) = 
        \begin{cases}
            1 + \frac{\exp ({-|x|})}{2(2-\exp({-|x|}))}
            &\quad \text{if $|x| \geq 2$} \,,\\[5pt]
            1 + \frac{1}{2}\frac{(|x|/2 - 1 + \exp(-|x|))^2}{1-x^2/4+2\exp(-|x|)(2-\exp(-|x|))}
            &\quad \text{if $ |x| \leq 2$}.
        \end{cases}
    \end{equation*}
    \begin{equation*}
    \text{Moreover,}\quad    R_\infty = R(\pm\sigma) = R_{1}(\pm 1)=
        \sqrt{1+\frac{1}{2(2e^2-1)}}.
    \end{equation*}
    \end{LemApp}
    \paragraph{Proof.} 
    In view of Lemma~\ref{lem-support}(ii,iii), it is sufficient to compare the values of $R_{1}(x)$ on $(-2, 0]$ with $R_{1}(-2)$. Then, the derivation of $R_{1}(x)$ is similar to the ones done in Lemma~\ref{lem-unif} and Lemma~\ref{lem-pareto}.
    Routine calculations show that $R_{1}$ has indeed a local maximum on $(-2,0)$, but it remains lower than  $R_1(-2)$ (see last column of Figure~\ref{fig:example_known_dist} for an illustration), which is therefore the value of $R_\infty$.

\end{appendices}

\bibliography{references}

\end{document}